\documentclass[11pt]{amsart}
\usepackage[foot]{amsaddr}
\usepackage{graphicx}
\graphicspath{images/}

\usepackage{amsmath}
\usepackage{amsfonts}
\usepackage{tikz-cd} 
\usepackage{amssymb}
\usepackage{amsthm}
\usepackage{tikzsymbols}
\usepackage{mathtools}
\usepackage{xcolor}
\usepackage{adjustbox}
\usepackage{blkarray}
\usepackage{makecell}
\usepackage{float}
\usepackage{caption}
\usepackage[font={small,it}]{caption}
\usepackage[margin=1cm]{caption}
\usepackage{longtable}

\usepackage{theoremref}
\usepackage{hyperref}

\theoremstyle{plain}
\newtheorem{theorem}{Theorem}[subsection]
\newtheorem{lemma}[theorem]{Lemma}
\newtheorem{corollary}{Corollary}[theorem]

\theoremstyle{definition}

\newtheorem{example}[theorem]{Example}
\newtheorem{observation}[theorem]{Observation}
\newtheorem{strategy}[theorem]{Strategy}

\usepackage{blindtext}

\usepackage[many]{tcolorbox}
\tcolorboxenvironment{lemma}{
  breakable,
  colback=black!10!white,
  boxrule=0pt,
  boxsep=1pt,
  left=2pt,right=2pt,top=2pt,bottom=2pt,
  oversize=2pt,
  sharp corners,
  before skip=\topsep,
  after skip=\topsep,
}

\tcolorboxenvironment{observation}{
  breakable,
  colback=black!10!white,
  boxrule=0pt,
  boxsep=1pt,
  left=2pt,right=2pt,top=2pt,bottom=2pt,
  oversize=2pt,
  sharp corners,
  before skip=\topsep,
  after skip=\topsep,
}

\tcolorboxenvironment{example}{
  breakable,
  colback=black!10!white,
  boxrule=0pt,
  boxsep=1pt,
  left=2pt,right=2pt,top=2pt,bottom=2pt,
  oversize=2pt,
  sharp corners,
  before skip=\topsep,
  after skip=\topsep,
}

\tcolorboxenvironment{strategy}{
  breakable,
  colback=black!10!white,
  boxrule=0pt,
  boxsep=1pt,
  left=2pt,right=2pt,top=2pt,bottom=2pt,
  oversize=2pt,
  sharp corners,
  before skip=\topsep,
  after skip=\topsep,
}

\tcolorboxenvironment{corollary}{
  breakable,
  colback=black!10!white,
  boxrule=0pt,
  boxsep=1pt,
  left=2pt,right=2pt,top=2pt,bottom=2pt,
  oversize=2pt,
  sharp corners,
  before skip=\topsep,
  after skip=\topsep,
}

\tcolorboxenvironment{theorem}{
  breakable,
  colback=blue!10!white,
  boxrule=0pt,
  boxsep=1pt,
  left=2pt,right=2pt,top=2pt,bottom=2pt,
  oversize=2pt,
  sharp corners,
  before skip=\topsep,
  after skip=\topsep,
}

\usepackage{subfiles} 

\title[$k$-Polymatroids with Binary $k$-Natural Matroids]{The Excluded Minors for $k$-Polymatroids with Binary $k$-Natural Matroids}
\author[F. Young]{Fiona Young}
\address{Department of Mathematics, Cornell University, Ithaca, New York, 14853}
\email{fy79@cornell.edu}
\date{March 17, 2023}

\begin{document}

\maketitle

ABSTRACT. If $\mathcal{C}$ is a minor-closed class of matroids, then the class $\widetilde{\mathcal{C}}'_k$ of $k$-polymatroids whose $k$-natural matroids are in $\mathcal{C}$ is also minor-closed. We investigate the following question: When $\mathcal{C}$ is the class of binary matroids, what are the excluded minors for $\widetilde{\mathcal{C}}'_k$? When $k = 1$, $\widetilde{\mathcal{C}}'_1$ is simply the class of binary matroids, which has $U_{2,4}$ as its only excluded minor. Joseph E. Bonin and Kevin Long answered the question for $k = 2$ and found that the set of excluded minors for $\widetilde{\mathcal{C}}'_2$ is infinite. We determine the sets of excluded minors for $\widetilde{\mathcal{C}}'_k$ when $k \geq 3$ and find that they are finite. There are $12$ excluded minors for $\widetilde{\mathcal{C}}'_3$ and when $k > 3$, there are $k+7$ excluded minors for $\widetilde{\mathcal{C}}'_k$.

\tableofcontents

\section{Introduction}
\label{section:introduction}

An \textit{integer polymatroid} is a pair $(E,\rho)$ where $E$ is a finite set and $\rho: 2^E \to \mathbb{Z}$ is a function satisfying:
\begin{enumerate}
    \item $\rho(\emptyset) = 0$.
    \item \textit{Monotonicity}: If $A \subseteq B \subseteq E$, then $\rho(A) \leq \rho(B)$.
    \item \textit{Submodularity}: For all $A, B \subseteq E$, $\rho(A) + \rho(B) \geq \rho(A \cap B) + \rho(A \cup B)$.
\end{enumerate}

In this paper, we refer to integer polymatroids as simply \textit{polymatroids}. If $\rho(e) \leq k$ for all $e \in E$, then $\rho$ is a \textit{k-polymatroid}. In particular, all $1$-polymatroids are matroids. Polymatroids generalize matroids by allowing elements of higher rank, and this has a nice geometric interpretation: while matroids consist only of loops (elements of rank $0$) and points (elements of rank $1$), $k$-polymatroids can contain lines (elements of rank $2$), planes (elements of rank $3$), etc.\ up to elements of rank $k$.

The minors of a matroid are crucial to illuminating its structural decomposition. Many important classes of matroids are \textit{minor-closed}, that is, they are closed under the operations of deletion and contraction. This allows us to characterize such a class $\mathcal{C}$ by its set of \textit{excluded minors}: matroids that are not in $\mathcal{C}$, whose proper minors are all in $\mathcal{C}$. One such minor-closed class is the class of $\mathbb{F}$-representable matroids, for a fixed field $\mathbb{F}$. Rota's conjecture states that if $\mathbb{F}$ is finite, then the set of excluded minors for $\mathbb{F}$-representable matroids is finite. Geelen, Gerards, and Whittle \cite{geelen} announced a proof of this conjecture in 2014.

It is natural to ask the same question for polymatroids. Characterizing the class of $\mathbb{F}$-representable $k$-polymatroids by its excluded minors appears to be quite a difficult task, even for the simplest nontrivial case when $k = 2$ and $\mathbb{F} = \mathbb{F}_2$, i.e. the class of \textit{binary} $2$-polymatroids. We consider a simpler variation of this problem by first assigning a unique matroid to each $k$-polymatroid, called its \textit{$k$-natural matroid}. Geometrically, the \textit{$k$-natural matroid} of a polymatroid $(E,\rho)$ is obtained by replacing each $e \in E$ with $k$ points lying freely in $e$. A similar notion is the \textit{natural matroid} of $(E,\rho)$, obtained by replacing each $e \in E$ by $\rho(e)$ points lying freely in $E$.

If $\mathcal{C}$ is a minor-closed class of matroids, then the following are also minor-closed:
\begin{enumerate}
    \item The class $\mathcal{C}'$ of polymatroids whose natural matroids are in $\mathcal{C}$.
    \item The class $\mathcal{C}'_k$ of $k$-polymatroids in $\mathcal{C}'$.
    \item The class $\widetilde{\mathcal{C}}'_k$ of $k$-polymatroids whose $k$-natural matroids are in $\mathcal{C}$.
\end{enumerate}
Let $\mathcal{C}$ be the class of binary matroids, which has a single excluded minor: $U_{2,4}$. In this paper, we investigate the question: for a fixed $k \geq 2$, what is the set of excluded minors for $\widetilde{\mathcal{C}}'_k$? (When $k = 1$, this is simply the class of binary matroids.) From now on, we will refer to $\widetilde{\mathcal{C}}'_k$ as $\mathcal{P}_{U_{2,4}}^k$. This notation is inspired by \cite{bonin long}, where authors Joseph E. Bonin and Kevin Long determined the set of excluded minors for $\mathcal{P}_{U_{2,4}}$, the class of $2$-polymatroids whose natural matroids are binary. They found an infinite sequence of excluded minors related by \textit{compression}, along with eight other excluded minors that do not belong to this sequence. Since polymatroids with parallel points cannot be excluded minors for $\mathcal{P}_{U_{2,4}}^2$, it must be that the set of excluded minors for $\mathcal{P}_{U_{2,4}}$ is equal to the set of excluded minors for $\mathcal{P}_{U_{2,4}}^2$. Thus, they have resolved the question for the case $k = 2$. For each $k \geq 3$, we find the set of excluded minors for $\mathcal{P}_{U_{2,4}}^k$ and show that it is finite.

In Section \ref{section:definitions and background}, we provide the necessary definitions and background and outline a general strategy for finding the set of excluded minors for $\mathcal{P}_{U_{2,4}}^k$, which involves first exhausting all possibilities on a ground set of size $4$ or smaller, and then analyzing \textit{decompressions} of excluded minors whose ground sets have size $4$. The former is done in Sections \ref{section:ground set 1} through \ref{section:ground set 4}, and the latter is done in Section \ref{section:decompressions}. Finally, in Appendix \ref{appendix: excluded minors}, we provide a list of the excluded minors and their main properties.

\section*{Acknowledgement}
The author is grateful to Joe Bonin for suggesting this problem and for many helpful discussions.

\section{Definitions and background}
\label{section:definitions and background}

Much of the exposition in this section comes directly from \cite{bonin long}; we have chosen to reprint it for the reader's convenience and for the sake of completeness. Some proofs are omitted and for those we encourage the reader to refer directly to \cite{bonin long}. We follow Oxley \cite{oxley} for standard matroid terminology and notation. Many matroid notions generalize nicely to polymatroids.

\subsection{Polymatroids}

As is the case with matroids, there are many \textit{cryptomorphic}\footnote{This adjective is used to describe two objects which are equivalent but not in an obvious manner. Although its definition applies generally, it occurs most frequently in the matroid-theoretic literature.} ways of defining polymatroids. In this paper, we primarily use the rank function definition given at the beginning of Section \ref{section:introduction}. To keep the notation less cumbersome, we commonly omit the use of curly brackets when referring to subsets of $E$. For example, we might refer to the set $\{e, f, g\}$ as simply $efg$. We also frequently pass to the geometric representation of a polymatroid \cite{oxley}, especially when it aids intuition. When the context is clear, we often refer to a polymatroid $(E,\rho)$ as simply $\rho$. We might also refer to the ground set of a polymatroid $\rho$ as $E(\rho)$. 

\begin{figure}[H]
\centering
    \includegraphics[width=\textwidth]{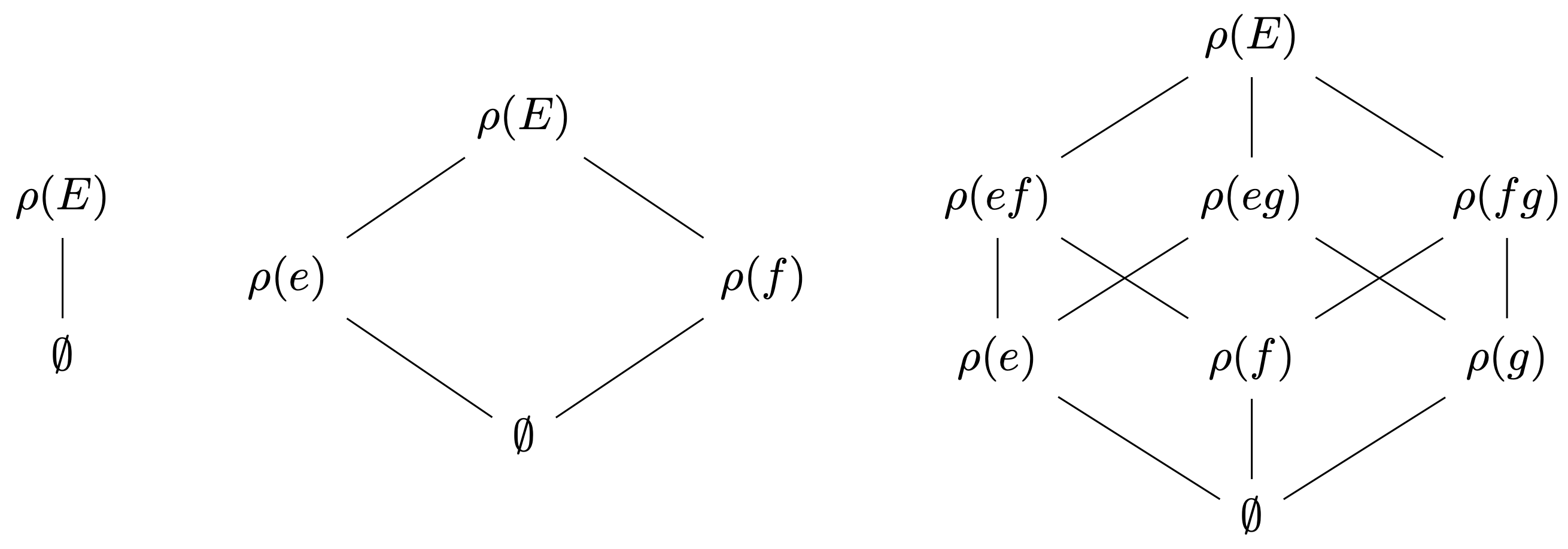}
    \caption{Generic Hasse diagrams displaying the values of $\rho$ on $2^E$ where, from left to right, $E = \{e\}$, $E = \{ef\}$, and $E = \{e, f, g\}$.}
    \label{figure: hasse diagrams 1 2 3}
\end{figure}

\begin{figure}[H]
\centering
\includegraphics[width=\textwidth]{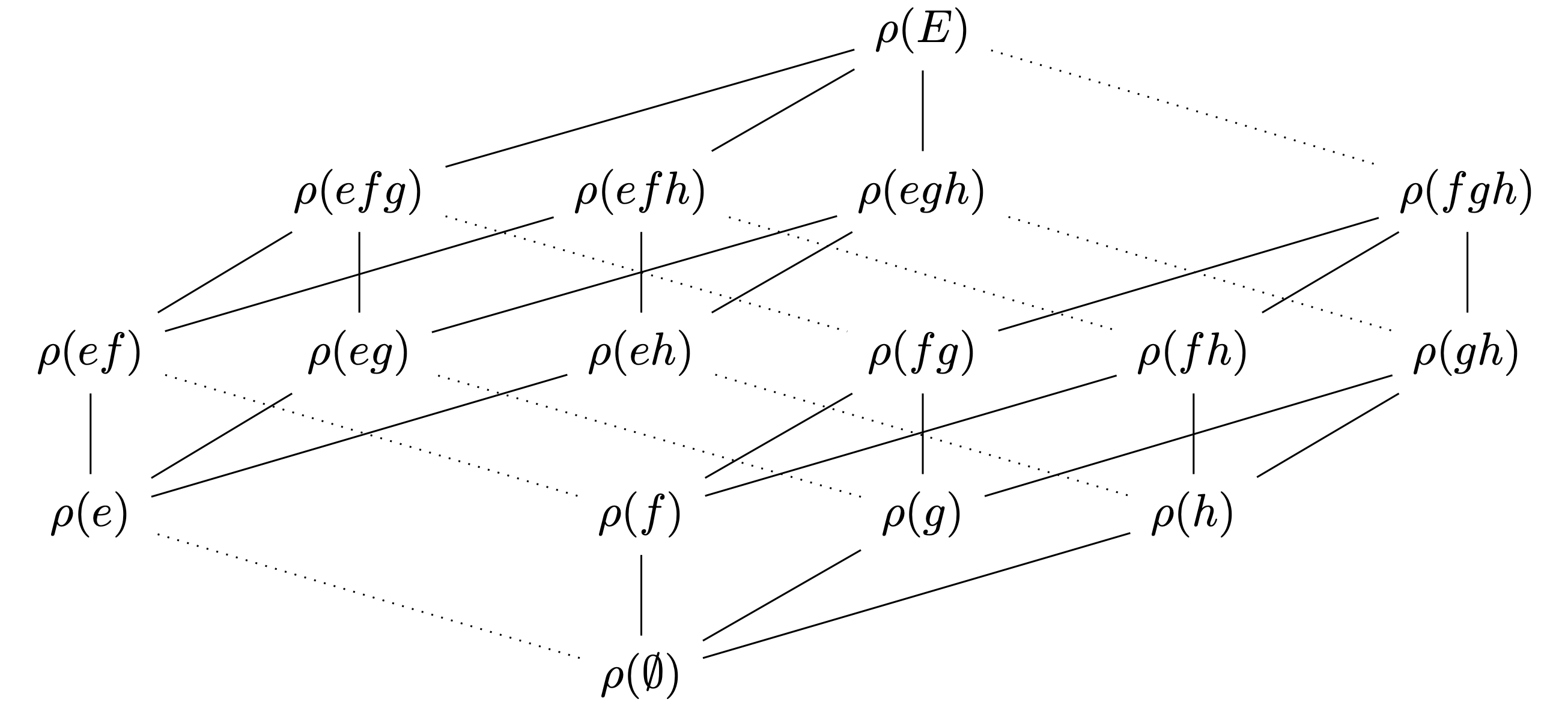}
\caption{A generic Hasse diagram displaying the values of $\rho$ on $2^E$ where $E = \{e, f, g, h\}$.}
    \label{figure: hasse diagrams 4}
\end{figure}

We will frequently use Hasse diagrams (see Figures \ref{figure: hasse diagrams 1 2 3} and \ref{figure: hasse diagrams 4}) to display the values of $\rho$ on $2^E$. For $X \subseteq E$, if there are multiple possibilities for the value of $\rho(X)$ based on the available analysis, the various options for $\rho(X)$ will be stacked vertically, in descending order from top to bottom.


\textbf{Minors.} Let $(E,\rho)$ be a polymatroid and let $X \subseteq E$. 
\begin{enumerate}
    \item The \textit{deletion} of $X$ is the polymatroid $(E-X, \rho\backslash X)$ where 
    \begin{displaymath}
        (\rho\backslash X)(Y) = \rho(Y)
    \end{displaymath}
    for all $Y \subseteq E-X$. The deletion of $X$ from $\rho$ is equivalent to the \textit{restriction} of $\rho$ to $E-X$, denoted $(E-X, \rho|_{E-X})$.
    \item The \textit{contraction} of $X$ is the polymatroid $(E-X, \rho/X)$ where
    \begin{displaymath}
        (\rho/X)(Y) = \rho(X \cup Y)-\rho(X)
    \end{displaymath}
    for all $Y \subseteq E-X$. 
\end{enumerate} 
The operations of deletion and contraction commute when applied to disjoint subsets; i.e if $X, Y \subseteq E$ and $X \cap Y = \emptyset$, then $(\rho\backslash X)/Y = (\rho/Y)\backslash X$. Furthermore, for such $X$ and $Y$, we have $(\rho\backslash X)\backslash Y = \rho\backslash (X\cup Y)$ and $(\rho/X)/Y = \rho/(X \cup Y)$. A polymatroid $\rho'$ is a \textit{minor} of $\rho$ if $\rho' = (\rho\backslash X)/Y$ for some disjoint subsets $X$ and $Y$ of $E$. That is, we obtain $\rho'$ from $\rho$ via a sequence of deletions and contractions. Depending on the context, we might emphasize that $\rho'$ is a \textit{polymatroid} minor of $\rho$. 

\textbf{Minor-closed classes.} A class of polymatroids is \textit{minor-closed} if for any $\rho$ in the class, all minors of $\rho$ are also in the class. We can completely characterize a minor-closed class of polymatroids by its set of \textit{excluded minors}: those polymatroids which are not in the class, whose proper minors are all in the class.

\textbf{Geometry.} We borrow terminology from projective geometry when describing polymatroids; this is consistent with the geometric representation view of matroids \cite{oxley} and polymatroids. \begin{enumerate}
    \item We say that $e, f \in E$ are \textit{parallel} if $0 < \rho(e) = \rho(f) = \rho(ef)$. In a matroid, there can be parallel points; in a $k$-polymatroid, we can have parallel points, parallel lines, parallel planes, etc. up to parallel elements of rank $k$. 
    \item For $X \subseteq E$, if $\rho(X) = 2$, then the elements of $X$ are \textit{collinear}, and if $\rho(X) = 3$, then the elements of $X$ are \textit{coplanar}.
    \item For $e, f \in E$, if $\rho(e),\rho(f)>0$ and $\rho(e)+\rho(f) = \rho(ef)$, then $e$ and $f$ are \textit{skew}.
    \item If $0 < \rho(e) < \rho(f)$ and $\rho(f) = \rho(ef)$, then $e$ \textit{lies on} $f$.
\end{enumerate}

\textbf{Connectedness.} Let $(E_1,\rho_1)$ and $(E_2,\rho_2)$ be polymatroids, with $E_1 \cap E_2 = \emptyset$. Their \textit{direct sum} is the polymatroid $(E_1 \cup E_2, \rho_1 \oplus \rho_2)$, where for $X \subseteq E_1 \cup E_2$,
\begin{displaymath}
    (\rho_1 \oplus \rho_2)(X) = \rho_1(X \cap E_1) + \rho_2(X \cap E_2).
\end{displaymath}
A polymatroid is \textit{connected} if it is not isomorphic to the direct sum of two nonempty polymatroids; otherwise it is \textit{disconnected}.

\textbf{Closure.} The \textit{closure} of $X \subseteq E$ is the set
\begin{displaymath}
    \text{cl}(X) := \{e \in E: \rho(X \cup e) = \rho(X)\}.
\end{displaymath}

\subsection{The \texorpdfstring{$k$}{}-natural matroid}
Let $(E,\rho)$ be a polymatroid. For each $e \in E$, let $X_e$ be a set of $k$ elements. If $e, f \in E$ are distinct, then we require $X_e \cap X_f = \emptyset$. For any subset $A \subseteq E$, we define 
        \begin{displaymath}
            X_A := \bigcup_{e \in A} X_e.
        \end{displaymath}
        The \textit{$k$-natural matroid} $M_\rho^k$ of $\rho$ is the matroid $(X_E, r)$ where 
        \begin{displaymath}
            r(X) := \min\{\rho(A) + |X-X_A|: A \subseteq E\}.
        \end{displaymath}

By a result of McDiarmid \cite{mcdiarmid}, we see that $M_\rho^k = (X_E, r)$ is indeed a matroid. For $x,y \in X_E$, consider the transposition $\varphi|_{(x,y)}$ which swaps $x$ and $y$ while fixing all other elements of $X_E$. If $\varphi|_{(x,y)}$ is an automorphism of $M_\rho^k$, then we say that $x$ and $y$ are \textit{clones}. Observe that for any $e \in E$, any pair of distinct elements of $X_e$ are clones. Hence, we say that $X_e$ is \textit{set of clones}. Geometrically, $M_\rho^k$ is obtained from $\rho$ by replacing each $e \in E$ with $k$ points lying freely in $e$. We state the following lemma (without proof) and observations from \cite{bonin long}:
\begin{lemma}\label{k-natural is unique}
    Let $\rho$, $E$, $X_e$, and $X_A$ be as above. A matroid $M$ on $X_E$ is $M_\rho^k$ if and only if each set $X_e$ is a set of clones and $r_M(X_A) = \rho(A)$ for all $A \subseteq E$.
\end{lemma}

\begin{observation}\label{k-natural observations}\leavevmode
\begin{enumerate}
    \item If $\rho$ and $\rho'$ are polymatroids on $E$ with $M_\rho^k = M_{\rho'}^k$ where for each $e \in E$ the corresponding set $Y_e$ is the same in both $k$-natural matroids, then $\rho = \rho'$. 
    \item For any polymatroid $\rho$, we have $M_{\rho\backslash e}^k = M_{\rho}^k\backslash X_e$ and $M_{\rho/e}^k = M_\rho^k/X_e$ for all $e \in E$.
\end{enumerate}
\end{observation}

By Observation \ref{k-natural observations}(2) above, if $\mathcal{C}$ is a minor-closed class of matroids, then $\widetilde{\mathcal{C}}'_k$, the class of $k$-polymatroids whose $k$-natural matroids are in $\mathcal{C}$, is also closed, since minors of $k$-polymatroids are $k$-polymatroids. Now, we will state and prove two lemmas which will greatly simplify our analysis in Sections \ref{section:ground set 1} through \ref{section:decompressions}.

\begin{lemma}\label{simplification}
Let $\mathcal{C}$ be a minor-closed class of matroids whose excluded minors are all simple. If the polymatroid $\rho$ is an excluded minor for $\widetilde{\mathcal{C}}'_k$, then $\rho$ cannot have any loops (elements of rank $0$) or nontrivial parallel classes of points (elements of rank $1$).
\end{lemma}
\begin{proof}
Assume $\rho$ is a polymatroid which contains a loop or a nontrivial parallel class of points, and is not in $\widetilde{\mathcal{C}}'_k$. This means $M_\rho^k$ contains an excluded minor $M'$ for $\mathcal{C}$, where $M'$ is simple by assumption. If $\rho$ has a loop $e$, then $X_e$ is a set of $k$ loops in $M_\rho^k$. Since $M'$ is simple, it must be that $E(M') \cap X_e = \emptyset$. Therefore, $M_{\rho\backslash e}^k = M_\rho^k\backslash X_e$ would also contain $M'$ as a minor, implying $\rho\backslash e$ is not in $\widetilde{\mathcal{C}}'_k$, so $\rho$ cannot be an excluded minor for $\widetilde{\mathcal{C}}_k'$.  

If $\rho$ contains a pair of parallel points, say $e$ and $f$, then $X_e \cup X_f$ would be a subset of some parallel class in $M_\rho^k$ containing at least $2k$ points (note that $k \geq 1$). Since $M'$ is simple, it must be that $|E(M') \cap (X_e \cup X_f)| \leq 1$. If $E(M') \cap (X_e \cup X_f) = \emptyset$, then the previous argument for loops applies. If $E(M') \cap (X_e \cup X_f)$ is a singleton $x$, then without loss of generality we can assume $x \in X_e$. Then $M_{\rho\backslash f}^k = M_\rho^k\backslash X_f$ would also contain $M'$ as a minor, implying $\rho\backslash f$ is not in $\widetilde{\mathcal{C}}'_k$, so $\rho$ cannot be an excluded minor for $\widetilde{\mathcal{C}}_k'$. 
\end{proof}

Let $(E,\rho)$ be a polymatroid. We will let $S(\rho)$ denote the \textit{simplification} of $\rho$, defined as the polymatroid obtained by deleting all loops of $\rho$, and in each nontrivial parallel class of points of $\rho$, deleting all points except for one. 

\begin{figure}[H]
    \centering
    \includegraphics[scale = 0.1]{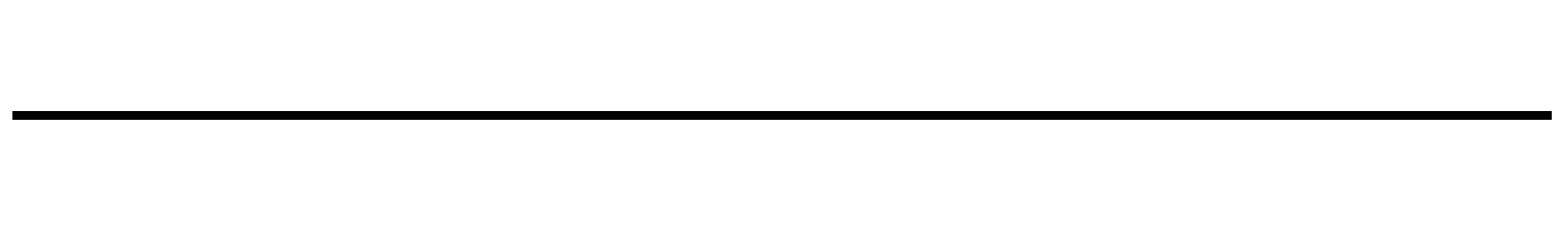}\hfill
    \includegraphics[scale = 0.1]{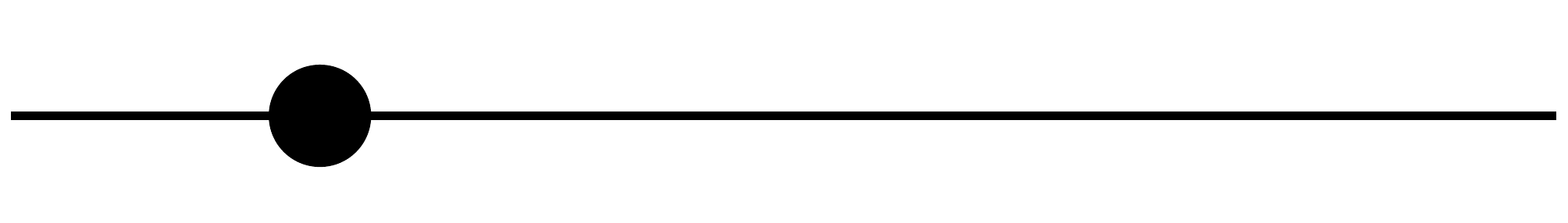}\hfill
    \includegraphics[scale = 0.1]{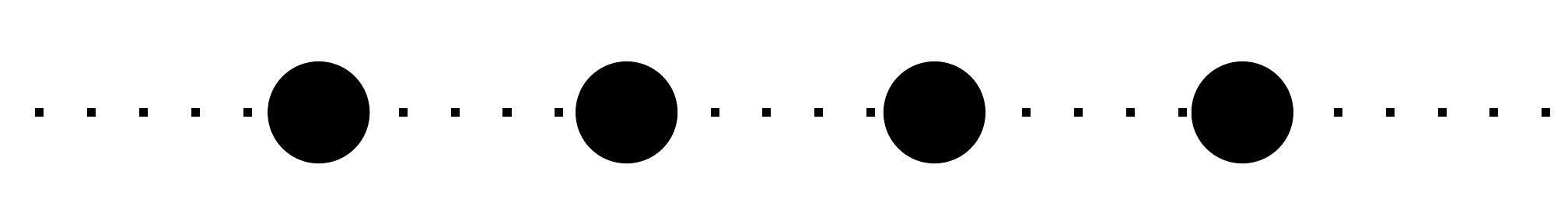}
    \caption{Left: A line, which is also the excluded minor $Ex^2$ for $\mathcal{P}_{U_{2,4}}^k$ where $k \geq 4$. Center: A point lying on a line, which is also the excluded minor $Ex_\gamma^2$ for $\mathcal{P}_{U_{2,4}}^3$. Right: $U_{2,4} = S(M_{U_{2,4}}^k)$, which is an excluded minor for $\mathcal{P}_{U_{2,4}}^k$ for all $k \geq 1$. }
    \label{fig: polymatroids with the same k-natural matroid up to simplification}
\end{figure}

\begin{example}\label{two polymatroids with same simplified k-natural matroid}
Different $k$-polymatroids may have the same $k$-natural matroid up to simplification of the latter. Let $\rho_1$ be the polymatroid consisting of a single line, and let $\rho_2$ be the polymatroid consisting of a point lying on a line. Then $M_{\rho_1}^k$ is equal to $S(M_{U_{2,4}}^k)$ if and only if $k = 4$, and $M_{\rho_2}^k$ is equal to $S(M_{U_{2,4}}^k)$ if and only if $k = 3$. We will see in Section \ref{section:ground set 1} that $\rho_1$ is the excluded minor $Ex^2$ for $\mathcal{P}_{U_{2,4}}^k$ where $k \geq 4$ and in Section \ref{section:ground set 2} that $\rho_2$ is the excluded minor $Ex_\gamma^2$ for $\mathcal{P}_{U_{2,4}}^3$. See Figure \ref{fig: polymatroids with the same k-natural matroid up to simplification}.
\end{example}

It is clear that $S(\rho)$ is a restriction of $\rho$. If $\rho = S(\rho)$, then we say that $\rho$ is \textit{simple}. If $\rho \neq S(\rho)$, then $S(\rho)$ is a proper minor of $\rho$, so if $S(\rho)$ is not in $\widetilde{\mathcal{C}}'_k$, then neither is $\rho$. Going forward, Lemma \ref{simplification} allows us to automatically rule out polymatroids containing loops or nontrivial parallel classes of points as candidates for excluded minors for $\mathcal{P}_{U_{2,4}}^k$. Furthermore, when analyzing matroid minors of $M_\rho^k$ to see if any of them are isomorphic to $U_{2,4}$, it will suffice to consider minors of $S(M_\rho^k)$ instead.

\begin{lemma}\label{connectedness}
Let $\mathcal{C}$ be a minor-closed class of matroids characterized exclusively by connected excluded minors. Let $\rho_1, \hdots, \rho_m$ be $k$-polymatroids where $\rho_i \in \widetilde{\mathcal{C}}'_k$ for all $i$. Then $\rho_1\oplus \hdots \oplus \rho_m\in \widetilde{\mathcal{C}}'_k$. (Equivalently, $\widetilde{\mathcal{C}}'_k$ is closed under the operation of direct sum.) \end{lemma}
\begin{proof}
    Let $m = 2$. Observe that $M_{\rho_1 \oplus \rho_2}^k = M_{\rho_1}^k \oplus M_{\rho_2}^k$. It suffices to show $M_{\rho_1}^k \oplus M_{\rho_2}^k \in \mathcal{C}$. Any minor $M'$ of $M_{\rho_1}^k \oplus M_{\rho_2}^k$ must be isomorphic to $M_1' \oplus M_2'$ where $M_i'$ is a minor of $M_{\rho_i}^k$. Now assume $M'$ is an excluded minor for $\mathcal{C}$. Since $\mathcal{C}$ is characterized exclusively by connected excluded minors, it must be that $M_1'$ is empty or $M_2'$ is empty. Without loss of generality, assume $M_2'$ is empty, so $M_1' \oplus M_2' = M_1'$. Since $M_1'$ is a minor of $M_{\rho_1}^k$, this implies $\rho_1 \notin \widetilde{\mathcal{C}}'_k$, a contradiction. Hence, no minor of $M_{\rho_1}^k \oplus M_{\rho_2}^k$ is an excluded minor for $\mathcal{C}$, implying $M_{\rho_1}^k \oplus M_{\rho_2}^k \in \mathcal{C}$. The result for $m > 2$ follows by induction.
\end{proof}

The class of binary matroids is characterized exclusively by the excluded minor $U_{2,4}$, which is a connected matroid. Thus, if a $k$-polymatroid $\rho$ is disconnected, Lemma \ref{connectedness} allows us to immediately rule out $\rho$ as an excluded minor for $\mathcal{P}_{U_{2,4}}^k$.

\textbf{$k$-Duality.} For a matroid $(E,r)$, its \textit{dual matroid} $(E, r^*)$ is characterized by 
\begin{displaymath}
r^*(X) = |X|+r(E-X) - r(E),
\end{displaymath}
and the map $r \mapsto r^*$ is the only involution on the class of matroids that swaps deletion and contraction \cite{kung}. That is, $(r\backslash e)^* = r^*/e$ and $(r/e)^* = r^*\backslash e$ for all $e \in E$. Inspired by the utility of this property in excluded minor problems, we desire a notion of $k$-polymatroid duality that behaves nicely with the operations of polymatroid deletion and contraction. As in \cite{bonin long}, we will use $k$-duality. For a $k$-polymatroid $(E,\rho)$, its $k$-dual is the $k$-polymatroid $(E, \rho^*)$ where
\begin{displaymath}
    \rho^*(X): = k|X| + \rho(E-X)-\rho(E).
\end{displaymath}
If $(E, \rho) \cong (E,\rho^*)$, we say $\rho$ is \textit{self-$k$-dual} and if $(E,\rho) = (E,\rho^*)$, we say $\rho$ is \textit{identically self-$k$-dual}.

Indeed, $\rho \mapsto \rho^*$ is the only involution on the class of $k$-polymatroids that swaps deletion and contraction \cite{whittle}. Therefore, if a minor-closed class of $k$-polymatroids is closed under $k$-duality, then so is its set of excluded minors.

We state the following lemmas from \cite{bonin long} without proof:
\begin{lemma}\label{k-dual of k-natural}
    Let $(E,\rho)$ be a $k$-polymatroid and $(E,\rho^*)$ be its $k$-dual. Then the $k$-natural matroid of $\rho$ is dual to that of $\rho^*$, i.e. $(M_\rho^k)^* = M_{\rho^*}^k$.
\end{lemma}

\begin{lemma}
    Let $(E,\rho)$ be a $k$-polymatroid and $(E,\rho^*)$ be its $k$-dual. If $e \in E$, then $M_{\rho\backslash e}^k = \left(M_{\rho^*/ e}^k\right)^*$ and $M_{\rho/ e}^k = \left(M_{\rho^*\backslash e}^k\right)^*$.
\end{lemma}

The next theorem is a straightforward generalization of its counterpart in \cite{bonin long}.

\begin{theorem}\label{kduality}
Let $\mathcal{C}$ be a minor-closed, dual-closed class of matroids. Then $\widetilde{\mathcal{C}}'_k$ is closed under $k$-duality, as is its set of excluded minors.
\end{theorem}
\begin{proof}
    Let $\rho$ be a $k$-polymatroid in $\widetilde{\mathcal{C}}'_k$, which implies $M_\rho^k \in \mathcal{C}$. Since $\mathcal{C}$ is dual-closed, $(M_\rho^k)^* \in \mathcal{C}$. By Lemma \ref{k-dual of k-natural}, $(M_\rho^k)^* = M_{\rho^*}^k$, so this implies $M_{\rho^*}^k \in \mathcal{C}$, giving us $\rho^* \in \widetilde{\mathcal{C}}'_k$ as desired for the first assertion. The second assertion follows.
\end{proof}

\begin{figure}[H]
    \centering
    \includegraphics[scale = 0.2]{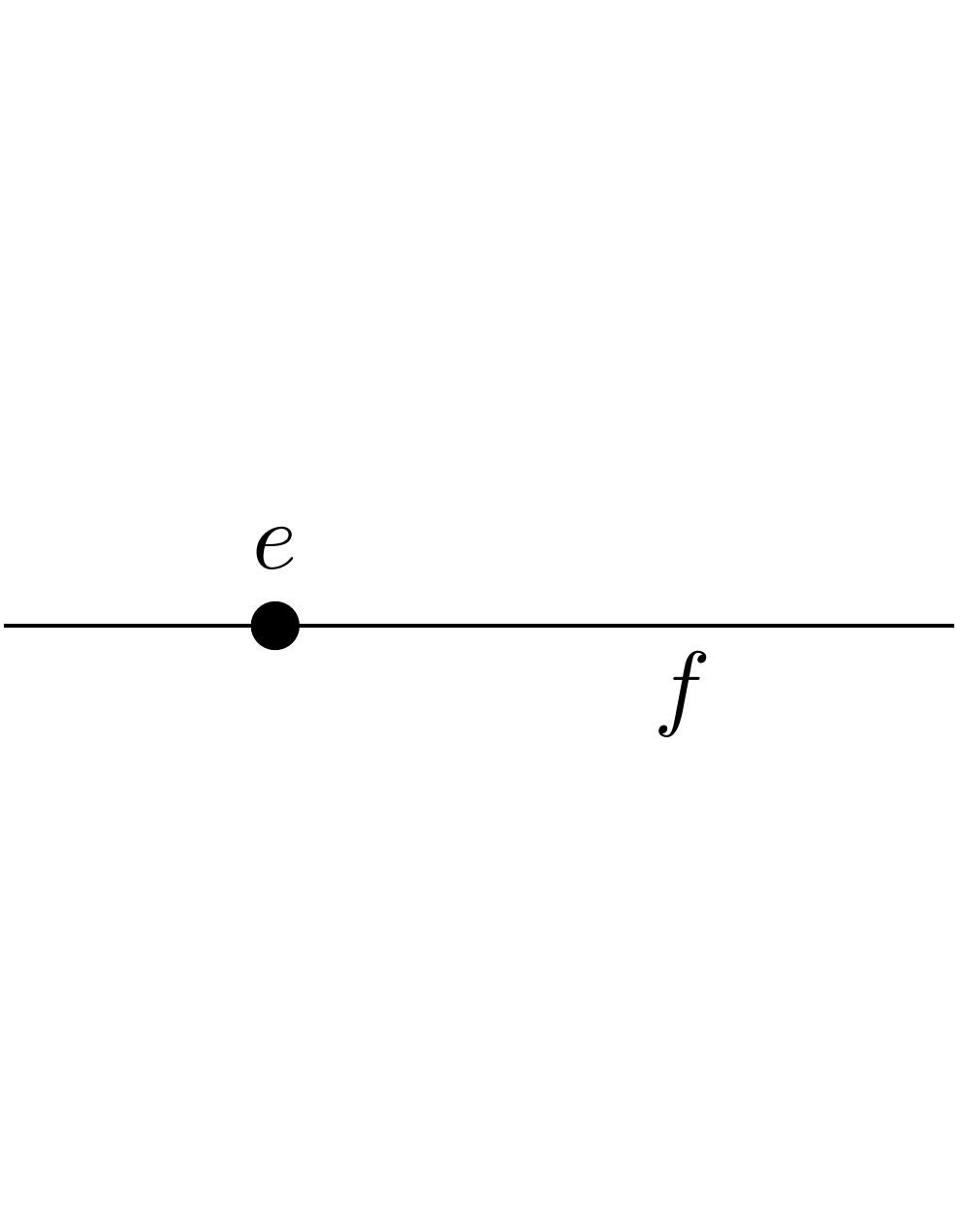}\hspace{10ex}
    \includegraphics[scale = 0.15]{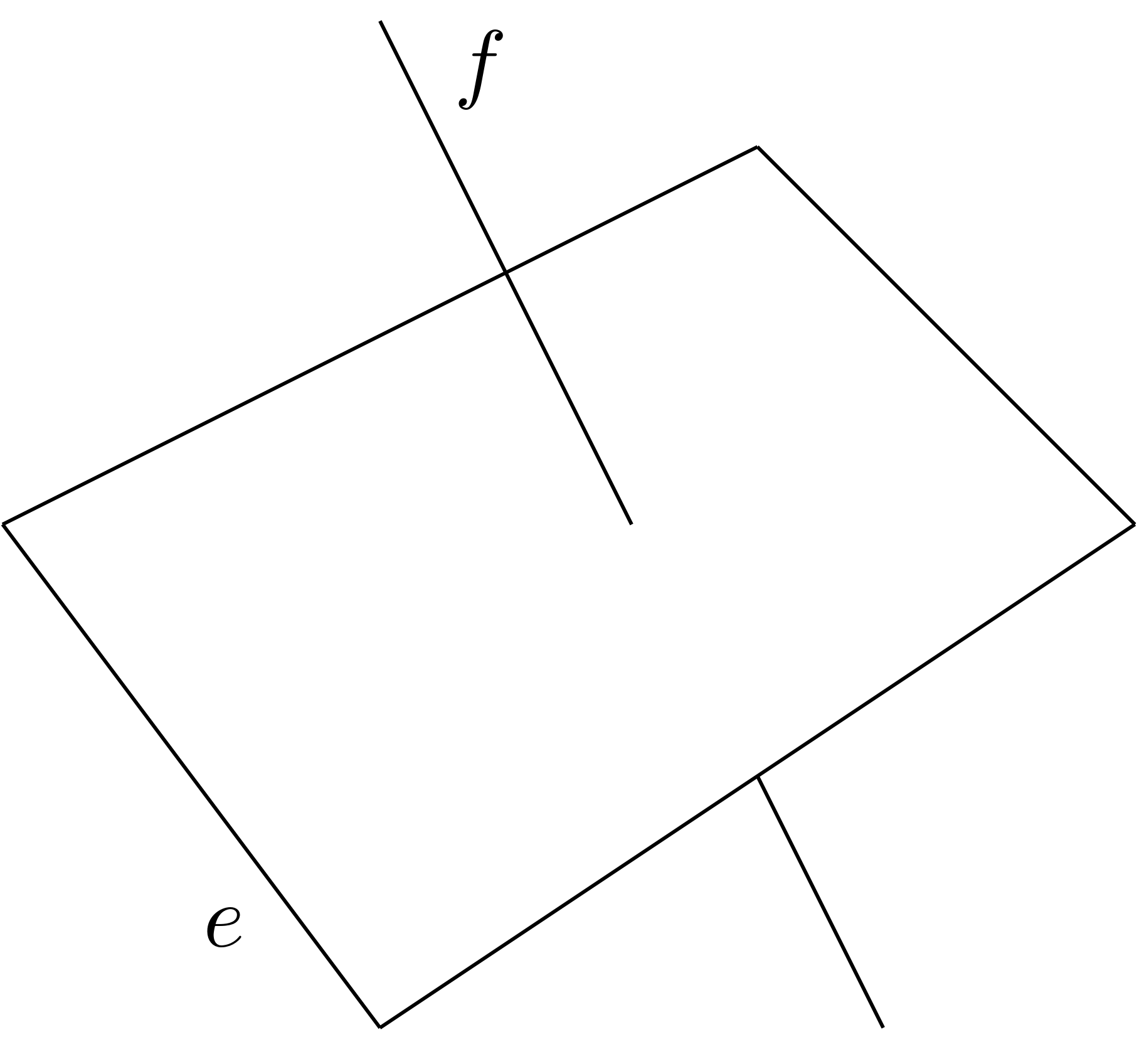}
    \caption{Left: The polymatroid $Ex_\gamma^2$ consisting of a point lying on a line. Right: The $3$-dual of $Ex_\gamma^2$, denoted $Ex_\epsilon^4$, where $e$ becomes a plane and $f$ remains a line.}
    \label{fig:kduality example}
\end{figure}

\begin{example}\label{kduality example}
    Consider $Ex_\gamma^2$ consisting of a point lying on a line, first introduced in Example \ref{two polymatroids with same simplified k-natural matroid}. Its $3$-dual consists of a line and a plane spanning rank $4$, and is also an excluded minor for $\mathcal{P}_{U_{2,4}}^3$, denoted $Ex_\epsilon^4$. See Figure \ref{fig:kduality example}.
\end{example}

\subsection{Compression}

Let $(E,\rho)$ be a $k$-polymatroid. For $1 \leq l \leq k-1$, the \textit{$l$-compression of $\rho$ by $e \in E$} is the polymatroid $\rho_{\downarrow e}^l$ obtained by freely adding $l$ points $e_1,\hdots, e_l$ to $e$, then contracting $e_1, \hdots, e_l$ and deleting $e$. If $e$ is a point, then $\rho_{\downarrow e}^l = \rho/e$ for any $l$. If $e$ is a loop, then $\rho_{\downarrow e}^l = \rho\backslash e$ for any $l$. If $\rho' = \rho_{\downarrow e}^l$ for some $e$ and $l$, then we say that $\rho'$ is a \textit{compression} of $\rho$ (and $\rho$ is a \textit{decompression} of $\rho'$). If $l$ is known, we can also say that $\rho'$ is an \textit{$l$-compression} of $\rho$ (and $\rho$ is an \textit{$l$-decompression} of $\rho'$).

\begin{figure}[H]
    \centering
    \includegraphics[scale = 0.25]{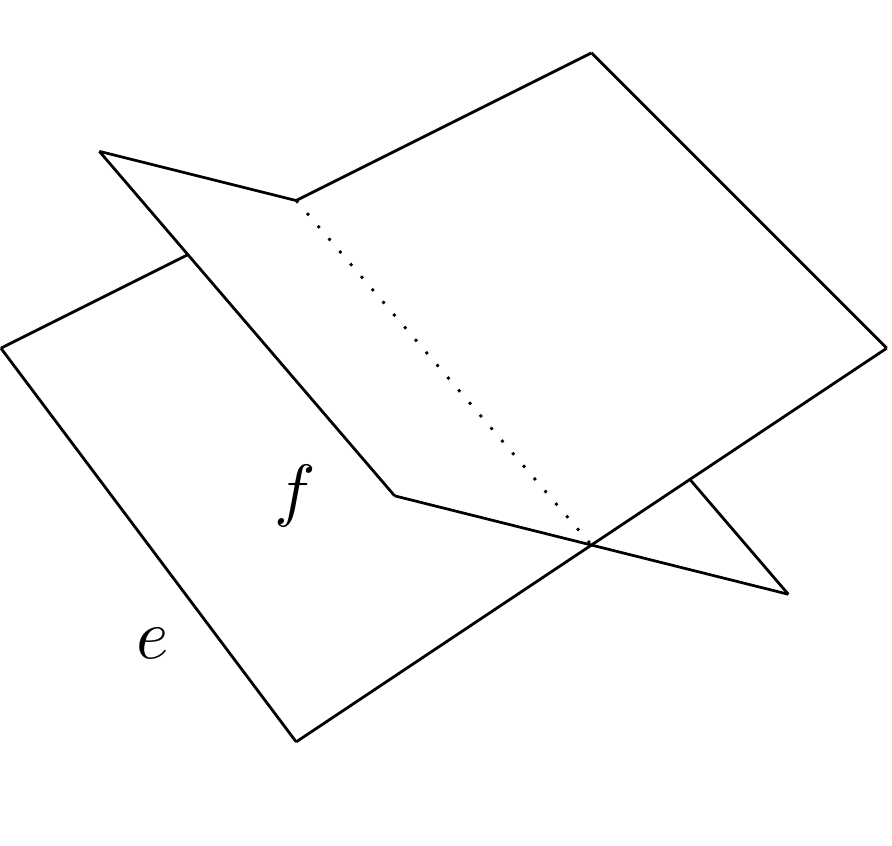}\hfill \includegraphics[scale = 0.25]{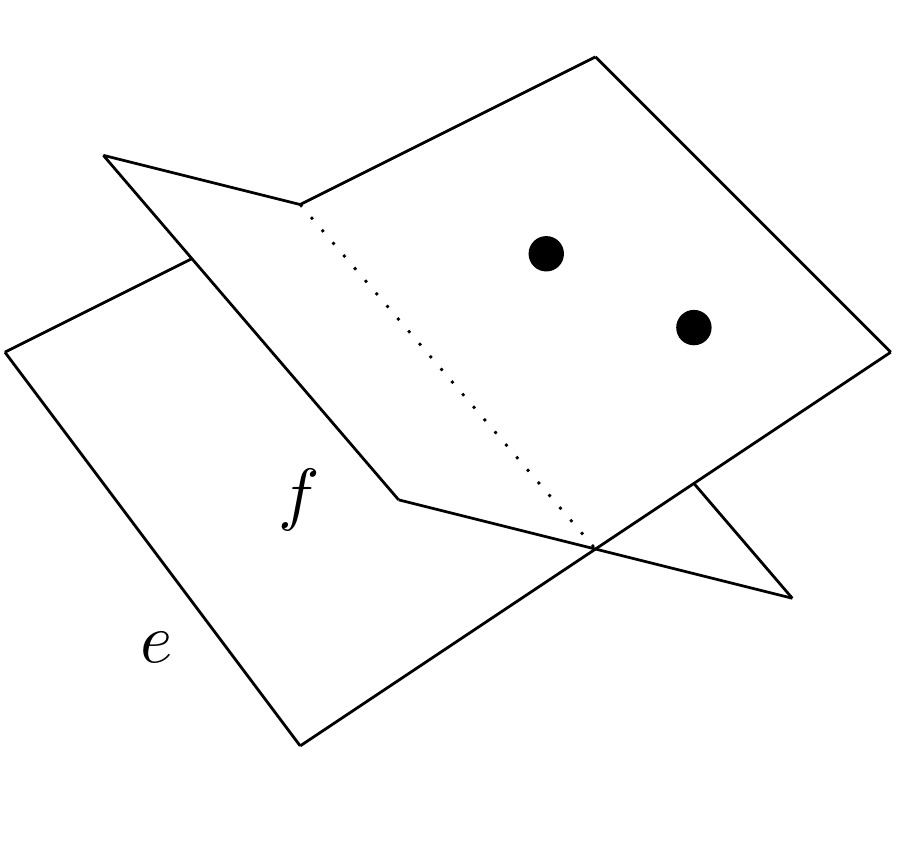}\hfill \includegraphics[scale = 0.15]{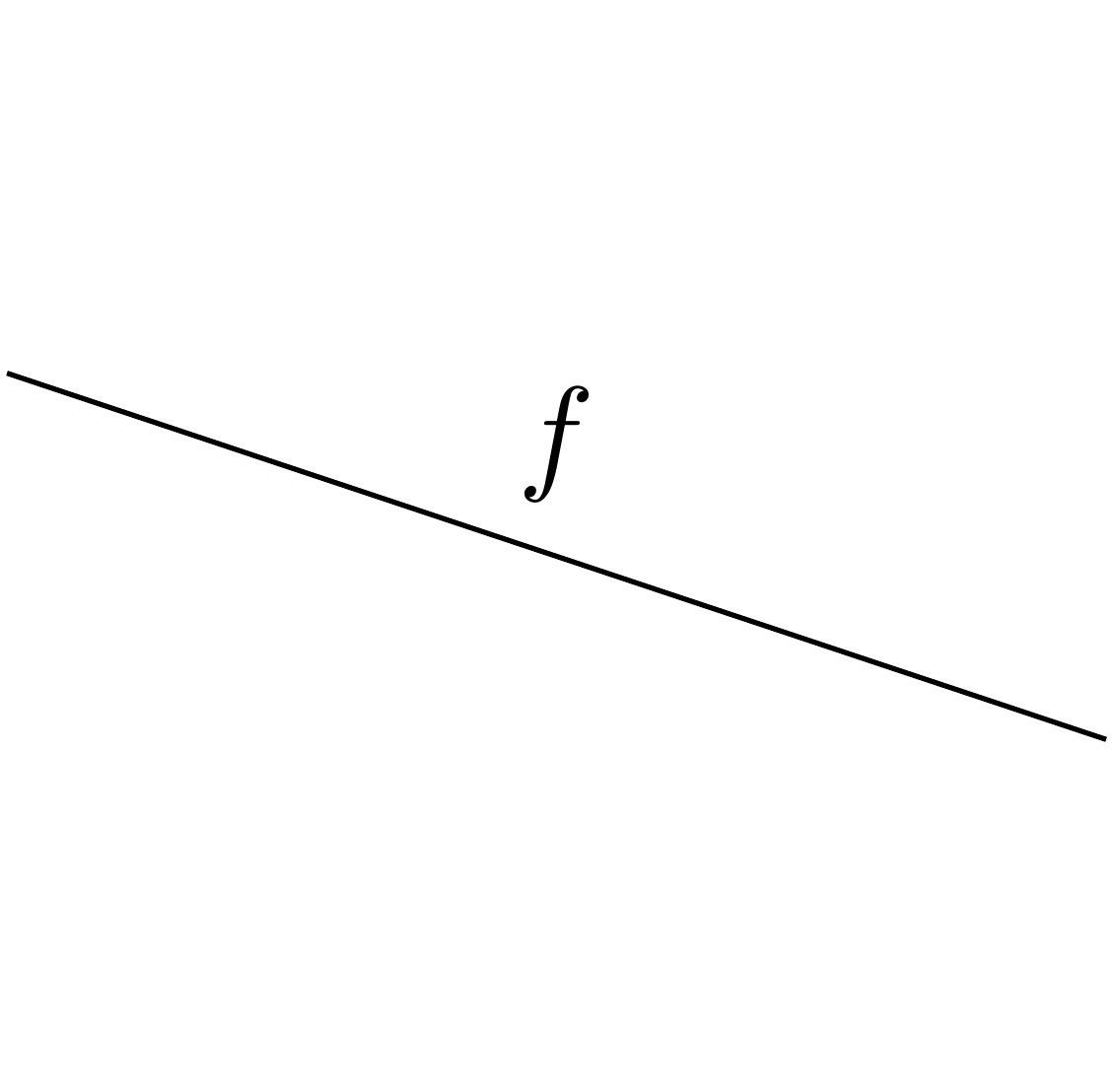}
    \caption{Left: The polymatroid $Ex_\epsilon^4$ consisting of two planes $e$ and $f$ spanning rank $4$. Center: To take the $2$-compression of $Ex_\epsilon^4$ by $e$, first we must freely add $2$ points to $e$. Right: After contracting on the two points and deleting $e$, we are left with a line.}
    \label{fig:compression example}
\end{figure}

\begin{example}\label{compression example}
    By definition, excluded minors cannot be related to one another via taking proper minors, but they \textit{can} be related via the operation of taking $l$-compressions. Let $\rho$ be the polymatroid consisting of two planes $e$ and $f$ spanning rank $4$. This is the excluded minor $Ex_\epsilon^4$ for $\mathcal{P}_{U_{2,4}}^4$. The $2$-compression of $\rho$ by $e$ gives a line, which is the excluded minor $Ex^2$ for $\mathcal{P}_{U_{2,4}}^4$. See Figure \ref{fig:compression example}.
\end{example}

The compression operation is crucial to finding the excluded minors for $\mathcal{P}_{U_{2,4}}^k$ as the size of the ground set increases. This will be illuminated in Lemma \ref{compression lemma} and the discussion following. The next few results generalize their counterparts in \cite{bonin long}.

\begin{lemma}\label{k-natural of l-compression}
    For a $k$-polymatroid $(E,\rho)$ and $1 \leq l \leq k-1$, fix $e \in E$ with $\rho(e) > 0$ and fix $e_1,\hdots,e_l \in X_e$. The $k$-natural matroid of $\rho_{\downarrow e}^l$ is $M_\rho^k/\{e_1,\hdots,e_l\}\backslash (X_e-\{e_1,\hdots,e_l\})$.
\end{lemma}
\begin{proof}
    Let $N_e^k = M_\rho^k/\{e_1,\hdots,e_l\}\backslash (X_e-\{e_1,\hdots,e_l\})$. By Lemma \ref{k-natural is unique}, it suffices to show that 
    \begin{enumerate}
        \item[i.] For each $f \in E-e$, the set $X_f$ is a set of clones of $N_e^k$, and
        \item[ii.] $r_{N_e^k}(X_A) = \rho_{\downarrow e}^l(A)$ for all $A \subseteq E-e$.
    \end{enumerate}
    Property (i) holds since $X_f$ is a set of clones of $M_\rho^k$ and clones in a matroid remain clones in each minor that contains them. Property (ii) follows from the definition of $l$-compression.
\end{proof}

\begin{lemma}\label{compression lemma} Let $\mathcal{C}$ be a minor-closed class of matroids. Let the $k$-polymatroid $(E, \rho)$ be an excluded minor for $\widetilde{\mathcal{C}}'_k$ and fix $e \in E$ with $\rho(e) \geq 2$. Fix any $1 \leq l \leq k-1$. Then $\rho_{\downarrow e}^l$ is an excluded minor for $\widetilde{\mathcal{C}}'_k$ if and only if $\rho_{\downarrow e}^l \notin \widetilde{\mathcal{C}}'_k$.
\end{lemma}
\begin{proof}
The forward implication is immediate. For the converse, assume $\rho_{\downarrow e}^l \notin \widetilde{\mathcal{C}}'_k$. Fix $f \in E-e$. To show that $(\rho_{\downarrow e}^l)\backslash f$ is in $\widetilde{\mathcal{C}}'_k$, it suffices to show that $M_{(\rho_{\downarrow e}^l)\backslash f}^k$ is a minor of $M_{\rho\backslash f}^k$, which is in $\mathcal{C}$. Then by Observation \ref{k-natural observations}(2) and Lemma \ref{k-natural of l-compression},
\begin{align*}
    M^k_{(\rho_{\downarrow e}^l)\backslash f} & = M_{(\rho_{\downarrow e}^l)}^k\backslash X_f\\
    & = M_\rho^k/\{e_1, \hdots, e_l\} \backslash (X_e-\{e_1, \hdots, e_l\}\backslash X_f \\
    & = M_{\rho\backslash f}^k / \{e_1, \hdots, e_l\} \backslash (X_e - \{e_1, \hdots, e_l\}).
\end{align*}

Similarly, we will show that $(\rho_{\downarrow e}^l)/ f$ is in $\widetilde{\mathcal{C}}'_k$ by showing that $M_{(\rho_{\downarrow e}^l)/ f}$  is a minor of $M_{\rho/f}$, which is in $\mathcal{C}$. By Observation \ref{k-natural observations}(2) and Lemma \ref{k-natural of l-compression},
\begin{align*}
    M^k_{(\rho_{\downarrow e}^l)/ f} & = M_{(\rho_{\downarrow e}^l)}^k/ X_f\\
    & = M_\rho^k/\{e_1, \hdots, e_l\} \backslash (X_e-\{e_1, \hdots, e_l\}/X_f\\
    & = M_{\rho/ f}^k / \{e_1, \hdots, e_l\} \backslash (X_e - \{e_1, \hdots, e_l\}).
\end{align*} 
\end{proof}

The next result makes precise the interaction between the operations of taking an $l$-compression and taking the $k$-dual of a $k$-polymatroid. This will be useful in Section \ref{section:decompressions} where we analyze decompressions.

\begin{lemma}
\label{compression commutes with dual}
    Let $\mathcal{P}^k$ be the set of all $k$-polymatroids. Then
    \begin{displaymath}
        \left(\rho^l_{\downarrow e}\right)^* = \left(\rho^*\right)^{k-l}_{\downarrow e}.
    \end{displaymath}
\end{lemma}
\begin{proof}
    By Observation \ref{k-natural observations}(1), it suffices to show $M_{\left(\rho^l_{\downarrow e}\right)^*}^k = M_{\left(\rho^*\right)^{k-l}_{\downarrow e}}^k$. By Lemmas \ref{k-dual of k-natural} and \ref{k-natural of l-compression},
    \begin{align*}
        M_{\left(\rho^l_{\downarrow e}\right)^*}^k & = \left(M_{\rho^l_{\downarrow e}}^k\right)^*\\
        & = \left(M_{\rho}^k/\{e_1, \hdots, e_l\}\backslash (X_e-\{e_1, \hdots, e_l\})\right)^*\\
        & = \left(M_{\rho}^k\right)^*\backslash \{e_1, \hdots, e_l\}/(X_e-\{e_1, \hdots, e_l\})\\
        & = M_{\rho^*}^k\backslash \{e_1, \hdots, e_l\}/(X_e-\{e_1, \hdots, e_l\})\\
        & = M_{\left(\rho^*\right)^{k-l}_{\downarrow e}}^k.
    \end{align*}
\end{proof}

Let $\mathcal{C}$ be a minor-closed class of matroids and consider any $k$-polymatroid $(E,\rho)$ which is an excluded minor for $\widetilde{\mathcal{C}}'_k$. Let $Z_\rho$ be the set of compressions of $\rho$, i.e. the set of $k$-polymatroids $\rho'$ such that $\rho' = \rho_{\downarrow e}^l$ for some $e \in E$ and $1 \leq l \leq k-1$. Then either (1) $Z_\rho \subseteq \widetilde{\mathcal{C}}'_k$, or (2) there exists some $\rho' \in Z_\rho$ such that $\rho' \notin \widetilde{\mathcal{C}}'_k$.

Let us denote by $\Gamma$ the set of excluded minors $\rho$ such that (1) applies to $Z_\rho$. If $\rho\notin\Gamma$, i.e. if $\rho$ is such that (2) applies, then by Lemma \ref{compression lemma}, some sequence of compressions of elements $e_1,\hdots, e_m \in E$ each of rank $2$ or higher starting from $\rho$ eventually yields $\rho' \in \Gamma$. This grants us a two-part strategy for finding the excluded minors for $\widetilde{\mathcal{C}}'_k$:

\begin{strategy}\label{strategy}\leavevmode
\begin{enumerate}
    \item Find $\Gamma$.
    \item For each $\rho \in \Gamma$, find all decompressions of $\rho$ and determine which (if any) are excluded minors for $\widetilde{\mathcal{C}}'_k$. Repeat this step on all newly discovered excluded minors until the process terminates or we have characterized all excluded minors for $\widetilde{\mathcal{C}}'_k$, including infinite families of such.
\end{enumerate}
\end{strategy}

\begin{figure}[H]
    \centering
\includegraphics[width = \textwidth]{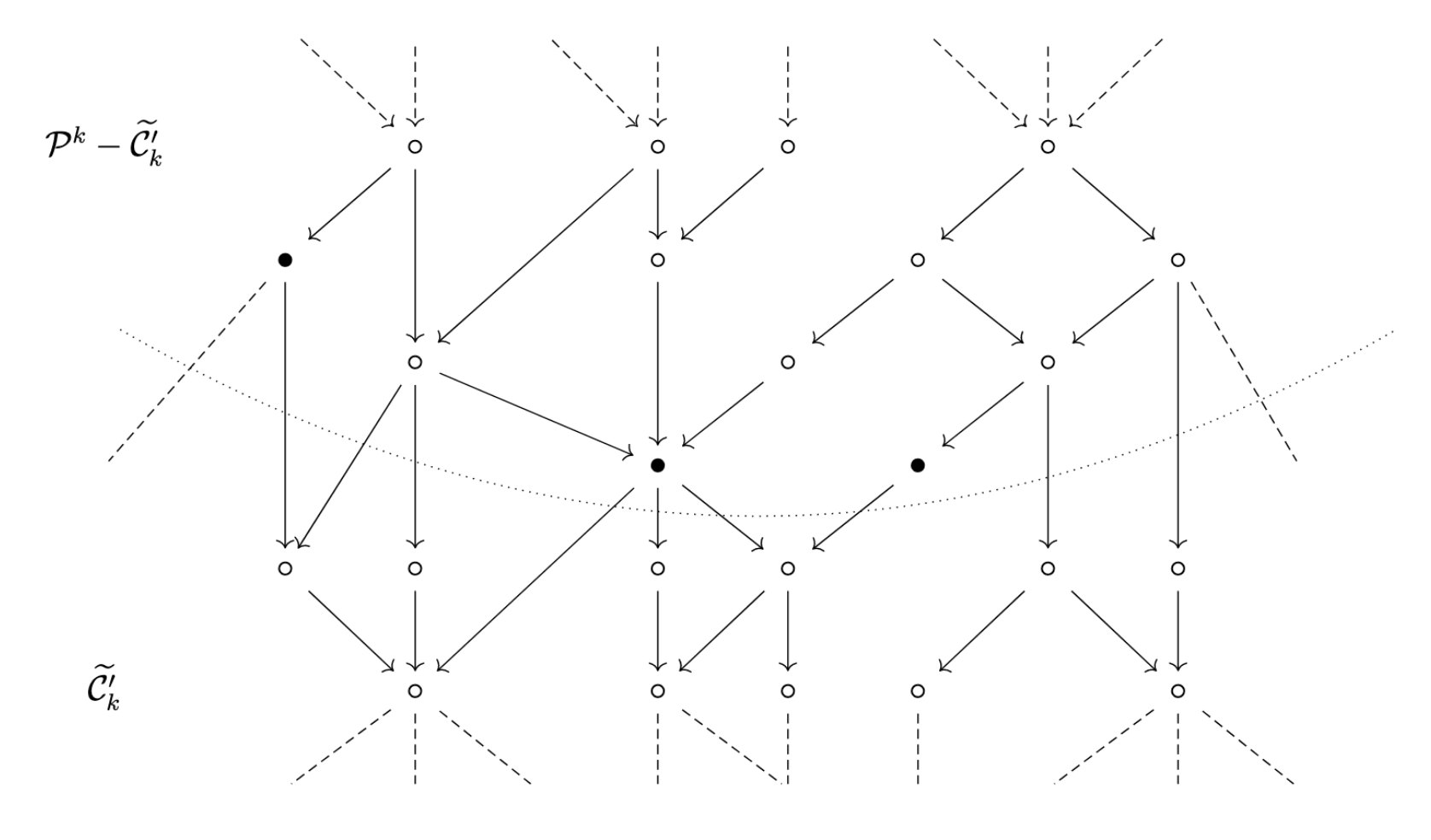}
    \caption{An intuitive sketch of Strategy \ref{strategy}. The upper region represents the $k$-polymatroids not in $\widetilde{\mathcal{C}}'_k$ and the lower region represents the the $k$-polymatroids that are in $\widetilde{\mathcal{C}}'_k$. Each node represents a $k$-polymatroid and the directed edges $\rho \to \rho'$ indicate that $\rho'$ is a compression of $\rho$. The black nodes represent the polymatroids in $\Gamma$ because all of their compressions lie in $\widetilde{\mathcal{C}}'_k$.}
    \label{figure: compressions}
\end{figure}

Now consider when $\mathcal{C}$ is the class of binary matroids, so $\widetilde{\mathcal{C}}'_k = \mathcal{P}_{U_{2,4}}^k$. Let $(E,\rho)$ be an excluded minor for $\mathcal{P}_{U_{2,4}}^k$, so $M_\rho^k$ contains a $U_{2,4}$-minor. For any $e \in E$, it must be that deleting the entirety of $X_e$ or contracting the entirety of $X_e$ from $M_\rho^k$ eliminates all $U_{2,4}$-minors of $M_\rho^k$. This is because $M_\rho^k\backslash X_e = M_{\rho\backslash e}^k$ and $M_\rho^k/X_e = M_{\rho/e}^k$, but $M_{\rho\backslash e}^k$ and $M_{\rho/e}^k$ are both binary since $\rho$ being an excluded minor for $\mathcal{P}_{U_{2,4}}^k$ implies $\rho/e$ and $\rho\backslash e$ are in $\mathcal{P}_{U_{2,4}}^k$. It must therefore be that if $\rho(e) = 1$, then at least one clone of $X_e = \{e_1,\hdots, e_k\}$ is in each $U_{2,4}$-minor of $M_\rho^k$, and if $\rho(e)\geq 2$, then to get a $U_{2,4}$-minor of $M_\rho^k$, we must do exactly one of the following:
\begin{enumerate}
    \item Contract at least one (but not all) of the clones of $X_e$ and delete the rest, or
    \item Have at least one of the $e_1, \hdots, e_k$ in the $U_{2,4}$-minor.
\end{enumerate}

Assume (1) applies to $e$, and $Y\subseteq X_e$ is the set of contracted clones, where $|Y| = l$ for some $1 \leq l \leq k-1$. Then the $k$-natural matroid of $\rho_{\downarrow e}^l$ is isomorphic to $M_\rho^k/Y\backslash X_e-Y$, which by assumption contains a $U_{2,4}$-minor. Therefore, it must be that $\rho\notin\Gamma$. If $\rho \in \Gamma$, then (1) cannot apply to any $e \in E$. Since (2) can occur for at most four elements $e$, it must be that $|E| \leq 4$. Therefore, to complete step (1) in our strategy for $\mathcal{P}_{U_{2,4}}^k$, we only need to iterate through all of the possibilities for $k$-polymatroids on a ground set of size at most $4$. For the remainder of this paper, we will only consider $k$-polymatroids where $k \geq 3$, and we will refer to \textit{excluded minors for $\mathcal{P}_{U_{2,4}}^k$} as simply \textit{excluded minors}.

\section{Excluded minors on \texorpdfstring{$|E| = 1$}{}}
\label{section:ground set 1}

\noindent
The following theorem tells us that there are $k-3$ excluded singletons for $\mathcal{P}_{U_{2,4}}^k$. The ability to rule out these singletons as elements of $E(\rho)$ for $\rho \in \mathcal{P}_{U_{2,4}}^k$ is crucial to the tractability of the problem as $k$ increases.

\begin{theorem}\label{singletons}
Consider the $k$-polymatroid $(\{e\}, \rho)$ of rank $m$.
\begin{enumerate}
    \item If $m \in \{0, 1, k-1, k\}$, then $\rho \in \mathcal{P}_{U_{2,4}}^k$ and hence is not an excluded minor.
    \item If $2 \leq m \leq k-2$, then $\rho$ is an excluded minor.
\end{enumerate}
\end{theorem}
\begin{proof} For any $k$-polymatroid $(\{e\}, \rho)$ of rank $m$, note that $|E(M_\rho^k)| = k$ and the rank of $M_\rho^k$ is $m$.
\begin{enumerate}
    \item If $m = 0$ or $m = 1$, then $\rho$ can be represented by a $1 \times 1$ zero matrix or identity matrix respectively, so $\rho \in \mathcal{P}_{U_{2,4}}^k$. Now, assume $m = k-1$ or $m = k$. To obtain a $U_{2,4}$-minor of $M_\rho^k$, we must lower the rank by $k-3$ or $k-2$ respectively and we need to decrease the size of the ground set by $k-4$. Since both $k-3$ and $k-2$ are strictly greater than $k-4$, this is not possible. So $\rho \in \mathcal{P}_{U_{2,4}}$.
    \item The only proper minor of $\rho$ is the empty matroid, which is in $\mathcal{P}_{U_{2,4}}$. Next, we show that $M_\rho^k$ contains a $U_{2,4}$-minor. Let $Y\subseteq X_e$ be a set of $m-2$ clones and let $Z\subseteq X_e-Y$ be a set of $k-m-2$ clones. Then $M_\rho^k/Y\backslash Z \cong U_{2,4}$.
\end{enumerate}
\end{proof}

Going forward, by Theorem \ref{singletons} and Lemma \ref{simplification}, if a $k$-polymatroid $(E,\rho)$ on $|E| \geq 2$ is an excluded minor, its singleton ranks must only take values in $\{1, k-1, k\}$. We define the \textit{type} of $\rho$ be 
\begin{displaymath}
T_\rho:= (a_1, a_{k-1},a_k)
\end{displaymath}
where $a_i$ represents the number of elements in $E(\rho)$ of rank $i$.

\section{Excluded minors on \texorpdfstring{$|E| = 2$}{}}
\label{section:ground set 2}

In this section, let $(E,\rho)$ be a $k$-polymatroid on $E = \{e, f\}$. For $\rho$ to be an excluded minor, its type $T_\rho$ must come from Table \ref{figure: types |E| = 2} (see below). Note that $(2, 0, 0)$ is  omitted from the table because $E(S(M_\rho^k))$ must contain at least $4$ elements in order to contain a $U_{2,4}$-minor.

\begin{table}[H]
\centering
\begin{tabular}{ |c|c|c|c|c|c| }
 \hline
\textbf{Name}&$\alpha$&$\beta$&$\gamma$&$\delta$&$\epsilon$ \\ 
\hline
\textbf{Type $T_\rho$}&$(0, 2, 0)$&$(0, 0, 2)$&$(1, 1, 0)$&$(1, 0, 1)$&$(0, 1, 1)$\\
 \hline
\end{tabular}
\caption{The possible types $T_\rho$ for an excluded minor $\rho$ on $|E| = 2$.}
\label{figure: types |E| = 2}
\end{table}

\begin{theorem}\label{excluded minors |E| = 2}
For a fixed $k \geq 4$, the following table lists the excluded minors for $\mathcal{P}_{U_{2,4}}^k$ on $|E| = 2$:
\vspace{-1ex}
\begin{center}
\def\arraystretch{1.5}
\adjustbox{max width=\textwidth}{\begin{tabular}{ |c|c|c|c| } 
 \hline
\textbf{Name(s)} & \textbf{Rank}  & \textbf{Geometric Description} & \textbf{$k$-Dual} \\ 
\hline
$Ex_\alpha^{k-1}$&  $k-1$ & two parallel rank-$(k-1)$ elements & $Ex_\beta^{k+1}$\\
 \hline
 $Ex_\alpha^k$ & $k$ & two rank-$(k-1)$ elements spanning rank $k$ & identically self-$k$-dual\\
 \hline
 $Ex_\beta^k$&$k$& two parallel rank-$k$ elements & identically self-$k$-dual\\
 \hline
 $Ex_\beta^{k+1}$& $k+1$& two rank-$k$ elements spanning rank $k+1$& $Ex_\alpha^{k-1}$\\
 \hline
 $Ex_\epsilon^k$&$k$& a rank-$(k-1)$ element lying on a rank-$k$ element&self-$k$-dual\\
\hline
\end{tabular}}
\end{center}
When $k = 3$, the previous table combined with the following table provide a complete list of excluded minors for $\mathcal{P}_{U_{2,4}}^3$ on $|E| = 2$. (See Figure \ref{fig:kduality example}.)
\begin{center}
\def\arraystretch{1.5}
\adjustbox{max width=\textwidth}{\begin{tabular}{ |c|c|c|c| } 
 \hline
\textbf{Name(s)} &\textbf{Rank}  & \textbf{Geometric Description} & \textbf{$3$-Dual} \\ 
\hline
 $Ex_\gamma^2$&$2$&a point lying on a line&$Ex_\epsilon^4$\\
\hline
 $Ex_\epsilon^4$& $4$&a line and a plane spanning rank $4$&$Ex_\gamma^2$\\
 \hline
\end{tabular}}
\end{center}
\end{theorem}

\begin{proof}
The justification for Theorem \ref{excluded minors |E| = 2} follows from the remaining lemmas in this section which individually address each type $T_\rho$ from Table \ref{figure: types |E| = 2}. For each lemma, it is assumed that the total rank $\rho(E)$ must be at least the rank of the largest singleton of $\rho$ and at most the sum of the singleton ranks of $\rho$. Note: we will not explicitly derive the $k$-duals of the excluded minors as they are straightforward to check.

\begin{lemma}[Type $\alpha$]\label{alpha}
Let $T_\rho = \alpha = (0, 2, 0)$.
\begin{enumerate}
    \item If $\rho(E) = k-1$ or $k$, then $\rho$ is an excluded minor. We will denote $\rho$ as $Ex_\alpha^{k-1}$ or $Ex_\alpha^k$ respectively.
    \item If $k+1 \leq \rho(E) \leq 2k-3$, then $\rho\notin \mathcal{P}_{U_{2,4}}^k$ but $\rho$ is not an excluded minor. (Note: When $k = 3$, $k+1 > 2k-3$ so there are no polymatroids to consider for $k = 3$.)
    \item If $\rho(E) = 2k-2$, then $\rho\in \mathcal{P}_{U_{2,4}}^k$.
\end{enumerate}
\end{lemma}
\begin{proof}
    We have $\rho(e) = \rho(f) = k-1$.
    \begin{enumerate}
        \item \begin{enumerate}
            \item Assume $\rho(E) = k-1$. Then $\rho$ must consist of two parallel rank-$(k-1)$ elements. $M_\rho^k$ consists of $2k$ points lying freely in rank $k-1$. Let $Y\subseteq X_e$ be a set of $k-3$ clones. Then $M_\rho^k/Y$ contains a $U_{2,4}$-restriction, so $\rho \notin \mathcal{P}_{U_{2,4}}^k$. It remains to show that all nonempty proper polymatroid minors of $\rho$ are in $\mathcal{P}_{U_{2,4}}^k$. Without loss of generality, it suffices to show that $\rho/e \in \mathcal{P}_{U_{2,4}}^k$ and $\rho\backslash e \in \mathcal{P}_{U_{2,4}}^k$. The rank of $\rho/e$ is $0$, and $\rho\backslash e$ is a singleton of rank $k-1$, so both are in $\mathcal{P}_{U_{2,4}}^k$ by Theorem \ref{singletons}.
            \item Assume $\rho(E) = k$. Then $\rho$ must consist of two rank-$(k-1)$ elements spanning rank $k$. Let $Y$ be a subset of $k-2$ clones of $X_e$. Then $M_\rho^k/Y$ contains a $U_{2,4}$-restriction so $\rho\notin\mathcal{P}_{U_{2,4}}^k$. It remains to show that all nonempty proper polymatroid minors of $\rho$ are in $\mathcal{P}_{U_{2,4}}^k$. Without loss of generality, it suffices to show that $\rho/e \in \mathcal{P}_{U_{2,4}}^k$ and $\rho\backslash e \in \mathcal{P}_{U_{2,4}}^k$. The rank of $\rho/e$ is $1$ and $\rho\backslash e$ is a singleton of rank $k-1$, so both are in $\mathcal{P}_{U_{2,4}}^k$ by Theorem \ref{singletons}.
        \end{enumerate}
        \item  Assume $k+1 \leq \rho(E) \leq 2k-3$. We have $(\rho/e)(f) = \rho(E)-(k-1)$, so $2 \leq (\rho/e)(f) \leq k-2$. Therefore, $\rho/e$ is an excluded minor by Theorem \ref{singletons}, which implies both that $\rho \notin \mathcal{P}_{U_{2,4}}^k$ and that $\rho$ is not an excluded minor.
        \item Assume $\rho(E) = 2k-2$. This implies $\rho = \rho|_e \oplus \rho|_f$. By Theorem \ref{singletons}, $\rho|_e, \rho|_f \in \mathcal{P}_{U_{2,4}}^k$, so by Lemma \ref{connectedness} we have $\rho \in \mathcal{P}_{U_{2,4}}^k$.
    \end{enumerate}
\end{proof}

\begin{corollary}\label{skew kminus1 elements}
Let $(E,\rho)$ be a $k$-polymatroid on $|E| \geq 3$ which is an excluded minor. Assume $\rho$ contains exactly $l$ elements of rank $k-1$. Then $\rho(E) \geq l(k-1)$.
\end{corollary}
\begin{proof}
Any pair of rank-$(k-1)$ elements $e$ and $f$ must together span rank $2k-2$ in $\rho$. Otherwise, by Lemma \ref{alpha}, $\{e,f\}$ would be a proper restriction of $\rho$ that is not in $\mathcal{P}_{U_{2,4}}^k$ and therefore $\rho$ cannot be an excluded minor of $\mathcal{P}_{U_{2,4}}^k$. We induct on $l$. Let $l = 2$. Then the two rank-$(k-1)$ elements must span rank $2k-2$, so $\rho(E) \geq 2(k-1)$. Now, given $\rho$ with $l+1$ rank-$(k-1)$ elements, fix one, call it $e$. Since $e$ is skew to each of the other $l$ rank-$(k-1)$ elements, contracting $e$ leaves the latter unaffected. By the induction hypothesis, $(\rho/e)(E-e) \geq l(k-1)$, so $\rho(E) \geq l(k-1)+(k-1) = (l+1)(k-1)$.
\end{proof}

\begin{lemma}[Type $\beta$]\label{beta}
Let $T_\rho = \beta = (0, 0, 2)$. 
\begin{enumerate}
    \item If $\rho(E) = k$ or $\rho(E) = k+1$, then $\rho$ is an excluded minor. We will denote $\rho$ as $Ex_\beta^k$ or $Ex_\beta^{k+1}$ respectively.
    \item If $k+2 \leq \rho(E) \leq 2k-2$, then $\rho \notin\mathcal{P}_{U_{2,4}}^k$ but $\rho$ is not an excluded minor. (Note: When $k = 3$, $k+2 > 2k-2$, so there are no polymatroids to consider for $k = 3$.)
    \item If $\rho(E) = 2k-1$ or $\rho(E) = 2k$, then $\rho \in\mathcal{P}_{U_{2,4}}^k$.
\end{enumerate}
\end{lemma}
\begin{proof} We have $\rho(e) = \rho(f) = k$. 
\begin{enumerate}
    \item \begin{enumerate}
        \item Assume $\rho(E) = k$. Then $\rho$ must consist of two parallel rank-$k$ elements. $M_\rho^k$ consists of $2k$ points lying freely in rank $k$. Let $Y$ be a subset of $k-2$ clones of $X_e$. Then $M_\rho^k/Y$ contains a $U_{2,4}$-restriction so $\rho\notin\mathcal{P}_{U_{2,4}}^k$. It remains to show that all nonempty proper polymatroid minors of $\rho$ are in $\mathcal{P}_{U_{2,4}}^k$. Without loss of generality, it suffices to show that $\rho/e \in \mathcal{P}_{U_{2,4}}^k$ and $\rho\backslash e \in \mathcal{P}_{U_{2,4}}^k$. The rank of $\rho/e$ is $0$, and $\rho\backslash e$ is a singleton of rank $k$, so both are in $\mathcal{P}_{U_{2,4}}^k$ by Theorem \ref{singletons}.
        \item Assume $\rho(E) = k+1$. Then $\rho$ must consist of two rank-$k$ elements spanning rank $k+1$. This is the $k$-dual of $Ex_\alpha^{k-1}$, which we have already shown to be an excluded minor by Lemma \ref{alpha}. Therefore, by Lemma \ref{kduality}, $\rho$ is also an excluded minor.
    \end{enumerate}
    \item Assume $k+2 \leq \rho(E) \leq 2k-2$. We have $(\rho/e)(f) = \rho(E)-k$, so $2 \leq (\rho/e)(f) \leq k-2$. Therefore, $\rho/e$ is an excluded minor by Theorem \ref{singletons}, which implies both that $\rho\notin \mathcal{P}_{U_{2,4}}^k$ and that $\rho$ is not an excluded minor.
    \item \begin{enumerate}
        \item Assume $\rho(E) = 2k-1$. Note that $|E(M_\rho^k)| = 2k$ and that the total rank of $M_\rho^k$ is $2k-1$. To get a $U_{2,4}$-minor of $M_\rho^k$, we must decrease the size of the ground set by $2k-4$ and decrease the rank by $2k-3$. Since $2k-4 < 2k-3$, this is not possible. Therefore, $\rho \in \mathcal{P}_{U_{2,4}}^k$.
        \item Assume $\rho(E) = 2k$. This implies $\rho = \rho|_e \oplus \rho|_f$. By Theorem \ref{singletons}, $\rho|_e, \rho|_f \in \mathcal{P}_{U_{2,4}}^k$, so by Lemma \ref{connectedness}, we have $\rho \in \mathcal{P}_{U_{2,4}}^k$.
    \end{enumerate}
\end{enumerate}
\end{proof}

\begin{lemma}[Type $\gamma$]\label{gamma} Let $T_\rho = \gamma = (1, 1, 0)$. 
\begin{enumerate}
    \item Let $\rho(E) = k-1$. Then $\rho\notin\mathcal{P}_{U_{2,4}}^k$. 
    \begin{enumerate}
        \item If $k = 3$, then $\rho$ is an excluded minor for $\mathcal{P}_{U_{2,4}}^3$. We will denote $\rho$ as $Ex_\gamma^2$.
        \item If $k\geq 4$, then $\rho$ is not an excluded minor.
    \end{enumerate}
    \item If $\rho(E) = k$, then $\rho\in\mathcal{P}_{U_{2,4}}^k$.
\end{enumerate}
\end{lemma}
\begin{proof}   Let $\rho(e) = 1$ and $\rho(f) = k-1$.
\begin{enumerate}
    \item Assume $\rho(E) = k-1$.
    \begin{enumerate}
        \item Let $k = 3$. We have $\rho(e) = 1$ and $\rho(f) = \rho(E) = 2$. Geometrically, $\rho$ consists of a point lying on a line. We have $S(M_{\rho}^3) \cong U_{2,4}$, so $\rho\notin \mathcal{P}_{U_{2,4}}^3$. It remains to show that all nonempty proper polymatroid minors of $\rho$ are in $\mathcal{P}_{U_{2,4}}^3$. Indeed, if we delete or contract $e$ or $f$ from $\rho$ to get a polymatroid minor $\rho'$, $M_{\rho'}^3$ will only contain $3$ points, and thus cannot contain a $U_{2,4}$-minor, so $\rho' \in \mathcal{P}_{U_{2,4}}^3$.
        \item Let $k \geq 4$. We have $(\rho/e)(f) = k-2$, which implies $\rho/e$ is an excluded minor by Theorem \ref{singletons}, so $\rho \notin \mathcal{P}_{U_{2,4}}^k$ and $\rho$ is not an excluded minor.
    \end{enumerate}
    \item Assume $\rho(E) = k$. This implies $\rho = \rho|_e \oplus \rho|_f$. By Theorem \ref{singletons}, $\rho|_e, \rho|_f \in \mathcal{P}_{U_{2,4}}^k$, so by Lemma \ref{connectedness}, we have $\rho \in \mathcal{P}_{U_{2,4}}^k$.
\end{enumerate}
\end{proof}

\begin{lemma}[Type $\delta$]\label{delta}
If $T_\rho = \delta = (1, 0, 1)$, then $\rho \in \mathcal{P}_{U_{2,4}}^k$.
\end{lemma}
\begin{proof}
    Let $\rho(e) = 1$ and $\rho(f) = k$. First, assume $\rho(E) = k$. Then $M_\rho^k$ consists of $k+1$ points freely placed in rank $k$. To obtain a $U_{2,4}$-minor of $M_\rho^k$, we must decrease the size of the ground set by $k-3$ and decrease the rank by $k-2$. Since $k-3 < k-2$, this is not possible. Therefore, $\rho \in \mathcal{P}_{U_{2,4}}^k$.

    Next, assume $\rho(E) = k+1$. This implies $\rho = \rho|_e \oplus \rho|_f$. By Theorem \ref{singletons}, $\rho|_e, \rho|_f \in \mathcal{P}_{U_{2,4}}^k$, so by Lemma \ref{connectedness}, we have $\rho \in \mathcal{P}_{U_{2,4}}^k$.
\end{proof}

\begin{lemma}[Type $\epsilon$]\label{epsilon}
Let $T_\rho = \epsilon = (0, 1, 1)$.
\begin{enumerate}
    \item If $\rho(E) = k$, then $\rho$ is an excluded minor. We will denote $\rho$ as $Ex_\epsilon^k$.
    \item Let $k+1 \leq \rho(E) \leq 2k-2$. Then $\rho\notin \mathcal{P}_{U_{2,4}}^k$. 
    \begin{enumerate}
        \item If $k = 3$, then $\rho$ is an excluded minor for $\mathcal{P}_{U_{2,4}}^3$. We will denote $\rho$ as $Ex_\epsilon^4$.
        \item If $k \geq 4$, then $\rho$ is not an excluded minor.
    \end{enumerate}
    \item If $\rho(E) = 2k-1$, then $\rho\in \mathcal{P}_{U_{2,4}}^k$.
\end{enumerate}
\end{lemma}
\begin{proof} Let $\rho(e) = k-1$ and $\rho(f) = k$.
\begin{enumerate}
    \item Assume $\rho(E) = k$. Then $\rho$ must consists of $e$ lying on $f$. Let $Y$ be a subset of $k-2$ clones of $X_e$. Then $M_\rho^k/Y$ contains a $U_{2,4}$-restriction. It remains to show that all nonempty proper polymatroid minors of $\rho$ are in $\mathcal{P}_{U_{2,4}}^k$. Indeed, the minors $\rho/f$, $\rho/e$, $\rho\backslash f$, and $\rho\backslash e$ are singletons of ranks $0, 1, k$, and $k-1$ respectively, so by Theorem \ref{singletons}, they are all in $\mathcal{P}_{U_{2,4}}^k$.
    \item Assume $k+1 \leq \rho(E) \leq 2k-2$.
    \begin{enumerate}
        \item Let $k = 3$. Then $\rho$ is the $k$-dual of $Ex_\gamma^2$, which we have already shown in Lemma \ref{gamma} to be an excluded minor. Thus, $\rho$ is also an excluded minor by Lemma \ref{kduality}.
        \item Let $k \geq 4$. First, let $k+1 \leq \rho(E) \leq 2k-3$. We have $(\rho/e)(f) = \rho(E)-(k-1)$, so $2 \leq (\rho/e)(f) \leq k-2$. Therefore, by Theorem \ref{singletons}, $(\rho/e)\notin \mathcal{P}_{U_{2,4}}^k$, which implies both that $\rho\notin \mathcal{P}_{U_{2,4}}^k$ and that $\rho$ is not an excluded minor. Finally, if $\rho(E) = 2k-2$, then $(\rho/f)(e) = k-2$, so the same argument holds as above for $\rho/f$.
    \end{enumerate}
    \item Assume $\rho(E) = 2k-1$. This implies $\rho = \rho|_e \oplus \rho|_f$. By Theorem \ref{singletons}, $\rho|_e, \rho|_f \in \mathcal{P}_{U_{2,4}}^k$, so by Lemma \ref{connectedness}, we have $\rho \in \mathcal{P}_{U_{2,4}}^k$.
\end{enumerate}
\end{proof}

This concludes the proof of Theorem \ref{excluded minors |E| = 2}.
\end{proof}

\section{Excluded minors on \texorpdfstring{$|E| = 3$}{}}
\label{section:ground set 3}

In this section, $(E,\rho)$ will always be a $k$-polymatroid on $E = \{e, f, g\}$.

\begin{theorem}\label{noE3} There are no excluded minors for $\mathcal{P}_{U_{2,4}}^k$ on a ground set of size $3$.
\end{theorem}

\begin{proof}
The justification for Theorem \ref{noE3} follows from the lemmas in this section addressing the corresponding types $T_\rho$ shown in Table \ref{figure: types |E| = 3}. As in the previous section, it will be assumed that the total rank $\rho(E)$ must be at least the rank of the largest singleton of $\rho$ and at most the sum of the singleton ranks of $\rho$. By Lemma \ref{kduality}, since $\rho^*(E) + \rho(E) = 3k$, we only need to show that there are no excluded minors $\rho$ on $|E| = 3$ for which $\rho(E) \leq \lfloor \frac{3k}{2} \rfloor$. Note that $T_\rho = (3, 0, 0)$ is omitted from the table because $E(S(M_\rho^k))$ must contain at least $4$ elements in order to contain a $U_{2,4}$-minor.

\begin{table}[H]
\centering
\footnotesize
\begin{tabular}{ |c|c|c|c|c|c|c| } 
 \hline
\textbf{Lemma} & \ref{030E3 021E3} & \ref{120E3}& \ref{111E3}&\ref{003E3 102E3 012E3}&\ref{210E3}& \ref{201E3}\\
\hline
\textbf{Types $T_\rho$}& $(0, 3, 0)$& $(1, 2, 0)$& \makecell{$(1, 1, 1)$\\$(0, 2, 1)$}& \makecell{$(0, 0, 3)$\\$(1, 0, 2)$\\$(0, 1, 2)$}& $(2, 1, 0)$& $(2, 0, 1)$\\
\hline
\end{tabular}
\caption{The possible types $T_\rho$ for an excluded minor $\rho$ on $|E| = 3$ and the corresponding lemmas in which they are addressed.}
\label{figure: types |E| = 3}
\end{table}

\begin{lemma}\label{030E3 021E3} There are no excluded minors $(E,\rho)$ of types $(0, 3, 0)$ or $(0, 2, 1)$ such that $\rho(E) \leq \lfloor \frac{3k}{2} \rfloor$.
\end{lemma}
\begin{proof}
    Let $\rho$ be an excluded minor of type $T_\rho = (0, 3, 0)$. By Corollary \ref{skew kminus1 elements}, $\rho(E) \geq 3k-3$. But we are requiring $\rho(E) \leq \lfloor \frac{3k}{2} \rfloor < 3k-3$, so the desired conclusion follows.
\end{proof}

\begin{lemma}\label{120E3}
Let $(E,\rho)$ be simple and of type $T_\rho = (1, 2, 0)$. 
\begin{enumerate}
    \item If $k-1 \leq \rho(E) \leq 2k-2$, then $\rho\notin\mathcal{P}_{U_{2,4}}^k$ and $\rho$ is not an excluded minor.
    \item If $\rho(E) = 2k-1$, then $\rho\in\mathcal{P}_{U_{2,4}}^k$ so $\rho$ cannot be an excluded minor.
\end{enumerate}
\end{lemma}
\begin{proof} Let $\rho(e) = 1$, and $\rho(f) = \rho(g) = k-1$. 
\begin{enumerate}
    \item 
    \begin{enumerate}
        \item Assume $k-1 \leq \rho(E) \leq 2k-3$. By Lemma \ref{alpha}, $(\rho\backslash e)\notin\mathcal{P}_{U_{2,4}}^k$, which implies $\rho\notin\mathcal{P}_{U_{2,4}}^k$ and that $\rho$ is not an excluded minor. 
        \item Assume $\rho(E) = 2k-2$ and that $\rho$ is an excluded minor. By the lemmas in Section \ref{section:ground set 2}, it must be that $\rho(ef) = \rho(eg) = k$, and $\rho(fg) = 2k-2$. Then $(\rho/e) \notin \mathcal{P}_{U_{2,4}}^k$ by Lemma \ref{alpha}, so $\rho$ cannot be an excluded minor.
    \end{enumerate}
    \item Assume $\rho(E) = 2k-1$. Then $\rho \cong \rho|_e \oplus \rho|_f \oplus \rho|_g$. By Theorem \ref{singletons}, $\rho|_e, \rho|_f, \rho|_g \in \mathcal{P}_{U_{2,4}}^k$, so $\rho \in \mathcal{P}_{U_{2,4}}^k$ by Lemma \ref{connectedness} and thus $\rho$ cannot be an excluded minor.
\end{enumerate}
\end{proof}

\begin{lemma}\label{111E3}
There are no excluded minors $(E,\rho)$ of types $(1, 1, 1)$ or $(0, 2, 1)$ such that $\rho(E) \leq \lfloor \frac{3k}{2} \rfloor$.
\end{lemma}
\begin{proof}
    Let $\rho$ be an excluded minor of type $(1, 1, 1)$ or $(0, 2, 1)$. Let $\rho(f) = k-1$ and $\rho(g) = k$. By Lemma \ref{epsilon}, $\rho(E) \geq \rho(fg) \geq 2k-1$. Since we are requiring that $\rho(E) \leq \lfloor \frac{3k}{2}\rfloor< 2k-1$, there is no such $\rho$.
\end{proof}

\begin{lemma}\label{003E3 102E3 012E3} \leavevmode
\begin{enumerate}
    \item There are no excluded minors $(E, \rho)$ of types $(0, 0, 3), (1, 0, 2)$, or $(0, 1, 2)$ such that $\rho(E) \leq \lfloor \frac{3k}{2} \rfloor$.
    \item If $\rho$ is of type $(0, 0, 3)$ and $\rho\in \mathcal{P}_{U_{2,4}}^k$, then $3k-2\leq \rho(E) \leq 3k$.
\end{enumerate}
\end{lemma}
\begin{proof}\leavevmode
\begin{enumerate}
    \item Let $\rho$ be an excluded minor of type $(0, 0, 3), (1, 0, 2)$, or $(0, 1, 2)$. Let $\rho(f)= \rho(g) = k$. By Lemma \ref{beta}, $\rho(E) \geq \rho(fg) \geq 2k-1$. Since we are requiring that $\rho(E) \leq \lfloor \frac{3k}{2} \rfloor< 2k-1$, there is no such $\rho$.
    \item Let $\rho(e) = \rho(f) = \rho(g) = k$. It must be that $\rho(E) \leq \rho(e) + \rho(f) + \rho(g) = 3k$. By Lemma \ref{beta}, $\rho(ef), \rho(eg), \rho(fg) \in \{2k-1, 2k\}$. Then, it must be that $\rho(E) \geq 3k-2$; otherwise by Lemmas \ref{alpha}, \ref{beta}, and \ref{epsilon}, $(\rho/e) \notin \mathcal{P}_{U_{2,4}}^k$.
\end{enumerate}
\end{proof}

\begin{lemma}\label{210E3} Let $(E,\rho)$ be simple and of type $(2, 1, 0)$.
\begin{enumerate}
    \item If $\rho(E) = k-1$ or $k$, then $\rho\notin \mathcal{P}_{U_{2,4}}$ and $\rho$ is not an excluded minor.
    \item If $\rho(E) = k+1$, then $\rho \in \mathcal{P}_{U_{2,4}}^k$ so $\rho$ cannot be an excluded minor.
\end{enumerate}
\end{lemma}
\begin{proof}
    Let $\rho(e) = \rho(f) = 1$, and $\rho(g) = k-1$.
    \begin{enumerate}
        \item \begin{enumerate}
            \item If $\rho(E) = k-1$, then $(\rho\backslash f) \notin \mathcal{P}_{U_{2,4}}^k$ by Lemma \ref{gamma}, which implies $\rho \notin \mathcal{P}_{U_{2,4}}^k$ and that $\rho$ is not an excluded minor.
            \item Assume $\rho(E) = k$ and $\rho$ is an excluded minor. By Lemmas \ref{simplification} and \ref{gamma}, $\rho(ef) = 2$ and $\rho(eg) = \rho(fg) = k$. Observe that $(\rho/e) \notin \mathcal{P}_{U_{2,4}}^k$ by Lemma \ref{gamma}. This implies $\rho\notin \mathcal{P}_{U_{2,4}}^k$ and that $\rho$ is not an excluded minor.
        \end{enumerate}
        \item Let $\rho(E) = k+1$. Then $\rho \cong \rho|_e \oplus \rho|_f \oplus \rho|_g$. By Theorem \ref{singletons}, $\rho|_e, \rho|_f, \rho|_g \in \mathcal{P}_{U_{2,4}}^k$, so $\rho \in \mathcal{P}_{U_{2,4}}^k$ by Lemma \ref{connectedness} and thus $\rho$ cannot be an excluded minor.
    \end{enumerate}
\end{proof}

\begin{lemma}\label{201E3} Let $(E,\rho)$ be simple and of type $(2, 0, 1)$.
\begin{enumerate}
    \item If $\rho(E) = k$, then $\rho\notin \mathcal{P}_{U_{2,4}}$ and $\rho$ is not an excluded minor.
    \item If $\rho(E) = k+1$ or $k+2$, then $\rho \in \mathcal{P}_{U_{2,4}}^k$ so $\rho$ cannot be an excluded minor.
\end{enumerate}
\end{lemma}
\begin{proof} Let $\rho(e) = \rho(f) = 1$ and $\rho(g) = k$.
\begin{enumerate}
    \item Assume $\rho(E) = k$ and that $\rho$ is an excluded minor. We know that $\rho(eg) = \rho(fg) = k$ because $\rho(g) = k$. We also know $\rho(ef) = 2$ by Lemma \ref{simplification}. Then $(\rho/e) \notin \mathcal{P}_{U_{2,4}}^k$ by Lemma \ref{gamma}. Thus $\rho \notin \mathcal{P}_{U_{2,4}}^k$ and $\rho$ is not an excluded minor.
    \item \begin{enumerate}
        \item If $\rho(E) = k+1$, then $S(M_\rho^k)$ would have $k+2$ points in rank $k+1$. To get a $U_{2,4}$-minor of $S(M_\rho^k)$, we must decrease the size of its ground set by $k-2$ and decrease the rank by $k-1$. Since $k-2 < k-1$, this is not possible. It must be that $\rho \in \mathcal{P}_{U_{2,4}}^k$, so $\rho$ is not an excluded minor.
        \item If $\rho(E) = k+2$, then $\rho \cong \rho|_e \oplus \rho|_f \oplus \rho|_g$. By Theorem \ref{singletons}, $\rho|_e, \rho|_f, \rho|_g \in \mathcal{P}_{U_{2,4}}^k$, so it must be that $\rho \in \mathcal{P}_{U_{2,4}}^k$ by Lemma \ref{connectedness} and thus $\rho$ cannot be an excluded minor.
    \end{enumerate}
\end{enumerate} 
\end{proof}
This concludes the proof of Theorem \ref{noE3}.
\end{proof}

\section{Excluded minors with \texorpdfstring{$|E| = 4$}{}}
\label{section:ground set 4}

In this section, $(E,\rho)$ will always be a $k$-polymatroid on $E = \{e, f, g, h\}$. 
We will frequently begin our analysis by assuming that $\rho$ is an excluded minor; from there we fill in as many values of $\rho$ as possible by repeatedly applying the following preliminary criteria (1) through (9) without mention:
\begin{enumerate}
    \item Since $\rho$ is a polymatroid, $\rho$ must be monotone and submodular.
    \item Since $\rho$ is connected, we require $\rho(E) < \rho(e) + \rho(f) + \rho(g) + \rho(h)$.
    \item By Lemma \ref{simplification}, $\rho$ must be simple. In particular, if $e$ and $f$ are points, then $\rho(ef) = 2$.
\end{enumerate}
For any $e, f, g \in E$,
\begin{enumerate}\addtocounter{enumi}{3}
    \item If $\rho(e) = \rho(f) = k-1$, then $\rho(ef) = 2k-2$ by Lemma \ref{alpha}.
    \item If $\rho(e) = \rho(f) = k$, then $\rho(ef) = 2k-1$ or $2k$ by Lemma \ref{beta}.
    \item If $\rho(e) = 1$ and $\rho(f) = k-1$, then $\rho(ef) = k$ by Lemma \ref{gamma}.
    \item If $\rho(e) = k-1$ and $\rho(f) = k$, then $\rho(ef) = 2k-1$ by Lemma \ref{epsilon}.
    \item If $T_{(\rho|_{efg})} = (2, 1, 0)$, then $\rho(efg) = k+1$ by Lemma \ref{210E3}.
    \item If $T_{(\rho|_{efg})} = (2, 0, 1)$, then $\rho(efg)$ cannot equal $k$ by Lemma \ref{201E3}. That is, two points cannot lie in a rank-$k$ element.
\end{enumerate}

Note: If $\rho(e) = 1$ for any $e \in E$, then $X_e$ is a parallel class in $M_\rho^k$ consisting of $k$ clones, so in $S(M_\rho^k)$, exactly one of these clones survives. We will refer to this surviving clone also as $e$.

\begin{theorem}\label{excluded minors |E| = 4}
The following is the list of excluded minors for $\mathcal{P}_{U_{2,4}}^k$ on $|E| = 4$: 
\begin{itemize}
    \item $Ex_{(4, 0, 0)}^2$, also known as $U_{2,4}$, consisting of four collinear points. Its $k$-dual is $Ex_{(0, 0, 4)}^{4k-2}$.
    \item $Ex_{(3, 0, 1)}^{k+1}$, with rank shown below. Its $k$-dual is $Ex_{(1, 0, 3)}^{3k-1}$.
    \begin{displaymath}
        \begin{tikzcd}[row sep = 10]&&&\underset{\color{blue}E\color{black}}{ k+1}\arrow[d,no head]\arrow[rrd,no head,dotted]&&\\&\underset{\color{blue}efg\color{black}}{ 2}\arrow[rru,no head]\arrow[d,no head]\arrow[rrd,no head,dotted]&\underset{\color{blue}efh\color{black}}{ k+1}\arrow[ru,no head]\arrow[d,no head]\arrow[rrd,no head,dotted]&\underset{\color{blue}egh\color{black}}{ k+1}\arrow[ld,no head]\arrow[rrd,no head,dotted]&&\underset{\color{blue}fgh\color{black}}{ k+1}\arrow[ld,no head]\arrow[d,no head]\\\underset{\color{blue}ef\color{black}}{ 2}\arrow[ru,no head]\arrow[rru,no head]&\underset{\color{blue}eg\color{black}}{ 2}\arrow[rru,no head]&\underset{\color{blue}eh\color{black}}{ k+1}\arrow[rrd,no head,dotted]&\underset{\color{blue}fg\color{black}}{ 2}\arrow[rru,no head]&\underset{\color{blue}fh\color{black}}{ k+1}\arrow[d,no head]&\underset{\color{blue}gh\color{black}}{ k+1}\\\underset{\color{blue}e\color{black}}{ 1}\arrow[u,no head]\arrow[rru,no head]\arrow[ru,no head]&&\underset{\color{blue}f\color{black}}{ 1}\arrow[ru,no head]\arrow[rru,no head]\arrow[llu,no head,dotted]&\underset{\color{blue}g\color{black}}{ 1}\arrow[u,no head]\arrow[rru,no head]\arrow[llu,no head,dotted]&\underset{\color{blue}h\color{black}}{ k}\arrow[ru,no head]&\\&&\underset{\color{blue}\emptyset\color{black}}{0}\arrow[u,no head]\arrow[ru,no head]\arrow[rru,no head]\arrow[llu,no head,dotted]&&&\end{tikzcd}
    \end{displaymath}
    \item $Ex_{(2, 0, 2)}^{2k}$, with rank shown below. It is self-$k$-dual. 
    \begin{displaymath}
        \begin{tikzcd}[row sep = 10]&&&\underset{\color{blue}E\color{black}}{ 2k}\arrow[d,no head]\arrow[rrd,no head,dotted]&&\\&\underset{\color{blue}efg\color{black}}{ k+1}\arrow[rru,no head]\arrow[d,no head]\arrow[rrd,no head,dotted]&\underset{\color{blue}efh\color{black}}{ k+1}\arrow[ru,no head]\arrow[d,no head]\arrow[rrd,no head,dotted]&\underset{\color{blue}egh\color{black}}{ 2k}\arrow[ld,no head]\arrow[rrd,no head,dotted]&&\underset{\color{blue}fgh\color{black}}{ 2k}\arrow[ld,no head]\arrow[d,no head]\\\underset{\color{blue}ef\color{black}}{ 2}\arrow[ru,no head]\arrow[rru,no head]&\underset{\color{blue}eg\color{black}}{ k+1}\arrow[rru,no head]&\underset{\color{blue}eh\color{black}}{ k+1}\arrow[rrd,no head,dotted]&\underset{\color{blue}fg\color{black}}{ k+1}\arrow[rru,no head]&\underset{\color{blue}fh\color{black}}{ k+1}\arrow[d,no head]&\underset{\color{blue}gh\color{black}}{ 2k}\\\underset{\color{blue}e\color{black}}{ 1}\arrow[u,no head]\arrow[rru,no head]\arrow[ru,no head]&&\underset{\color{blue}f\color{black}}{ 1}\arrow[ru,no head]\arrow[rru,no head]\arrow[llu,no head,dotted]&\underset{\color{blue}g\color{black}}{ k}\arrow[u,no head]\arrow[rru,no head]\arrow[llu,no head,dotted]&\underset{\color{blue}h\color{black}}{ k}\arrow[ru,no head]&\\&&\underset{\color{blue}\emptyset\color{black}}{0}\arrow[u,no head]\arrow[ru,no head]\arrow[rru,no head]\arrow[llu,no head,dotted]&&&\end{tikzcd}
    \end{displaymath}
    \item $Ex_{(1, 0, 3)}^{3k-1}$, with rank shown below. Its $k$-dual is $Ex_{(3, 0, 1)}^{k+1}$.
    \begin{displaymath}
        \begin{tikzcd}[row sep = 10]&&&\underset{\color{blue}E\color{black}}{ 3k-1}\arrow[d,no head]\arrow[rrd,no head,dotted]&&\\&\underset{\color{blue}efg\color{black}}{ 2k}\arrow[rru,no head]\arrow[d,no head]\arrow[rrd,no head,dotted]&\underset{\color{blue}efh\color{black}}{ 2k}\arrow[ru,no head]\arrow[d,no head]\arrow[rrd,no head,dotted]&\underset{\color{blue}egh\color{black}}{ 2k}\arrow[ld,no head]\arrow[rrd,no head,dotted]&&\underset{\color{blue}fgh\color{black}}{ 3k-1}\arrow[ld,no head]\arrow[d,no head]\\\underset{\color{blue}ef\color{black}}{ k+1}\arrow[ru,no head]\arrow[rru,no head]&\underset{\color{blue}eg\color{black}}{ k+1}\arrow[rru,no head]&\underset{\color{blue}eh\color{black}}{ k+1}\arrow[rrd,no head,dotted]&\underset{\color{blue}fg\color{black}}{ 2k}\arrow[rru,no head]&\underset{\color{blue}fh\color{black}}{ 2k}\arrow[d,no head]&\underset{\color{blue}gh\color{black}}{ 2k}\\\underset{\color{blue}e\color{black}}{ 1}\arrow[u,no head]\arrow[rru,no head]\arrow[ru,no head]&&\underset{\color{blue}f\color{black}}{ k}\arrow[ru,no head]\arrow[rru,no head]\arrow[llu,no head,dotted]&\underset{\color{blue}g\color{black}}{ k}\arrow[u,no head]\arrow[rru,no head]\arrow[llu,no head,dotted]&\underset{\color{blue}h\color{black}}{ k}\arrow[ru,no head]&\\&&\underset{\color{blue}\emptyset\color{black}}{0}\arrow[u,no head]\arrow[ru,no head]\arrow[rru,no head]\arrow[llu,no head,dotted]&&&\end{tikzcd}
    \end{displaymath}
    \item $Ex_{(0, 0, 4)}^{4k-2}$, with rank shown below. Its $k$-dual is $U_{2,4}$.
    \begin{displaymath}
        \begin{tikzcd}[row sep = 10]&&&\underset{\color{blue}E\color{black}}{ 4k-2}\arrow[d,no head]\arrow[rrd,no head,dotted]&&\\&\underset{\color{blue}efg\color{black}}{ 3k-1}\arrow[rru,no head]\arrow[d,no head]\arrow[rrd,no head,dotted]&\underset{\color{blue}efh\color{black}}{ 3k-1}\arrow[ru,no head]\arrow[d,no head]\arrow[rrd,no head,dotted]&\underset{\color{blue}egh\color{black}}{ 3k-1}\arrow[ld,no head]\arrow[rrd,no head,dotted]&&\underset{\color{blue}fgh\color{black}}{ 3k-1}\arrow[ld,no head]\arrow[d,no head]\\\underset{\color{blue}ef\color{black}}{ 2k}\arrow[ru,no head]\arrow[rru,no head]&\underset{\color{blue}eg\color{black}}{ 2k}\arrow[rru,no head]&\underset{\color{blue}eh\color{black}}{2k}\arrow[rrd,no head,dotted]&\underset{\color{blue}fg\color{black}}{ 2k}\arrow[rru,no head]&\underset{\color{blue}fh\color{black}}{ 2k}\arrow[d,no head]&\underset{\color{blue}gh\color{black}}{ 2k}\\\underset{\color{blue}e\color{black}}{ k}\arrow[u,no head]\arrow[rru,no head]\arrow[ru,no head]&&\underset{\color{blue}f\color{black}}{ k}\arrow[ru,no head]\arrow[rru,no head]\arrow[llu,no head,dotted]&\underset{\color{blue}g\color{black}}{ k}\arrow[u,no head]\arrow[rru,no head]\arrow[llu,no head,dotted]&\underset{\color{blue}h\color{black}}{ k}\arrow[ru,no head]&\\&&\underset{\color{blue}\emptyset\color{black}}{0}\arrow[u,no head]\arrow[ru,no head]\arrow[rru,no head]\arrow[llu,no head,dotted]&&&\end{tikzcd}
    \end{displaymath}
\end{itemize}
\end{theorem}

The justification for Theorem \ref{excluded minors |E| = 4} follows from the lemmas in this section addressing the corresponding types $T_\rho$ shown in Table \ref{figure: types |E| = 4}. 

\begin{table}[H]
\begin{center}
\def\arraystretch{1.5}
\adjustbox{max width=\textwidth}{\begin{tabular}{|c|c|c|c|c|c|c|c|c|c|} 
 \hline
Lemma \ref{U24} & \ref{301E4} & \ref{202E4} &\ref{040E4 031E4 130E4}&\ref{310E4}& \ref{220E4} & \ref{004E4 103E4 013E4}&\ref{022E4 112E4}& \ref{121E4} & \ref{211E4}\\
\hline
$T_\rho = (4, 0, 0)$& $(3, 0, 1)$& $(2, 0, 2)$& \makecell{$(0, 4, 0)$\\$(0, 3, 1)$\\$(1, 3, 0)$}& $(3, 1, 0)$& $(2, 2, 0)$&\makecell{$(0, 0, 4)$\\$(1, 0, 3)$\\$(0, 1, 3)$}&\makecell{$(0, 2, 2)$\\$(1, 1, 2)$}& $(1, 2, 1)$&$(2, 1, 1)$\\
\hline
\end{tabular}}
\end{center}
\caption{The possible types $T_\rho$ for an excluded minor $\rho$ on $|E| = 4$ and the corresponding lemmas in which they are addressed.}
\label{figure: types |E| = 4}
\end{table}

By Lemma \ref{kduality}, since $\rho^*(E) + \rho(E) = 4k$, we only need to consider $\rho$ such that $\rho(E) \leq 2k$.

\vspace{2ex}

    \subsection{Type \texorpdfstring{$(4, 0, 0)$}{}}~\\
        
\begin{lemma}\label{U24}
The excluded minor of type $T_\rho = (4, 0, 0)$ is $U_{2,4}$. \end{lemma}
\begin{proof}
    $U_{2,4}$ is obviously an excluded minor of type $T_\rho = (4, 0, 0)$. If the four points of $\rho$ were in any other simple configuration, then $\rho(E) \geq 3$. To get a minor of $S(M_\rho^k)$ which is isomorphic to $U_{2,4}$, we would have to contract or delete at least one of its points in order to lower the rank to $2$, but then $|E(S(M_\rho^k))| \leq 3 < |E(U_{2,4})|$. Therefore, $U_{2,4}$, also denoted $Ex_{(4, 0, 0)}^2$, is the only excluded minor of type $(4, 0, 0)$ on $|E| = 4$.
\end{proof}

\begin{corollary}\label{U24 kdual}
    The $k$-dual of $U_{2,4}$, denoted $Ex_{(0, 0, 4)}^{4k-2}$, is an excluded minor. \end{corollary}
\vspace{2ex}

    \subsection{Type \texorpdfstring{$(3, 0, 1)$}{}}~\\
        
\begin{lemma}\label{301E4}
    The polymatroid consisting of three collinear points and a rank-$k$ element collectively spanning rank $k+1$, with no point lying on the rank-$k$ element, is the only excluded minor of type $T_\rho = (3, 0, 1)$. We denote it as $Ex_{(3, 0, 1)}^{k+1}$.
\end{lemma}
\begin{proof}
    Let $(E,\rho)$ be an excluded minor of type $T_\rho = (3, 0, 1)$. Let $\rho(e) = \rho(f) = \rho(g) = 1$ and $\rho(h) = k$. Because $\rho(h) = k$ and $\rho$ must be connected, we have $k \leq \rho(E) \leq k+2$. It must also be that $\rho$ is simple. Thus, $S(M_\rho^k)$ has $k+3$ points.
    
    It must be that $\rho(efh) = \rho(egh) = \rho(fgh) = k+1$; to see this, without loss of generality, let $\rho(fgh) = k$, then $\rho|_{fgh} \notin \mathcal{P}_{U_{2,4}}^k$ by Lemma \ref{201E3}. This also forces $\rho(E) \geq k+1$ by monotonicity. If $\rho(E)\geq k+2$, then to get a $U_{2,4}$-minor of $S(M_\rho^k)$, we would have to decrease the size of the ground set by $k-1$ while decreasing the rank by at least $k$. Since $k-1 < k$, this is not possible. Therefore, $\rho(E)$ must equal $k+1$. Thus far, we have for $\rho$:
    
    \[\begin{tikzcd}[row sep = 10]&&&\underset{\color{blue}E\color{black}}{ k+1}\arrow[d,no head]\arrow[rrd,no head,dotted]&&\\&\underset{\color{blue}efg\color{black}}{ \substack{3\\2}}\arrow[rru,no head]\arrow[d,no head]\arrow[rrd,no head,dotted]&\underset{\color{blue}efh\color{black}}{k+1}\arrow[ru,no head]\arrow[d,no head]\arrow[rrd,no head,dotted]&\underset{\color{blue}egh\color{black}}{k+1}\arrow[ld,no head]\arrow[rrd,no head,dotted]&&\underset{\color{blue}fgh\color{black}}{k+1}\arrow[ld,no head]\arrow[d,no head]\\\underset{\color{blue}ef\color{black}}{ 2}\arrow[ru,no head]\arrow[rru,no head]&\underset{\color{blue}eg\color{black}}{ 2}\arrow[rru,no head]&\underset{\color{blue}eh\color{black}}{ \substack{k+1\\k}}\arrow[rrd,no head,dotted]&\underset{\color{blue}fg\color{black}}{ 2}\arrow[rru,no head]&\underset{\color{blue}fh\color{black}}{ \substack{k+1\\k}}\arrow[d,no head]&\underset{\color{blue}gh\color{black}}{ \substack{k+1\\k}}\\\underset{\color{blue}e\color{black}}{ 1}\arrow[u,no head]\arrow[rru,no head]\arrow[ru,no head]&&\underset{\color{blue}f\color{black}}{ 1}\arrow[ru,no head]\arrow[rru,no head]\arrow[llu,no head,dotted]&\underset{\color{blue}g\color{black}}{ 1}\arrow[u,no head]\arrow[rru,no head]\arrow[llu,no head,dotted]&\underset{\color{blue}h\color{black}}{ k}\arrow[ru,no head]&\\&&\underset{\color{blue}\emptyset\color{black}}{0}\arrow[u,no head]\arrow[ru,no head]\arrow[rru,no head]\arrow[llu,no head,dotted]&&&\end{tikzcd}\]

    From here, we deduce that $\rho(efg) = 2$, because if $\rho(efg) = 3$, then $(\rho/ef) \notin \mathcal{P}_{U_{2,4}}^k$ by Lemma \ref{gamma}.

    Next, by submodularity of $\rho$, we require
    \begin{align*}
    \rho(eh) + \rho(fh) &\geq \rho(h) + \rho(efh)\\
    \rho(eh) + \rho(gh) &\geq \rho(h) + \rho(egh)\\
    \rho(fh) + \rho(gh) &\geq \rho(h) + \rho(fgh)
    \end{align*}
    implying that at least two of $\rho(eh)$, $\rho(fh)$, and $\rho(gh)$ must be equal to $k+1$. 

    \begin{enumerate}
        \item Assume exactly two of $\rho(eh)$, $\rho(fh)$, and $\rho(gh)$ are equal to $k+1$. Without loss of generality, let $\rho(eh) = k$ and $\rho(fh) = \rho(gh) = k+1$. This completely determines $\rho$:
        
            \begin{tikzcd}[row sep = 10]&&&\underset{\color{blue}E\color{black}}{ k+1}\arrow[d,no head]\arrow[rrd,no head,dotted]&&\\&\underset{\color{blue}efg\color{black}}{ 2}\arrow[rru,no head]\arrow[d,no head]\arrow[rrd,no head,dotted]&\underset{\color{blue}efh\color{black}}{ k+1}\arrow[ru,no head]\arrow[d,no head]\arrow[rrd,no head,dotted]&\underset{\color{blue}egh\color{black}}{ k+1}\arrow[ld,no head]\arrow[rrd,no head,dotted]&&\underset{\color{blue}fgh\color{black}}{ k+1}\arrow[ld,no head]\arrow[d,no head]\\\underset{\color{blue}ef\color{black}}{ 2}\arrow[ru,no head]\arrow[rru,no head]&\underset{\color{blue}eg\color{black}}{ 2}\arrow[rru,no head]&\underset{\color{blue}eh\color{black}}{ k}\arrow[rrd,no head,dotted]&\underset{\color{blue}fg\color{black}}{ 2}\arrow[rru,no head]&\underset{\color{blue}fh\color{black}}{ k+1}\arrow[d,no head]&\underset{\color{blue}gh\color{black}}{ k+1}\\\underset{\color{blue}e\color{black}}{ 1}\arrow[u,no head]\arrow[rru,no head]\arrow[ru,no head]&&\underset{\color{blue}f\color{black}}{ 1}\arrow[ru,no head]\arrow[rru,no head]\arrow[llu,no head,dotted]&\underset{\color{blue}g\color{black}}{ 1}\arrow[u,no head]\arrow[rru,no head]\arrow[llu,no head,dotted]&\underset{\color{blue}h\color{black}}{ k}\arrow[ru,no head]&\\&&\underset{\color{blue}\emptyset\color{black}}{0}\arrow[u,no head]\arrow[ru,no head]\arrow[rru,no head]\arrow[llu,no head,dotted]&&&\end{tikzcd}
            
See Figure \ref{fig: three points and a plane} (Left). To get a $U_{2,4}$-minor of $S(M_\rho^k)$, we must remove a set $Y\subseteq X_h$ of $k-1$ clones using a series of deletion or contraction operations. Let $h'$ be the remaining clone in $X_{h} - Y$. Let $Z \subseteq Y$ be the set of clones that is deleted and $Y-Z$ be the set that is contracted. See Figure \ref{fig: three points and a plane} (Center). Let $m = |Z|$. First, assume $1 \leq m \leq k-1$. That is, a nonempty set of clones is deleted. The total rank of $S(M_\rho^k)\backslash Z$ is $(k-m)+2$. This is because, by deleting $Z$, the remaining clones of $X_{h}$ have rank $k-m$; if we take the union of these remaining clones and $\{e,f,g\}$, the total rank is brought up to $(k-m)+2$. Then, the total rank of $S(M_\rho^k)\backslash Z/(Y-Z)$ is $((k-m)+2) - ((k-1)-m) = 3$. Therefore, $S(M_\rho^k)\backslash Z/(Y-Z)$ cannot be isomorphic to $U_{2,4}$. Finally, we consider $m = 0$, i.e. we contract all of $Y$ from $S(M_\rho^k)$. The total rank of $S(M_\rho^k)/Y$ is now $2$, but the rank of $\{e,h'\}$ in $S(M_\rho^k)/Y$ is equal to $k-(k-1) = 1$, which means $S(S(M_\rho^k)/Y)$ only contains $3$ points, and therefore cannot be isomorphic to $U_{2,4}$.

    \begin{figure}[H]
        \centering
        \includegraphics[scale = 0.25]{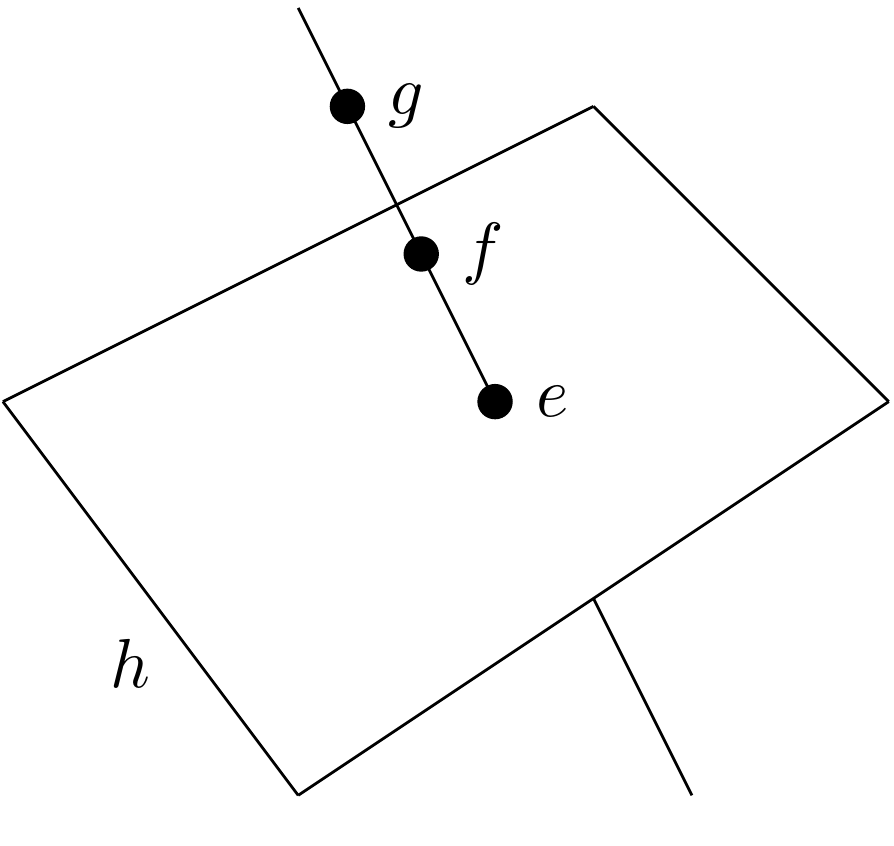}\hfill\includegraphics[scale = 0.25]{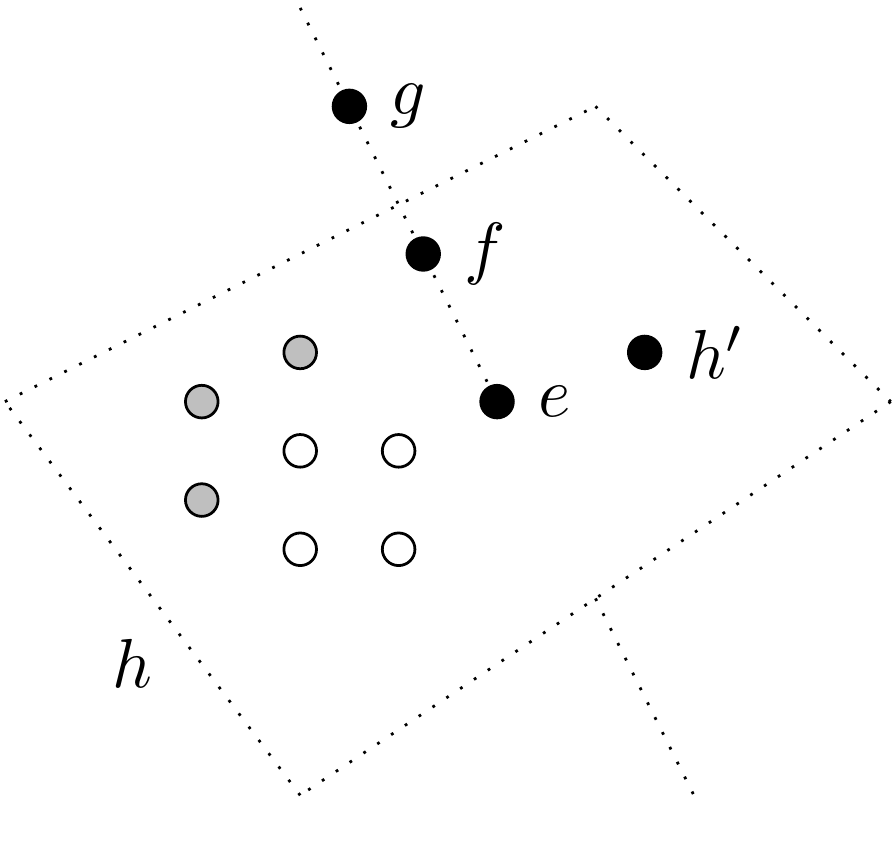}
        \hfill
        \includegraphics[scale = 0.25]{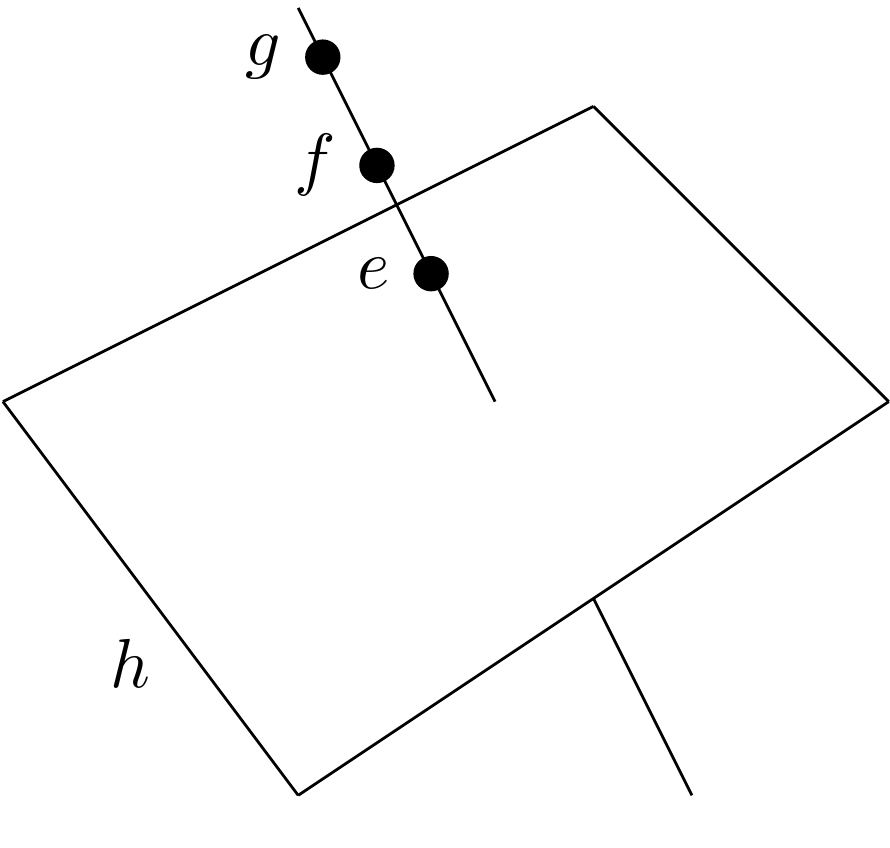}
        \caption{Left: The polymatroid $\rho_1$ consisting of three collinear points $e, f$, and $g$ which together with the rank-$k$ element $h$ span rank $k+1$, with $e$ (but neither $f$ nor $g$) lying on $h$. Center: An intuitive depiction of $S(M_{\rho_1}^k)$. The grey and white clones represent $Y$. The grey clones represent $Z$, the set of deleted clones, and the white clones represent $Y-Z$, the set of contracted clones. The only clone from $X_h$ which remains is $h'$. Right: The polymatroid $\rho_2$ consisting of three collinear points $e, f, g$, which together with the rank-$k$ element $h$ span rank $k+1$, with none of the points $e, f$ or $g$ lying on $h$.}
        \label{fig: three points and a plane}
    \end{figure}
    \item Assume $\rho(eh) = \rho(fh) = \rho(gh) = k+1$. This is $Ex_{(3, 0, 1)}^{k+1}$. See Figure \ref{fig: three points and a plane} (Right). $M_\rho^k$ has $k+3$ elements. We claim that $\rho$ is an excluded minor. To see this, first note that if we contract $k-1$ of the clones in $X_{h}$, then $S(M_\rho^k) \cong U_{2,4}$. It remains to show that every proper polymatroid minor of $\rho$ is in $\mathcal{P}_{U_{2,4}}^k$. We only need to consider deletions or contractions on one element of $E$. If we delete or contract $h$, then there will not be enough points in $S(M_{\rho\backslash h}^k)$ or $S(M_{\rho/h}^k)$ to contain a $U_{2,4}$ minor. If we delete any of the $e, f$, or $g$, we are left with a simple polymatroid of type $(2, 0, 1)$ in rank $k+1$, which we've already shown in Lemma \ref{201E3} must be in $\mathcal{P}_{U_{2,4}}^k$. If we contract any of the $e, f, g$ -- without loss of generality, say we contract on $e$ -- then $(\rho/e)(fg) = 1$. Thus, $S(M_{\rho/e}^k)$ consists of $k+1$ points lying freely in rank $k$. We must decrease the size of the ground set by $k-3$ and decrease the total rank by $k-2$. Since $k-3<k-2$, this is not possible. Therefore, we cannot create a $U_{2,4}$-minor of $M_{\rho/e}^k$, implying $(\rho/e) \in \mathcal{P}_{U_{2,4}}^k$.
    \end{enumerate}
\end{proof}

\begin{corollary}\label{other excluded minor E4 kdual}
    The $k$-dual of $Ex_{(3, 0, 1)}^{k+1}$, denoted $Ex_{(1, 0, 3)}^{3k-1}$ is an excluded minor. \end{corollary}
\vspace{2ex}

    \subsection{Type \texorpdfstring{$(2, 0, 2)$}{}}~\\
        
\begin{lemma}\label{202E4} 
The only excluded minor of type $T_\rho = (2, 0, 2)$ such that $\rho(E) \leq 2k$ is $Ex_{(2, 0, 2)}^{2k}$ which is identically self-$k$-dual.
\end{lemma}
\begin{proof}
Assume $\rho$ is an excluded minor of type $T_\rho = (2, 0, 2)$. Then $\rho$ must be simple. Let $\rho(e) = \rho(f) = 1$ and $\rho(g) = \rho(h) = k$. By Lemma \ref{beta}, $\rho(E) \geq \rho(gh) \geq 2k-1$, so we only need to consider two cases for $\rho(E)$: $\rho(E) = 2k-1$ and $\rho(E) = 2k$.
\begin{enumerate}
    \item Assume $\rho(E) = 2k-1$. We have:
    
    \begin{tikzcd}[row sep = 10]&&&\underset{\color{blue}E\color{black}}{ 2k-1}\arrow[d,no head]\arrow[rrd,no head,dotted]&&\\&\underset{\color{blue}efg\color{black}}{ \substack{k+2\\k+1}}\arrow[rru,no head]\arrow[d,no head]\arrow[rrd,no head,dotted]&\underset{\color{blue}efh\color{black}}{\substack{k+2\\k+1}}\arrow[ru,no head]\arrow[d,no head]\arrow[rrd,no head,dotted]&\underset{\color{blue}egh\color{black}}{ 2k-1}\arrow[ld,no head]\arrow[rrd,no head,dotted]&&\underset{\color{blue}fgh\color{black}}{ 2k-1}\arrow[ld,no head]\arrow[d,no head]\\\underset{\color{blue}ef\color{black}}{ 2}\arrow[ru,no head]\arrow[rru,no head]&\underset{\color{blue}eg\color{black}}{ \substack{k+1\\k}}\arrow[rru,no head]&\underset{\color{blue}eh\color{black}}{ \substack{k+1\\k}}\arrow[rrd,no head,dotted]&\underset{\color{blue}fg\color{black}}{\substack{k+1\\k}}\arrow[rru,no head]&\underset{\color{blue}fh\color{black}}{ \substack{k+1\\k}}\arrow[d,no head]&\underset{\color{blue}gh\color{black}}{ 2k-1}\\\underset{\color{blue}e\color{black}}{ 1}\arrow[u,no head]\arrow[rru,no head]\arrow[ru,no head]&&\underset{\color{blue}f\color{black}}{ 1}\arrow[ru,no head]\arrow[rru,no head]\arrow[llu,no head,dotted]&\underset{\color{blue}g\color{black}}{ k}\arrow[u,no head]\arrow[rru,no head]\arrow[llu,no head,dotted]&\underset{\color{blue}h\color{black}}{ k}\arrow[ru,no head]&\\&&\underset{\color{blue}\emptyset\color{black}}{0}\arrow[u,no head]\arrow[ru,no head]\arrow[rru,no head]\arrow[llu,no head,dotted]&&&\end{tikzcd}
    
It must be that $\rho(eg) = \rho(eh) = \rho(fg) = \rho(fh) = k$, otherwise $\rho\backslash e/f$ or $\rho\backslash f/e$ would not be in $\mathcal{P}_{U_{2,4}}^k$ by Lemmas \ref{beta} or \ref{epsilon}. But then submodularity of $\rho$ is violated: $\rho(eg) + \rho(fg) = 2k$ while $\rho(g) + \rho(efg) \geq 2k+1$. Therefore, $\rho(E)$ cannot equal $2k-1$.
    \item Now we consider $\rho(E) = 2k$. We have:
    
    \begin{tikzcd}[row sep = 10]&&&\underset{\color{blue}E\color{black}}{ 2k}\arrow[d,no head]\arrow[rrd,no head,dotted]&&\\&\underset{\color{blue}efg\color{black}}{ \substack{k+2\\k+1}}\arrow[rru,no head]\arrow[d,no head]\arrow[rrd,no head,dotted]&\underset{\color{blue}efh\color{black}}{ \substack{k+2\\k+1}}\arrow[ru,no head]\arrow[d,no head]\arrow[rrd,no head,dotted]&\underset{\color{blue}egh\color{black}}{ \substack{2k\\2k-1}}\arrow[ld,no head]\arrow[rrd,no head,dotted]&&\underset{\color{blue}fgh\color{black}}{ \substack{2k\\2k-1}}\arrow[ld,no head]\arrow[d,no head]\\\underset{\color{blue}ef\color{black}}{ 2}\arrow[ru,no head]\arrow[rru,no head]&\underset{\color{blue}eg\color{black}}{ \substack{k+1\\k}}\arrow[rru,no head]&\underset{\color{blue}eh\color{black}}{ \substack{k+1\\k}}\arrow[rrd,no head,dotted]&\underset{\color{blue}fg\color{black}}{ \substack{k+1\\k}}\arrow[rru,no head]&\underset{\color{blue}fh\color{black}}{ \substack{k+1\\k}}\arrow[d,no head]&\underset{\color{blue}gh\color{black}}{ \substack{2k\\2k-1}}\\\underset{\color{blue}e\color{black}}{ 1}\arrow[u,no head]\arrow[rru,no head]\arrow[ru,no head]&&\underset{\color{blue}f\color{black}}{ 1}\arrow[ru,no head]\arrow[rru,no head]\arrow[llu,no head,dotted]&\underset{\color{blue}g\color{black}}{ k}\arrow[u,no head]\arrow[rru,no head]\arrow[llu,no head,dotted]&\underset{\color{blue}h\color{black}}{ k}\arrow[ru,no head]&\\&&\underset{\color{blue}\emptyset\color{black}}{0}\arrow[u,no head]\arrow[ru,no head]\arrow[rru,no head]\arrow[llu,no head,dotted]&&&\end{tikzcd}

    We can further refine this by noting that we must have $\rho(efg) = \rho(efh) = k+1$; otherwise $(\rho/ef) \notin\mathcal{P}_{U_{2,4}}^k$ by Lemmas \ref{beta} and \ref{epsilon}:
    
    \begin{tikzcd}[row sep = 10]&&&\underset{\color{blue}E\color{black}}{ 2k}\arrow[d,no head]\arrow[rrd,no head,dotted]&&\\&\underset{\color{blue}efg\color{black}}{ k+1}\arrow[rru,no head]\arrow[d,no head]\arrow[rrd,no head,dotted]&\underset{\color{blue}efh\color{black}}{ k+1}\arrow[ru,no head]\arrow[d,no head]\arrow[rrd,no head,dotted]&\underset{\color{blue}egh\color{black}}{ \substack{2k\\2k-1}}\arrow[ld,no head]\arrow[rrd,no head,dotted]&&\underset{\color{blue}fgh\color{black}}{ \substack{2k\\2k-1}}\arrow[ld,no head]\arrow[d,no head]\\\underset{\color{blue}ef\color{black}}{ 2}\arrow[ru,no head]\arrow[rru,no head]&\underset{\color{blue}eg\color{black}}{ \substack{k+1\\k}}\arrow[rru,no head]&\underset{\color{blue}eh\color{black}}{ \substack{k+1\\k}}\arrow[rrd,no head,dotted]&\underset{\color{blue}fg\color{black}}{ \substack{k+1\\k}}\arrow[rru,no head]&\underset{\color{blue}fh\color{black}}{ \substack{k+1\\k}}\arrow[d,no head]&\underset{\color{blue}gh\color{black}}{ \substack{2k\\2k-1}}\\\underset{\color{blue}e\color{black}}{ 1}\arrow[u,no head]\arrow[rru,no head]\arrow[ru,no head]&&\underset{\color{blue}f\color{black}}{ 1}\arrow[ru,no head]\arrow[rru,no head]\arrow[llu,no head,dotted]&\underset{\color{blue}g\color{black}}{ k}\arrow[u,no head]\arrow[rru,no head]\arrow[llu,no head,dotted]&\underset{\color{blue}h\color{black}}{ k}\arrow[ru,no head]&\\&&\underset{\color{blue}\emptyset\color{black}}{0}\arrow[u,no head]\arrow[ru,no head]\arrow[rru,no head]\arrow[llu,no head,dotted]&&&\end{tikzcd}

Consider the three possibilities for $(\rho(egh), \rho(fgh))$ up to isomorphism:
    \begin{enumerate}
        \item Assume $(\rho(egh), \rho(fgh)) = (2k-1, 2k-1)$. This forces $\rho(gh) = 2k-1$ by monotonicity. Then $\rho(eg) = \rho(eh) = \rho(fg) = \rho(fh) = k$, otherwise $\rho\backslash e/f$ or $\rho\backslash f/e$ would not be in $\mathcal{P}_{U_{2,4}}^k$ by Lemmas \ref{beta} and \ref{epsilon}. This gives the following:
        
        \begin{tikzcd}[column sep = 23, row sep = 10]&&&\underset{\color{blue}E\color{black}}{ 2k}\arrow[d,no head]\arrow[rrd,no head,dotted]&&\\&\underset{\color{blue}efg\color{black}}{ k+1}\arrow[rru,no head]\arrow[d,no head]\arrow[rrd,no head,dotted]&\underset{\color{blue}efh\color{black}}{ k+1}\arrow[ru,no head]\arrow[d,no head]\arrow[rrd,no head,dotted]&\underset{\color{blue}egh\color{black}}{ 2k-1}\arrow[ld,no head]\arrow[rrd,no head,dotted]&&\underset{\color{blue}fgh\color{black}}{ 2k-1}\arrow[ld,no head]\arrow[d,no head]\\\underset{\color{blue}ef\color{black}}{ 2}\arrow[ru,no head]\arrow[rru,no head]&\underset{\color{blue}eg\color{black}}{ k}\arrow[rru,no head]&\underset{\color{blue}eh\color{black}}{ k}\arrow[rrd,no head,dotted]&\underset{\color{blue}fg\color{black}}{ k}\arrow[rru,no head]&\underset{\color{blue}fh\color{black}}{ k}\arrow[d,no head]&\underset{\color{blue}gh\color{black}}{ 2k-1}\\\underset{\color{blue}e\color{black}}{ 1}\arrow[u,no head]\arrow[rru,no head]\arrow[ru,no head]&&\underset{\color{blue}f\color{black}}{ 1}\arrow[ru,no head]\arrow[rru,no head]\arrow[llu,no head,dotted]&\underset{\color{blue}g\color{black}}{ k}\arrow[u,no head]\arrow[rru,no head]\arrow[llu,no head,dotted]&\underset{\color{blue}h\color{black}}{ k}\arrow[ru,no head]&\\&&\underset{\color{blue}\emptyset\color{black}}{0}\arrow[u,no head]\arrow[ru,no head]\arrow[rru,no head]\arrow[llu,no head,dotted]&&&\end{tikzcd}
        
But now $\rho$ violates submodularity: we have $\rho(egh) + \rho(fgh) = 4k-2$ while $\rho(gh)+\rho(E) = 4k-1$.
        \item Assume $(\rho(egh), \rho(fgh)) = (2k-1, 2k)$. This implies $\rho(gh) = 2k-1$ by monotonicity of $\rho$. It also must be that $\rho(eg) = \rho(eh) = k$; otherwise $\rho\backslash f/e$ would not be in $\mathcal{P}_{U_{2,4}}^k$ by Lemmas \ref{beta} and \ref{epsilon}. It must also be that $\rho(fg) = \rho(fh) = k+1$, otherwise submodularity of $\rho$ would be violated. To see this, assume $\rho(fg) = k$. Then $\rho(fg) + \rho(gh) < \rho(g) + \rho(fgh)$, so submodularity would be violated. Similarly, if $\rho(fh) = k$, then $\rho(fh) + \rho(gh) < \rho(h) + \rho(fgh)$, again violating submodularity. We have completely determined $\rho$:

        \begin{tikzcd}[column sep = 20, row sep = 10]&&&\underset{\color{blue}E\color{black}}{ 2k}\arrow[d,no head]\arrow[rrd,no head,dotted]&&\\&\underset{\color{blue}efg\color{black}}{ k+1}\arrow[rru,no head]\arrow[d,no head]\arrow[rrd,no head,dotted]&\underset{\color{blue}efh\color{black}}{ k+1}\arrow[ru,no head]\arrow[d,no head]\arrow[rrd,no head,dotted]&\underset{\color{blue}egh\color{black}}{ 2k-1}\arrow[ld,no head]\arrow[rrd,no head,dotted]&&\underset{\color{blue}fgh\color{black}}{ 2k}\arrow[ld,no head]\arrow[d,no head]\\\underset{\color{blue}ef\color{black}}{ 2}\arrow[ru,no head]\arrow[rru,no head]&\underset{\color{blue}eg\color{black}}{ k}\arrow[rru,no head]&\underset{\color{blue}eh\color{black}}{ k}\arrow[rrd,no head,dotted]&\underset{\color{blue}fg\color{black}}{ k+1}\arrow[rru,no head]&\underset{\color{blue}fh\color{black}}{ k+1}\arrow[d,no head]&\underset{\color{blue}gh\color{black}}{ 2k-1}\\\underset{\color{blue}e\color{black}}{ 1}\arrow[u,no head]\arrow[rru,no head]\arrow[ru,no head]&&\underset{\color{blue}f\color{black}}{ 1}\arrow[ru,no head]\arrow[rru,no head]\arrow[llu,no head,dotted]&\underset{\color{blue}g\color{black}}{ k}\arrow[u,no head]\arrow[rru,no head]\arrow[llu,no head,dotted]&\underset{\color{blue}h\color{black}}{ k}\arrow[ru,no head]&\\&&\underset{\color{blue}\emptyset\color{black}}{0}\arrow[u,no head]\arrow[ru,no head]\arrow[rru,no head]\arrow[llu,no head,dotted]&&&\end{tikzcd}

Notice that $\rho \cong \rho|_f \oplus \rho|_{egh}$, so $\rho$ cannot be an excluded minor for $\rho$ by Lemma \ref{connectedness}.
        \item Assume $(\rho(egh), \rho(fgh)) = (2k, 2k)$. We consider two subcases: (i) $\rho(gh) = 2k-1$, and (ii) $\rho(gh) = 2k$.
        \begin{enumerate}
            \item Assume $\rho(gh) = 2k-1$. It must be that $\rho(eg) = \rho(eh) = \rho(fg) = \rho(fh) = k+1$. Otherwise, without loss of generality, let $\rho(eg) = k$. Then submodularity of $\rho$ is violated because $\rho(eg) + \rho(gh) = 3k-1$ while $\rho(g) + \rho(egh) = 3k$. We have completely determined $\rho$:

            \begin{tikzcd}[column sep = 17, row sep = 10]&&&\underset{\color{blue}E\color{black}}{ 2k}\arrow[d,no head]\arrow[rrd,no head,dotted]&&\\&\underset{\color{blue}efg\color{black}}{ k+1}\arrow[rru,no head]\arrow[d,no head]\arrow[rrd,no head,dotted]&\underset{\color{blue}efh\color{black}}{ k+1}\arrow[ru,no head]\arrow[d,no head]\arrow[rrd,no head,dotted]&\underset{\color{blue}egh\color{black}}{ 2k}\arrow[ld,no head]\arrow[rrd,no head,dotted]&&\underset{\color{blue}fgh\color{black}}{ 2k}\arrow[ld,no head]\arrow[d,no head]\\\underset{\color{blue}ef\color{black}}{ 2}\arrow[ru,no head]\arrow[rru,no head]&\underset{\color{blue}eg\color{black}}{ k+1}\arrow[rru,no head]&\underset{\color{blue}eh\color{black}}{ k+1}\arrow[rrd,no head,dotted]&\underset{\color{blue}fg\color{black}}{ k+1}\arrow[rru,no head]&\underset{\color{blue}fh\color{black}}{ k+1}\arrow[d,no head]&\underset{\color{blue}gh\color{black}}{ 2k-1}\\\underset{\color{blue}e\color{black}}{ 1}\arrow[u,no head]\arrow[rru,no head]\arrow[ru,no head]&&\underset{\color{blue}f\color{black}}{ 1}\arrow[ru,no head]\arrow[rru,no head]\arrow[llu,no head,dotted]&\underset{\color{blue}g\color{black}}{ k}\arrow[u,no head]\arrow[rru,no head]\arrow[llu,no head,dotted]&\underset{\color{blue}h\color{black}}{ k}\arrow[ru,no head]&\\&&\underset{\color{blue}\emptyset\color{black}}{0}\arrow[u,no head]\arrow[ru,no head]\arrow[rru,no head]\arrow[llu,no head,dotted]&&&\end{tikzcd}

            This polymatroid is in $\mathcal{P}_{U_{2,4}}^k$. We leave it to the reader to check that the following $2k \times (2k+2)$ matrix is a binary representation of $S(M_\rho^k)$:
            
            \begin{align*}\hspace{20ex}
                \begin{blockarray}{cccccccccc}
            e&f&g_1&g_2&\hdots&g_k&h_1&h_2&\hdots&h_k \\\begin{block}{(c|c|cccc|cccc)}
            1 & 1 &  &  &  &  &  &  &  &  \\
            \cline{1-10}
             & 1 & 1 &  &  &  & 1 &  &  &  \\
             \cline{1-10}
             &  & 1 & 1 &  &  &  &  &  &  \\
             &  &  &  & \ddots &  &  &  &  &  \\
             &  &  &  &  & 1 &  &  &  &  \\
             \cline{1-10}
             &  &  &  &  &  & 1 & 1 &  &  \\
             &  &  &  &  &  &  &  & \ddots &  \\
             &  &  &  &  &  &  &  &  & 1 \\
            \end{block}\end{blockarray}\end{align*}
            \item $\rho(gh) = 2k$: At least one of $\rho(eg)$ and $\rho(eh)$ must equal $k+1$ and at least one of $\rho(fg)$ and $\rho(fh)$ must equal $k+1$. To see this, assume $\rho(eg) = \rho(eh) = k$. Then submodularity of $\rho$ would be violated: $\rho(eg)+\rho(eh) = 2k$ while $\rho(e)+\rho(egh) = 2k+1$. The argument for $\rho(fg)$ and $\rho(fh)$ is analogous. Up to isomorphism of $\rho$, there are four sub-cases to consider for $(\rho(eg),\rho(eh), \rho(fg),\rho(fh))$.
            \begin{enumerate}
                \item $(\rho(eg),\rho(eh), \rho(fg),\rho(fh))=  (k+1, k, k+1, k)$:  
                
                \begin{tikzcd}[column sep = 17, row sep = 10]&&&\underset{\color{blue}E\color{black}}{ 2k}\arrow[d,no head]\arrow[rrd,no head,dotted]&&\\&\underset{\color{blue}efg\color{black}}{ k+1}\arrow[rru,no head]\arrow[d,no head]\arrow[rrd,no head,dotted]&\underset{\color{blue}efh\color{black}}{ k+1}\arrow[ru,no head]\arrow[d,no head]\arrow[rrd,no head,dotted]&\underset{\color{blue}egh\color{black}}{ 2k}\arrow[ld,no head]\arrow[rrd,no head,dotted]&&\underset{\color{blue}fgh\color{black}}{ 2k}\arrow[ld,no head]\arrow[d,no head]\\\underset{\color{blue}ef\color{black}}{ 2}\arrow[ru,no head]\arrow[rru,no head]&\underset{\color{blue}eg\color{black}}{ k+1}\arrow[rru,no head]&\underset{\color{blue}eh\color{black}}{ k}\arrow[rrd,no head,dotted]&\underset{\color{blue}fg\color{black}}{ k+1}\arrow[rru,no head]&\underset{\color{blue}fh\color{black}}{ k}\arrow[d,no head]&\underset{\color{blue}gh\color{black}}{ 2k}\\\underset{\color{blue}e\color{black}}{ 1}\arrow[u,no head]\arrow[rru,no head]\arrow[ru,no head]&&\underset{\color{blue}f\color{black}}{ 1}\arrow[ru,no head]\arrow[rru,no head]\arrow[llu,no head,dotted]&\underset{\color{blue}g\color{black}}{ k}\arrow[u,no head]\arrow[rru,no head]\arrow[llu,no head,dotted]&\underset{\color{blue}h\color{black}}{ k}\arrow[ru,no head]&\\&&\underset{\color{blue}\emptyset\color{black}}{0}\arrow[u,no head]\arrow[ru,no head]\arrow[rru,no head]\arrow[llu,no head,dotted]&&&\end{tikzcd}
            
                This violates submodularity:  $\rho(eh)+\rho(fh) = 2k$ while $\rho(h)+\rho(efh) = 2k+1$.
                \item $(\rho(eg),\rho(eh), \rho(fg),\rho(fh)) = (k+1, k, k, k+1)$:
                
                \begin{tikzcd}[column sep = 17, row sep = 10]&&&\underset{\color{blue}E\color{black}}{ 2k}\arrow[d,no head]\arrow[rrd,no head,dotted]&&\\&\underset{\color{blue}efg\color{black}}{ k+1}\arrow[rru,no head]\arrow[d,no head]\arrow[rrd,no head,dotted]&\underset{\color{blue}efh\color{black}}{ k+1}\arrow[ru,no head]\arrow[d,no head]\arrow[rrd,no head,dotted]&\underset{\color{blue}egh\color{black}}{ 2k}\arrow[ld,no head]\arrow[rrd,no head,dotted]&&\underset{\color{blue}fgh\color{black}}{ 2k}\arrow[ld,no head]\arrow[d,no head]\\\underset{\color{blue}ef\color{black}}{ 2}\arrow[ru,no head]\arrow[rru,no head]&\underset{\color{blue}eg\color{black}}{ k+1}\arrow[rru,no head]&\underset{\color{blue}eh\color{black}}{ k}\arrow[rrd,no head,dotted]&\underset{\color{blue}fg\color{black}}{ k}\arrow[rru,no head]&\underset{\color{blue}fh\color{black}}{ k+1}\arrow[d,no head]&\underset{\color{blue}gh\color{black}}{ 2k}\\\underset{\color{blue}e\color{black}}{ 1}\arrow[u,no head]\arrow[rru,no head]\arrow[ru,no head]&&\underset{\color{blue}f\color{black}}{ 1}\arrow[ru,no head]\arrow[rru,no head]\arrow[llu,no head,dotted]&\underset{\color{blue}g\color{black}}{ k}\arrow[u,no head]\arrow[rru,no head]\arrow[llu,no head,dotted]&\underset{\color{blue}h\color{black}}{ k}\arrow[ru,no head]&\\&&\underset{\color{blue}\emptyset\color{black}}{0}\arrow[u,no head]\arrow[ru,no head]\arrow[rru,no head]\arrow[llu,no head,dotted]&&&\end{tikzcd}
            
Then $\rho \cong \rho|_{eh} \oplus \rho|_{fg}$ so $\rho$ is not an excluded minor by Lemma \ref{connectedness}.
                \item $(\rho(eg),\rho(eh), \rho(fg),\rho(fh)) = (k+1, k+1, k+1, k)$:
                
                \begin{tikzcd}[column sep = 17, row sep = 10]&&&\underset{\color{blue}E\color{black}}{ 2k}\arrow[d,no head]\arrow[rrd,no head,dotted]&&\\&\underset{\color{blue}efg\color{black}}{ k+1}\arrow[rru,no head]\arrow[d,no head]\arrow[rrd,no head,dotted]&\underset{\color{blue}efh\color{black}}{ k+1}\arrow[ru,no head]\arrow[d,no head]\arrow[rrd,no head,dotted]&\underset{\color{blue}egh\color{black}}{ 2k}\arrow[ld,no head]\arrow[rrd,no head,dotted]&&\underset{\color{blue}fgh\color{black}}{ 2k}\arrow[ld,no head]\arrow[d,no head]\\\underset{\color{blue}ef\color{black}}{ 2}\arrow[ru,no head]\arrow[rru,no head]&\underset{\color{blue}eg\color{black}}{ k+1}\arrow[rru,no head]&\underset{\color{blue}eh\color{black}}{ k+1}\arrow[rrd,no head,dotted]&\underset{\color{blue}fg\color{black}}{ k+1}\arrow[rru,no head]&\underset{\color{blue}fh\color{black}}{ k}\arrow[d,no head]&\underset{\color{blue}gh\color{black}}{ 2k}\\\underset{\color{blue}e\color{black}}{ 1}\arrow[u,no head]\arrow[rru,no head]\arrow[ru,no head]&&\underset{\color{blue}f\color{black}}{ 1}\arrow[ru,no head]\arrow[rru,no head]\arrow[llu,no head,dotted]&\underset{\color{blue}g\color{black}}{ k}\arrow[u,no head]\arrow[rru,no head]\arrow[llu,no head,dotted]&\underset{\color{blue}h\color{black}}{ k}\arrow[ru,no head]&\\&&\underset{\color{blue}\emptyset\color{black}}{0}\arrow[u,no head]\arrow[ru,no head]\arrow[rru,no head]\arrow[llu,no head,dotted]&&&\end{tikzcd}

                This polymatroid is in $\mathcal{P}_{U_{2,4}}^k$. We leave it to the reader to check that the following $(2k) \times (2k+2)$ matrix is a binary representation of $S(M_\rho^k)$:
                \begin{align*}\hspace{25ex}\begin{blockarray}{cccccccc}
e&f&g_1&\hdots&g_k&h_1&\hdots&h_k \\\begin{block}{(c|c|ccc|ccc)}
1 &  & 1 &  &  &  &  &  \\
\vdots &  &  & \ddots &  &  &  &  \\
1 &  &  &  & 1 &  &  &  \\
\cline{1-8}
1 & 1 &  &  &  & 1 &  &  \\
\vdots & \vdots &  &  &  &  & \ddots &  \\
1 & 1 &  &  &  &  &  & 1 \\
\end{block}\end{blockarray}\end{align*}

                \item $(\rho(eg),\rho(eh), \rho(fg),\rho(fh)) = (k+1, k+1, k+1, k+1)$:

                \begin{tikzcd}[column sep = 13,  row sep = 10]&&&\underset{\color{blue}E\color{black}}{ 2k}\arrow[d,no head]\arrow[rrd,no head,dotted]&&\\&\underset{\color{blue}efg\color{black}}{ k+1}\arrow[rru,no head]\arrow[d,no head]\arrow[rrd,no head,dotted]&\underset{\color{blue}efh\color{black}}{ k+1}\arrow[ru,no head]\arrow[d,no head]\arrow[rrd,no head,dotted]&\underset{\color{blue}egh\color{black}}{ 2k}\arrow[ld,no head]\arrow[rrd,no head,dotted]&&\underset{\color{blue}fgh\color{black}}{ 2k}\arrow[ld,no head]\arrow[d,no head]\\\underset{\color{blue}ef\color{black}}{ 2}\arrow[ru,no head]\arrow[rru,no head]&\underset{\color{blue}eg\color{black}}{ k+1}\arrow[rru,no head]&\underset{\color{blue}eh\color{black}}{ k+1}\arrow[rrd,no head,dotted]&\underset{\color{blue}fg\color{black}}{ k+1}\arrow[rru,no head]&\underset{\color{blue}fh\color{black}}{ k+1}\arrow[d,no head]&\underset{\color{blue}gh\color{black}}{ 2k}\\\underset{\color{blue}e\color{black}}{ 1}\arrow[u,no head]\arrow[rru,no head]\arrow[ru,no head]&&\underset{\color{blue}f\color{black}}{ 1}\arrow[ru,no head]\arrow[rru,no head]\arrow[llu,no head,dotted]&\underset{\color{blue}g\color{black}}{ k}\arrow[u,no head]\arrow[rru,no head]\arrow[llu,no head,dotted]&\underset{\color{blue}h\color{black}}{ k}\arrow[ru,no head]&\\&&\underset{\color{blue}\emptyset\color{black}}{0}\arrow[u,no head]\arrow[ru,no head]\arrow[rru,no head]\arrow[llu,no head,dotted]&&&\end{tikzcd}

                We claim that this is an excluded minor. To get a $U_{2,4}$-minor $M$ of $M_\rho^k$, contract $k-1$ clones each from $X_g$ and $X_h$. It remains to show that every proper polymatroid minor of $\rho$ is in $\mathcal{P}_{U_{2,4}}^k$. We only need to consider deletions and contractions on one element. By symmetry, we only need to consider deletions and contractions on $e$ and $g$. If we delete $e$, then the following is a binary representation for $S(M_{\rho\backslash e}^k)$:

                \begin{align*}\hspace{25ex}\begin{blockarray}{ccccccc}
f&g_1&\hdots&g_k&h_1&\hdots&h_k \\\begin{block}{(c|ccc|ccc)}
1  & 1 &  &  &  &  &  \\
\vdots  &  & \ddots &  &  &  &  \\
1  &  &  & 1 &  &  &  \\
\cline{1-7}
1 &   &  &  & 1 &  &  \\
\vdots &  &  &  &  & \ddots &  \\
1 &   &  &  &  &  & 1 \\
\end{block}\end{blockarray}\end{align*}

                If we contract $e$, then the following is a binary representation for $S(M_{\rho/e}^k)$:

                \begin{align*}\hspace{25ex}\begin{blockarray}{ccccccccc}
f&g_1&g_2&\hdots&g_k&h_1&h_2&\hdots&h_k \\\begin{block}{(c|cccc|cccc)}
 1 & 1 &  &  &  & 1 &  &  &  \\
 \cline{1-9}
& 1 & 1 &  &  &  &  &  &  \\
&  &  & \ddots &  &  &  &  &  \\
&  &  &  & 1 &  &  &  &  \\
 \cline{1-9}
&  &  &  &  & 1 & 1 &  &  \\
 &  &  &  &  &  &  & \ddots &  \\
&  &  &  &  &  &  &  & 1 \\
\end{block}\end{blockarray}\end{align*}

                If we delete $g$, then $\rho\backslash g$ is a simple polymatroid of type $(2, 0, 1)$ in rank $k+1$ which by Lemma \ref{201E3} is in $\mathcal{P}_{U_{2,4}}^k$. Finally, if we contract $g$, then the following is a binary representation for $S(M_{\rho/g}^k)$:
                \begin{align*}\hspace{25ex}
        \begin{blockarray}{ccccc}
                e&f&h_1&\hdots&h_k\\\begin{block}{(c|c|ccc)}
1 & 1 & 1 &  &  \\
\vdots & \vdots &  & \ddots &  \\
1 & 1 &  &  & 1 \\
\end{block}\end{blockarray}
                \end{align*}
                
            \end{enumerate}
        \end{enumerate}
    \end{enumerate}
\end{enumerate}
\end{proof}
\vspace{2ex}

    \subsection{Remaining types}~\\
        
\begin{lemma}\label{040E4 031E4 130E4} There are no excluded minors $(E,\rho)$ for $\mathcal{P}_{U_{2,4}}^k$ of types $(0, 4, 0), (0, 3, 1)$, or $(1, 3, 0)$ such that $\rho(E) \leq 2k$.
\end{lemma}
\begin{proof} Let $\rho$ be an excluded minor. The type $(0, 4, 0)$ contains $4$ elements of rank $k-1$, so by Corollary \ref{skew kminus1 elements}, $\rho(E) \geq 4k-4$. Since $4k-4 > 2k$, we can eliminate this type as a possibility for $\rho$. The types $(0, 3, 1)$ and $(1, 3, 0)$ each contain $3$ elements of rank $k-1$, so by Corollary \ref{skew kminus1 elements}, $\rho(E) \geq 3k-3$. When $k \geq 4$, we have $3k-3>2k$, so we only need to investigate these two types when $k = 3$. If $k = 3$, we only need to consider $\rho(E) = 3k-3 = 2k = 6$.
\begin{enumerate}
    \item Assume $T_\rho = (0, 3, 1)$. Let $\rho(e) = \rho(f) = \rho(g) = 2$ and $\rho(h) = 3$. Then, by Corollary \ref{skew kminus1 elements} and Lemma \ref{epsilon}, we have for $\rho\backslash e$:

    \begin{minipage}{\linewidth}
        $$\begin{tikzcd}[row sep = 10]&\underset{\color{blue}fgh\color{black}}{ \substack{6\\5}}&\\\underset{\color{blue}fg\color{black}}{ 4}\arrow[ru,no head]&\underset{\color{blue}fh\color{black}}{ 5}\arrow[u,no head]&\underset{\color{blue}gh\color{black}}{ 5}\arrow[lu,no head]\\\underset{\color{blue}f\color{black}}{ 2}\arrow[u,no head]\arrow[ru,no head]&\underset{\color{blue}g\color{black}}{ 2}\arrow[lu,no head]\arrow[ru,no head]&\underset{\color{blue}h\color{black}}{ 3}\arrow[lu,no head]\arrow[u,no head]\\&\underset{\color{blue}\emptyset\color{black}}{0}\arrow[lu,no head]\arrow[u,no head]\arrow[ru,no head]&\end{tikzcd}$$
    \end{minipage}
        
    If $\rho(fgh) = 5$, then $\rho\backslash e/f \cong Ex_\epsilon^3$. If $\rho(fgh) = 6$, then $\rho\backslash e/f \cong Ex_\epsilon^4$. Either way, $\rho$ cannot be an excluded minor for $\mathcal{P}_{U_{2,4}}^3$.
    \item Assume $T_\rho = (1, 3, 0)$. Let $\rho(e) = 1$ and $\rho(f) = \rho(g) = \rho(h) = k-1$. We have:

    \begin{minipage}{\linewidth}
    $$\begin{tikzcd}[row sep = 10]&&&\underset{\color{blue}E\color{black}}{ 6}\arrow[d,no head]\arrow[rrd,no head,dotted]&&\\&\underset{\color{blue}efg\color{black}}{ \substack{5\\4}}\arrow[rru,no head]\arrow[d,no head]\arrow[rrd,no head,dotted]&\underset{\color{blue}efh\color{black}}{ \substack{5\\4}}\arrow[ru,no head]\arrow[d,no head]\arrow[rrd,no head,dotted]&\underset{\color{blue}egh\color{black}}{ \substack{5\\4}}\arrow[ld,no head]\arrow[rrd,no head,dotted]&&\underset{\color{blue}fgh\color{black}}{ 6}\arrow[ld,no head]\arrow[d,no head]\\\underset{\color{blue}ef\color{black}}{ 3}\arrow[ru,no head]\arrow[rru,no head]&\underset{\color{blue}eg\color{black}}{ 3}\arrow[rru,no head]&\underset{\color{blue}eh\color{black}}{ 3}\arrow[rrd,no head,dotted]&\underset{\color{blue}fg\color{black}}{ 4}\arrow[rru,no head]&\underset{\color{blue}fh\color{black}}{ 4}\arrow[d,no head]&\underset{\color{blue}gh\color{black}}{ 4}\\\underset{\color{blue}e\color{black}}{ 1}\arrow[u,no head]\arrow[rru,no head]\arrow[ru,no head]&&\underset{\color{blue}f\color{black}}{ 2}\arrow[ru,no head]\arrow[rru,no head]\arrow[llu,no head,dotted]&\underset{\color{blue}g\color{black}}{ 2}\arrow[u,no head]\arrow[rru,no head]\arrow[llu,no head,dotted]&\underset{\color{blue}h\color{black}}{ 2}\arrow[ru,no head]&\\&&\underset{\color{blue}\emptyset\color{black}}{0}\arrow[u,no head]\arrow[ru,no head]\arrow[rru,no head]\arrow[llu,no head,dotted]&&&\end{tikzcd}$$
    \end{minipage}

If $\rho(efg) = 4$ or $\rho(efh) = 4$, then $\rho\backslash h/f \cong Ex_\gamma^2$ or $\rho\backslash g/f \cong Ex_\gamma^2$ respectively. This gives $\rho(efg) = \rho(efh) = 5$. But then $\rho/ef \cong Ex_\alpha^3$. Thus, $\rho$ cannot be an excluded minor for $\mathcal{P}_{U_{2,4}}^3$.
\end{enumerate}
\end{proof}

\begin{lemma}\label{310E4} There are no excluded minors $(E,\rho)$ for $\mathcal{P}_{U_{2,4}}^k$ such that $T_\rho = (3, 1, 0)$ and $\rho(E) \leq 2k$.
\end{lemma}
\begin{proof} Assume $\rho$ is an excluded minor. Then, $\rho$ must be simple and connected. Let $\rho(e) = \rho(f) = \rho(g) = 1$ and $\rho(h) = k-1$. We must have $k-1 \leq \rho(E) \leq k+1$.
\begin{enumerate}
    \item Assume $\rho(E) = k-1$, which implies $\rho(eh) = k-1$. Then $\rho|_{eh} \notin \mathcal{P}_{U_{2,4}}^k$ by Lemma \ref{gamma}, so $\rho$ cannot be an excluded minor.
    \item Assume $\rho(E) = k$. Then $\rho(fg) = 2$ and $\rho(fh) = \rho(gh) = \rho(fgh) = k$. Observe that $\rho\backslash e/f \notin \mathcal{P}_{U_{2,4}}^k$ by Lemma \ref{gamma}.
    \item Assume $\rho(E) = k+1$. We have:

    \begin{minipage}{\linewidth}
    $$\begin{tikzcd}[row sep = 10]&&&\underset{\color{blue}E\color{black}}{ k+1}\arrow[d,no head]\arrow[rrd,no head,dotted]&&\\&\underset{\color{blue}efg\color{black}}{ \substack{3\\2}}\arrow[rru,no head]\arrow[d,no head]\arrow[rrd,no head,dotted]&\underset{\color{blue}efh\color{black}}{ \substack{k+1\\k}}\arrow[ru,no head]\arrow[d,no head]\arrow[rrd,no head,dotted]&\underset{\color{blue}egh\color{black}}{ \substack{k+1\\k}}\arrow[ld,no head]\arrow[rrd,no head,dotted]&&\underset{\color{blue}fgh\color{black}}{ \substack{k+1\\k}}\arrow[ld,no head]\arrow[d,no head]\\\underset{\color{blue}ef\color{black}}{ 2}\arrow[ru,no head]\arrow[rru,no head]&\underset{\color{blue}eg\color{black}}{ 2}\arrow[rru,no head]&\underset{\color{blue}eh\color{black}}{ k}\arrow[rrd,no head,dotted]&\underset{\color{blue}fg\color{black}}{ 2}\arrow[rru,no head]&\underset{\color{blue}fh\color{black}}{ k}\arrow[d,no head]&\underset{\color{blue}gh\color{black}}{ k}\\\underset{\color{blue}e\color{black}}{ 1}\arrow[u,no head]\arrow[rru,no head]\arrow[ru,no head]&&\underset{\color{blue}f\color{black}}{ 1}\arrow[ru,no head]\arrow[rru,no head]\arrow[llu,no head,dotted]&\underset{\color{blue}g\color{black}}{ 1}\arrow[u,no head]\arrow[rru,no head]\arrow[llu,no head,dotted]&\underset{\color{blue}h\color{black}}{ k-1}\arrow[ru,no head]&\\&&\underset{\color{blue}\emptyset\color{black}}{0}\arrow[u,no head]\arrow[ru,no head]\arrow[rru,no head]\arrow[llu,no head,dotted]&&&\end{tikzcd}$$
    \end{minipage}

    If $\rho(A) = k$ where $A = efh, egh$, or $fgh$, then $\rho|_A \notin \mathcal{P}_{U_{2,4}}^k$ by Lemma \ref{210E3}. If $\rho(efh) =\rho(egh) = \rho(fgh) = k+1$ and $\rho(efg) = 2$, then $\rho \cong \rho|_{efg} \oplus \rho|_h$ so $\rho$ cannot be an excluded minor by Lemma \ref{connectedness}. If $\rho(efh) =\rho(egh) = \rho(fgh) = k+1$ and $\rho(efg) = 3$, then $(\rho/ef) \notin \mathcal{P}_{U_{2,4}}^k$ by Lemma \ref{gamma}. 
\end{enumerate}
\end{proof}

\begin{lemma}\label{220E4} There are no excluded minors $(E,\rho)$ for $\mathcal{P}_{U_{2,4}}^k$ such that $T_\rho = (2, 2, 0)$ and $\rho(E) \leq 2k$.
\end{lemma}
\begin{proof} By Lemma \ref{connectedness} and Corollary \ref{skew kminus1 elements}, we only need to consider $\rho(E) = 2k-2$ or $\rho(E) = 2k-1$. We have:

$$\begin{tikzcd}[row sep = 10]&&&\underset{\color{blue}E\color{black}}{ \substack{2k-1\\2k-2}}\arrow[d,no head]\arrow[rrd,no head,dotted]&&\\&\underset{\color{blue}efg\color{black}}{ k+1}\arrow[rru,no head]\arrow[d,no head]\arrow[rrd,no head,dotted]&\underset{\color{blue}efh\color{black}}{ k+1}\arrow[ru,no head]\arrow[d,no head]\arrow[rrd,no head,dotted]&\underset{\color{blue}egh\color{black}}{ \substack{2k-1\\2k-2}}\arrow[ld,no head]\arrow[rrd,no head,dotted]&&\underset{\color{blue}fgh\color{black}}{ \substack{2k-1\\2k-2}}\arrow[ld,no head]\arrow[d,no head]\\\underset{\color{blue}ef\color{black}}{ 2}\arrow[ru,no head]\arrow[rru,no head]&\underset{\color{blue}eg\color{black}}{ k}\arrow[rru,no head]&\underset{\color{blue}eh\color{black}}{ k}\arrow[rrd,no head,dotted]&\underset{\color{blue}fg\color{black}}{ k}\arrow[rru,no head]&\underset{\color{blue}fh\color{black}}{ k}\arrow[d,no head]&\underset{\color{blue}gh\color{black}}{ 2k-2}\\\underset{\color{blue}e\color{black}}{ 1}\arrow[u,no head]\arrow[rru,no head]\arrow[ru,no head]&&\underset{\color{blue}f\color{black}}{ 1}\arrow[ru,no head]\arrow[rru,no head]\arrow[llu,no head,dotted]&\underset{\color{blue}g\color{black}}{ k-1}\arrow[u,no head]\arrow[rru,no head]\arrow[llu,no head,dotted]&\underset{\color{blue}h\color{black}}{ k-1}\arrow[ru,no head]&\\&&\underset{\color{blue}\emptyset\color{black}}{0}\arrow[u,no head]\arrow[ru,no head]\arrow[rru,no head]\arrow[llu,no head,dotted]&&&\end{tikzcd}$$

It must be that $\rho(egh) = \rho(fgh) = 2k-1$. To see this, assume without loss of generality that $\rho(fgh) = 2k-2$. Then $(\rho\backslash e/f)\notin\mathcal{P}_{U_{2,4}}^k$ by Lemma \ref{alpha}. By monotonicity of $\rho$, it must be that $\rho(E) = 2k-1$, but then $(\rho/ef) \notin \mathcal{P}_{U_{2,4}}^k$ by Lemma \ref{alpha}. 
\end{proof}

The remaining types each contain a proper restriction of type $T_\rho = (0, 1, 1)$ or $(0, 0, 2)$. Therefore, by Lemmas \ref{beta} and \ref{epsilon}, we only need to consider $\rho(E) \geq 2k-1$.

\begin{lemma}\label{004E4 103E4 013E4} There are no excluded minors $(E,\rho)$ for $\mathcal{P}_{U_{2,4}}^k$ of types $(0, 0, 4), (1, 0, 3)$, or $(0, 1, 3)$ such that $\rho(E) \leq 2k$.
\end{lemma}
\begin{proof}
Assume $\rho$ is an excluded minor of type $(0, 0, 4), (1, 0, 3)$, or $(0, 1, 3)$. Let $\rho(f) = \rho(g) = \rho(h) = k$. 
\begin{enumerate}
    \item Assume $\rho(E) = 2k-1$. By Lemma \ref{beta}, $\rho\backslash e$ must be the following:

\begin{minipage}{\linewidth}
$$\begin{tikzcd}[row sep = 10]&\underset{\color{blue}fgh\color{black}}{ 2k-1}&\\\underset{\color{blue}fg\color{black}}{ 2k-1}\arrow[ru,no head]&\underset{\color{blue}fh\color{black}}{ 2k-1}\arrow[u,no head]&\underset{\color{blue}gh\color{black}}{ 2k-1}\arrow[lu,no head]\\\underset{\color{blue}f\color{black}}{ k}\arrow[u,no head]\arrow[ru,no head]&\underset{\color{blue}g\color{black}}{ k}\arrow[lu,no head]\arrow[ru,no head]&\underset{\color{blue}h\color{black}}{ k}\arrow[lu,no head]\arrow[u,no head]\\&\underset{\color{blue}\emptyset\color{black}}{0}\arrow[lu,no head]\arrow[u,no head]\arrow[ru,no head]&\end{tikzcd}$$
\end{minipage}

We observe that $\rho\backslash e/f \cong Ex_\alpha^{k-1}$.
    \item Assume $\rho(E) = 2k$. Then $\rho\backslash e$ must be the following:

\begin{minipage}{\linewidth}
$$\begin{tikzcd}[row sep = 10]&\underset{\color{blue}fgh\color{black}}{\substack{2k\\2k-1}}&\\\underset{\color{blue}fg\color{black}}{ \substack{2k\\2k-1}}\arrow[ru,no head]&\underset{\color{blue}fh\color{black}}{ \substack{2k\\2k-1}}\arrow[u,no head]&\underset{\color{blue}gh\color{black}}{ \substack{2k\\2k-1}}\arrow[lu,no head]\\\underset{\color{blue}f\color{black}}{ k}\arrow[u,no head]\arrow[ru,no head]&\underset{\color{blue}g\color{black}}{ k}\arrow[lu,no head]\arrow[ru,no head]&\underset{\color{blue}h\color{black}}{ k}\arrow[lu,no head]\arrow[u,no head]\\&\underset{\color{blue}\emptyset\color{black}}{0}\arrow[lu,no head]\arrow[u,no head]\arrow[ru,no head]&\end{tikzcd}$$
\end{minipage}

Depending on what we choose for $\rho(fg), \rho(fh)$, and $\rho(fgh)$, we observe that $\rho\backslash e/f$ would be isomorphic to one of the following: $Ex_\alpha^{k-1}$, $Ex_\alpha^k$, $Ex_\beta^k$, $Ex_\epsilon^k$. Hence $\rho$ cannot be an excluded minor.
\end{enumerate}
\end{proof}

\begin{lemma}\label{022E4 112E4} There are no excluded minors $(E,\rho)$ for $\mathcal{P}_{U_{2,4}}^k$ of types $(0, 2, 2)$ or $(1, 1, 2)$ such that $\rho(E) \leq 2k$.
\end{lemma}
\begin{proof} Assume $\rho$ is an excluded minor of type $(0, 2, 2)$ or $(1, 1, 2)$. Let $\rho(f) = k-1$, and $\rho(g) = \rho(h) = k$. 
\begin{enumerate}
    \item Assume $\rho(E) = 2k-1$. By Lemmas \ref{beta} and \ref{epsilon}, $\rho\backslash e$ must be the following:

\begin{minipage}{\linewidth}
$$\begin{tikzcd}[row sep = 10]&\underset{\color{blue}fgh\color{black}}{ 2k-1}&\\\underset{\color{blue}fg\color{black}}{ 2k-1}\arrow[ru,no head]&\underset{\color{blue}fh\color{black}}{ 2k-1}\arrow[u,no head]&\underset{\color{blue}gh\color{black}}{ 2k-1}\arrow[lu,no head]\\\underset{\color{blue}f\color{black}}{ k-1}\arrow[u,no head]\arrow[ru,no head]&\underset{\color{blue}g\color{black}}{ k}\arrow[lu,no head]\arrow[ru,no head]&\underset{\color{blue}h\color{black}}{ k}\arrow[lu,no head]\arrow[u,no head]\\&\underset{\color{blue}\emptyset\color{black}}{0}\arrow[lu,no head]\arrow[u,no head]\arrow[ru,no head]&\end{tikzcd}$$
\end{minipage}

We observe that $\rho\backslash e/f \cong Ex_\beta^k$.
\item Assume $\rho(E) = 2k$. By Lemmas \ref{beta} and \ref{epsilon}, $\rho\backslash e$ must be the following:

\begin{minipage}{\linewidth}
$$\begin{tikzcd}[row sep = 10]&\underset{\color{blue}fgh\color{black}}{\substack{2k\\2k-1}}&\\\underset{\color{blue}fg\color{black}}{ 2k-1}\arrow[ru,no head]&\underset{\color{blue}fh\color{black}}{ 2k-1}\arrow[u,no head]&\underset{\color{blue}gh\color{black}}{ \substack{2k\\2k-1}}\arrow[lu,no head]\\\underset{\color{blue}f\color{black}}{ k-1}\arrow[u,no head]\arrow[ru,no head]&\underset{\color{blue}g\color{black}}{ k}\arrow[lu,no head]\arrow[ru,no head]&\underset{\color{blue}h\color{black}}{ k}\arrow[lu,no head]\arrow[u,no head]\\&\underset{\color{blue}\emptyset\color{black}}{0}\arrow[lu,no head]\arrow[u,no head]\arrow[ru,no head]&\end{tikzcd}$$
\end{minipage}

If $\rho(fgh) = 2k-1$, then $\rho\backslash e/f \cong Ex_\beta^k$. If $\rho(fgh) = 2k$, then $\rho \backslash e/f \cong Ex_\beta^{k+1}$. Hence, there is no such $\rho$.
\end{enumerate}
\end{proof}

\begin{lemma}\label{121E4} There are no excluded minors $(E,\rho)$ for $\mathcal{P}_{U_{2,4}}^k$ such that $T_\rho = (1, 2, 1)$ and $\rho(E) \leq 2k$.
\end{lemma}
\begin{proof} Assume $\rho$ is an excluded minor of type $(1, 2, 1)$. Let $\rho(f) = \rho(g) = k-1$ and $\rho(h) = k$. Then $\rho\backslash e$ is the following:

$$\begin{tikzcd}[row sep = 10]&\underset{\color{blue}fgh\color{black}}{\substack{2k\\2k-1}}&\\\underset{\color{blue}fg\color{black}}{ 2k-2}\arrow[ru,no head]&\underset{\color{blue}fh\color{black}}{ 2k-1}\arrow[u,no head]&\underset{\color{blue}gh\color{black}}{ 2k-1}\arrow[lu,no head]\\\underset{\color{blue}f\color{black}}{ k-1}\arrow[u,no head]\arrow[ru,no head]&\underset{\color{blue}g\color{black}}{ k-1}\arrow[lu,no head]\arrow[ru,no head]&\underset{\color{blue}h\color{black}}{ k}\arrow[lu,no head]\arrow[u,no head]\\&\underset{\color{blue}\emptyset\color{black}}{0}\arrow[lu,no head]\arrow[u,no head]\arrow[ru,no head]&\end{tikzcd}$$

If $\rho(fgh) = 2k-1$, we observe that $\rho\backslash e/f \cong Ex_\epsilon^k$. If $\rho(fgh) = 2k$, then we observe that $\rho\backslash e /f \notin \mathcal{P}_{U_{2,4}}^k$ by Lemma \ref{epsilon}. Hence, there is no excluded minor of type $T_\rho = (1, 2, 1)$.
\end{proof}

\begin{lemma}\label{211E4} There are no excluded minors $(E,\rho)$ for $\mathcal{P}_{U_{2,4}}^k$ such that $T_\rho = (2, 1, 1)$ and $\rho(E) \leq 2k$.
\end{lemma}
\begin{proof} Assume $\rho$ is an excluded minor of type $(2, 1, 1)$. Let $\rho(e) = \rho(f) = 1$, $\rho(g) = k-1$, and $\rho(h) = k$. 
\begin{enumerate}
    \item Assume $\rho(E) = 2k-1$. We have:

\begin{minipage}{\linewidth}
    $$\begin{tikzcd}[column sep = 23, row sep = 10]&&&\underset{\color{blue}E\color{black}}{ 2k-1}\arrow[d,no head]\arrow[rrd,no head,dotted]&&\\&\underset{\color{blue}efg\color{black}}{ k+1}\arrow[rru,no head]\arrow[d,no head]\arrow[rrd,no head,dotted]&\underset{\color{blue}efh\color{black}}{ k+1}\arrow[ru,no head]\arrow[d,no head]\arrow[rrd,no head,dotted]&\underset{\color{blue}egh\color{black}}{ 2k-1}\arrow[ld,no head]\arrow[rrd,no head,dotted]&&\underset{\color{blue}fgh\color{black}}{ 2k-1}\arrow[ld,no head]\arrow[d,no head]\\\underset{\color{blue}ef\color{black}}{ 2}\arrow[ru,no head]\arrow[rru,no head]&\underset{\color{blue}eg\color{black}}{ k}\arrow[rru,no head]&\underset{\color{blue}eh\color{black}}{ \substack{k+1\\k}}\arrow[rrd,no head,dotted]&\underset{\color{blue}fg\color{black}}{ k}\arrow[rru,no head]&\underset{\color{blue}fh\color{black}}{ \substack{k+1\\k}}\arrow[d,no head]&\underset{\color{blue}gh\color{black}}{ 2k-1}\\\underset{\color{blue}e\color{black}}{ 1}\arrow[u,no head]\arrow[rru,no head]\arrow[ru,no head]&&\underset{\color{blue}f\color{black}}{ 1}\arrow[ru,no head]\arrow[rru,no head]\arrow[llu,no head,dotted]&\underset{\color{blue}g\color{black}}{ k-1}\arrow[u,no head]\arrow[rru,no head]\arrow[llu,no head,dotted]&\underset{\color{blue}h\color{black}}{ k}\arrow[ru,no head]&\\&&\underset{\color{blue}\emptyset\color{black}}{0}\arrow[u,no head]\arrow[ru,no head]\arrow[rru,no head]\arrow[llu,no head,dotted]&&&\end{tikzcd}$$
\end{minipage}

Observe that $(\rho/eg) \notin \mathcal{P}_{U_{2,4}}^k$ by Lemma \ref{gamma}, so $\rho$ cannot be an excluded minor.
\item If $\rho(E) = 2k$, we have:

\begin{minipage}{\linewidth}
$$\begin{tikzcd}[column sep = 23, row sep = 10]&&&\underset{\color{blue}E\color{black}}{ 2k}\arrow[d,no head]\arrow[rrd,no head,dotted]&&\\&\underset{\color{blue}efg\color{black}}{ k+1}\arrow[rru,no head]\arrow[d,no head]\arrow[rrd,no head,dotted]&\underset{\color{blue}efh\color{black}}{ k+1}\arrow[ru,no head]\arrow[d,no head]\arrow[rrd,no head,dotted]&\underset{\color{blue}egh\color{black}}{ \substack{2k\\2k-1}}\arrow[ld,no head]\arrow[rrd,no head,dotted]&&\underset{\color{blue}fgh\color{black}}{ \substack{2k\\2k-1}}\arrow[ld,no head]\arrow[d,no head]\\\underset{\color{blue}ef\color{black}}{ 2}\arrow[ru,no head]\arrow[rru,no head]&\underset{\color{blue}eg\color{black}}{ k}\arrow[rru,no head]&\underset{\color{blue}eh\color{black}}{ \substack{k+1\\k}}\arrow[rrd,no head,dotted]&\underset{\color{blue}fg\color{black}}{ k}\arrow[rru,no head]&\underset{\color{blue}fh\color{black}}{ \substack{k+1\\k}}\arrow[d,no head]&\underset{\color{blue}gh\color{black}}{ 2k-1}\\\underset{\color{blue}e\color{black}}{ 1}\arrow[u,no head]\arrow[rru,no head]\arrow[ru,no head]&&\underset{\color{blue}f\color{black}}{ 1}\arrow[ru,no head]\arrow[rru,no head]\arrow[llu,no head,dotted]&\underset{\color{blue}g\color{black}}{ k-1}\arrow[u,no head]\arrow[rru,no head]\arrow[llu,no head,dotted]&\underset{\color{blue}h\color{black}}{ k}\arrow[ru,no head]&\\&&\underset{\color{blue}\emptyset\color{black}}{0}\arrow[u,no head]\arrow[ru,no head]\arrow[rru,no head]\arrow[llu,no head,dotted]&&&\end{tikzcd}$$
\end{minipage}

If $\rho(eh) = \rho(fh) = k$, then $\rho(eh) + \rho(fh) < \rho(h) + \rho(efh)$, violating submodularity. We investigate the remaining possible options for $(\rho(eh), \rho(fh), \rho(egh), \rho(fgh))$ up to isomorphism of $\rho$.
    \begin{enumerate}
        \item $\rho(eh) = k$ and $\rho(fh) = k+1$:
        \begin{enumerate}
            \item Let $\rho(fgh) = 2k-1$. Then $(\rho\backslash e/f) \notin \mathcal{P}_{U_{2,4}}^k$ by Lemma \ref{epsilon} so $\rho$ cannot be an excluded minor.
            \item Let $\rho(fgh) = 2k$. 
            \begin{enumerate}
                \item If $\rho(egh) = 2k-1$, then $\rho \cong \rho|_f \oplus \rho|_{egh}$ so $\rho$ cannot be an excluded minor by Lemma \ref{connectedness}.
                \item If $\rho(egh) = 2k$, then $\rho(eh) + \rho(gh) < \rho(h) + \rho(egh)$, violating submodularity.
            \end{enumerate}
        \end{enumerate}
        \item $\rho(eh) = \rho(fh) = k+1$: 
        \begin{enumerate}
            \item Let $\rho(fgh) = 2k-1$. Then $(\rho\backslash e/f) \notin \mathcal{P}_{U_{2,4}}^k$ by Lemma \ref{epsilon}, so $\rho$ cannot be an excluded minor.
            \item Let $\rho(fgh) = 2k$. 
            \begin{enumerate}
                \item If $\rho(egh) = 2k-1$, then $(\rho\backslash f/e) \notin \mathcal{P}_{U_{2,4}}^k$ by Lemma \ref{epsilon}, so $\rho$ cannot be an excluded minor.
                \item If $\rho(egh) = 2k$, then $\rho \cong \rho|_g \oplus \rho|_{efh}$, so $\rho$ cannot be an excluded minor by Lemma \ref{connectedness}.
            \end{enumerate}
        \end{enumerate}
    \end{enumerate}
\end{enumerate}
\end{proof}
This concludes the proof of Theorem \ref{excluded minors |E| = 4}. \hspace{28ex} $\square$

\section{Decompressions}
\label{section:decompressions}

In this section, we will complete step (2) of Strategy \ref{strategy} by proving Theorem \ref{no decompressions}. It suffices to only consider decompressions of the excluded minors $\rho$ such that $|E(\rho)| = 4$. This is because decompressions $\rho'$ of the excluded minors $\rho$ such that $|E(\rho)| = 2$ have $|E(\rho')| = 3$, and by Theorem \ref{noE3} there are no excluded minors on a ground set of size $3$.

\begin{theorem}\label{no decompressions}
If $(E,\rho)$ is an excluded minor for $\mathcal{P}_{U_{2,4}}^k$ on a ground set of size $4$, then no decompression of $\rho$ is an excluded minor.
\end{theorem}
\begin{proof}
    By Lemmas \ref{kduality} and \ref{compression commutes with dual}, we only need to consider decompressions of $U_{2,4}$, $Ex_{(3, 0, 1)}^{k+1}$, and $Ex_{(2, 0, 2)}^{2k}$. The justification for Theorem \ref{no decompressions} follows from the lemmas in this section.
\begin{lemma}\label{decompressionofpt}
    Let $(E, \rho)$ be a polymatroid which is an excluded minor on $|E| \geq 4$. If $e \in E$ has rank $k$ or $k-1$, and $\rho_{\downarrow e}^l$ contains points $p_1, \hdots, p_n$, then the decompression of $p_i$ must be a point in $\rho$ for every $1 \leq i \leq n$. 
\end{lemma}
\begin{proof}
    Let the decompression of $p_i$ be $f$ in $\rho$. Since  $|E| \geq 4$, $\rho|_{f}$ would be a proper restriction of $\rho$ and therefore $\rho(f) \in \{1, k-1, k\}$. For the following cases, let $Y = \{e_1, \hdots, e_l\}$ be the subset of clones of $X_e$ that is contracted and $Z = X_e-Y$ be the subset of clones that is deleted from $M_\rho^k$ to get $\rho_{\downarrow e}^l$. Let $r$ be the rank function of $M_\rho^k$. By definition of $l$-compression, we have $\rho_{\downarrow e}^l = r\backslash Z / Y$, so 
    \begin{displaymath}
        (\rho_{\downarrow e}^l)(p_i) = (r\backslash Z)(Y \cup X_f) - (r\backslash Z)(Y). 
    \end{displaymath}
    By definition of deletion, the above equals
    \begin{displaymath}
         = r(Y \cup X_f) - r(Y). 
    \end{displaymath}
    Since $p_i$ is a point, $(\rho_{\downarrow e}^l)(p_i) = 1$, so we have $r(Y \cup X_f) - r(Y) = 1$. Since $1 \leq |Y| \leq k-1$, and $Y$ is a set of $l$ clones lying freely in $e$, which is assumed to have rank $k$ or $k-1$, it must be that $l = |Y| = r(Y)$. So $r(Y \cup X_f) = l+1$. The bounds on $l$ imply $2 \leq r(Y \cup X_f) \leq k$. 
    \begin{enumerate}
        \item If $\rho(f) = k$, then $r(Y \cup X_f) = k$. This implies $l = k-1$, so $Y$ consists of $k-1$ points lying freely in the span of $X_f$, which is rank $k$. This also implies $|Z| = 1$. Hence, $r(Y \cup X_f \cup Z) = \rho(ef) = k$ or $k+1$. If $\rho(e) = k$, then $\rho|_{ef}$ is of type $\beta$, so $\rho|_{ef}$ would be $Ex_\beta^k$ or $Ex_\beta^{k+1}$. If $\rho(e) = k-1$, then $\rho|_{ef}$ is of type $\epsilon$, so $\rho|_{ef}\notin\mathcal{P}_{U_{2,4}}^k$ by Lemma \ref{epsilon}. This contradicts $\rho$ being an excluded minor.
        \item If $\rho(f) = k-1$, then $k-1 \leq r(Y \cup X_f) \leq k$, so $l = k-1$ or $k-2$.
        \begin{enumerate}
            \item If $l = k-1$, then $Y$ consists of $k-1$ points lying freely in the span of $X_f$, which is rank $k-1$. This implies $|Z| = 1$. Hence, $r(Y \cup X_f \cup Z) = \rho(ef) = k-1$ or $k$. If $\rho(e) = k$, then $\rho(ef)$ must equal $k$. Since $\rho|_{ef}$ is of type $\epsilon$, it must be that $\rho|_{ef} = Ex_\epsilon^k$. If $\rho(e) = k-1$, then $\rho|_{ef}$ is of type $\alpha$, and $\rho|_{ef}$ would equal $Ex_\alpha^{k-1}$ or $Ex_\alpha^k$. This case contradicts $\rho$ being an excluded minor.
            \item If $l = k-2$, then $Y$ consists of $k-2$ points lying freely in the span of $X_f$ which is rank $k-1$. This implies $|Z| = 2$. Hence, $r(Y \cup X_f \cup Z) = \rho(ef) = k-1$, $k$, or $k+1$. If $\rho(e) = k$, then $\rho(ef)$ must equal $k$ or $k+1$. Since $\rho|_{ef}$ is of type $\epsilon$, it must be that $\rho|_{ef}\notin \mathcal{P}_{U_{2,4}}^k$ by Lemma \ref{epsilon}. If $\rho(e) = k-1$, then $\rho|_{ef}$ is of type $\alpha$. It must be that $\rho|_{ef}\notin \mathcal{P}_{U_{2,4}}^k$ by Lemma \ref{alpha}. This case contradicts $\rho$ being an excluded minor.
        \end{enumerate}
    \end{enumerate}
    We conclude that $\rho(f) = 1$. 
\end{proof}

For the rest of this section, let $(E,\rho')$ be a $k$-polymatroid on $\{d, e, f, g, h\}$ such that some $l$-compression of $\rho'$ by $d$ is equal to $\rho$, that is, $(\rho')_{\downarrow d}^l = \rho$. Equivalently, $\rho'$ is some $l$-decompression of $\rho$ by $d$. Let $r'$ be the rank function of $M_{\rho'}^k$ and $r$ be the rank function of $M_\rho^k$. The $l$-compression of a point is simply a deletion or contraction of that point, so we will only need to consider the cases $\rho'(d) = k$ and $\rho'(d) = k-1$ in each of our proofs. Finally, for any rank-$k$ element of $\rho$, we only need to consider when its decompression in $\rho'$ also has rank $k$, because $k$-polymatroids cannot have singletons whose ranks are higher than $k$.

\begin{lemma}\label{nodecomps1}
    No decompression $(E,\rho')$ of $U_{2,4}$ is an excluded minor.
\end{lemma}
\begin{proof}
    Assume $\rho'$ is an excluded minor. By Lemma \ref{decompressionofpt}, $\rho'(e) = \rho'(f) = \rho'(g) = \rho'(h) = 1$. Let $Y = \{e, f, g, h\}$.
    \begin{enumerate}
        \item Let $\rho'(d) = k$. By Lemmas \ref{connectedness} and \ref{201E3}, it must be that $k+1 \leq \rho'(E) \leq k+3$. Note that $|E(S(M_{\rho'}^k))| = k+4$.
        \begin{enumerate}
            \item Let $\rho'(E) = k+3$. To get a $U_{2,4}$-minor of $M_{\rho'}^k$ we must decrease the size of the ground set by $k$ and decrease the total rank by $k+1$. Since $k < k+1$, this is not possible.
            \item Let $\rho'(E) = k+2$. It must be that $\rho'(Y)$ is equal to $3$ or $4$ because if $\rho'(Y) =2$, then $\rho'|_Y \cong U_{2,4}$. In the compression from $\rho'$ to $\rho$, let $Z_1 \subseteq X_d$ be the contracted set and let $Z_2\subseteq X_d$ be the deleted set. $Z_1$ and $Z_2$ are disjoint and nonempty. This gives $1 \leq |Z_1|, |Z_2| \leq k-1$. Let $|Z_1| = l$. We know $r'(Z_1) = l$ since $Z_1$ consists of $l$ points lying freely in rank $l$.
            \begin{enumerate}
                \item If $\rho'(Y) = 3$, then $(r'\backslash Z_2)(X_Y \cup Z_1) = 3+l$. Then 
                \begin{minipage}{\linewidth}
                \begin{align*}
                    r(X_Y) &= (r'/Z_1\backslash Z_2)(X_Y)\\
                    &= (r'\backslash Z_2)(X_Y \cup Z_1) - (r'\backslash Z_2)(Z_1)\\
                    &= (3+l)-l = 3.
                \end{align*} 
                \end{minipage}
                But $r(X_Y) = 2$ because $\rho\cong U_{2,4}$, contradiction.
                \item If $\rho'(Y) = 4$, then if $1 \leq l \leq k-2$, we have $(r'\backslash Z_2)(X_Y \cup Z_1)= 4+l$, so
                
                \begin{minipage}{\linewidth}
                    \begin{align*}
                        r(X_Y) = (4+l)-l = 4.
                    \end{align*}
                \end{minipage}
                
                If $l = k-1$, then we have $(r'\backslash Z_2)(X_Y \cup Z_1)= k+2$, so
                
                \begin{minipage}{\linewidth}
                    \begin{align*}
                        r(X_Y) &= (k+2)-(k-1) = 3.
                    \end{align*}
                \end{minipage}
                
                But $r(X_Y) = 2$ because $\rho\cong U_{2,4}$, contradiction.
            \end{enumerate}
            \item Let $\rho'(E) = k+1$. It must be that $\rho'(ef) = 2$, and it cannot be that both $\rho'(efg)$ and $\rho'(efh) = 2$; otherwise we would have $\rho'(Y) = 2$ by submodularity and monotonicity, and then $\rho'|_Y \cong U_{2,4}$. Assume without loss of generality that $\rho'(efg) = 3$. Then $(\rho'/ef)(g) = 1$. By Lemma \ref{201E3} and the fact that $\rho'(E) = k+1$, it must be that $\rho'(def) = k+1$. Thus $(\rho'/ef)(d) = k-1$. Finally, it must also be that $\rho'(defg) = k+1$, so $(\rho'/ef)(dg) = k-1$. By Lemma \ref{gamma}, $\rho'\backslash h/ef \notin \mathcal{P}_{U_{2,4}}^k$, contradiction.
            \end{enumerate}
        \item Let $\rho'(d) = k-1$. Then by Lemmas \ref{connectedness} and \ref{210E3}, $\rho'(E)$ can be equal to $k+1$ or $k+2$.
        \begin{enumerate}
            \item Let $\rho'(E) = k+2$. It must be that $\rho'(Y)$ is equal to $3$ or $4$; otherwise if $\rho'(Y) = 2$, then $\rho'|_Y \cong U_{2,4}$, contradiction.
            \begin{enumerate}
                \item If $\rho'(Y) = 4$, then we have a contradiction by the same reasoning as in part (1)(b)(ii) of this proof.
                \item If $\rho'(Y) = 3$, then $\rho' = (\rho')|_Y \oplus (\rho')|_d$, contradiction by Lemma \ref{connectedness}.
            \end{enumerate}
            \item Let $\rho'(E) = k+1$. It must be that $\rho'(ef) = 2$ and it cannot be that both $\rho'(efg)$ and $\rho'(efh) = 2$; otherwise we would have $\rho'(Y) = 2$ by submodularity and $\rho'|_Y \cong U_{2,4}$.  Assume without loss of generality that $\rho'(efg) = 3$. Then $(\rho'/ef)(g) = 1$. By Lemma \ref{201E3} and the fact that $\rho'(E) = k+1$, it must be that $\rho'(def) = \rho'(defg) = k+1$. This gives $(\rho'/ef)(d) = (\rho'/ef)(dg) = k-1$. By Lemma \ref{gamma}, $\rho'\backslash h/ef \notin \mathcal{P}_{U_{2,4}}^k$ so $\rho'$ cannot be an excluded minor.
        \end{enumerate}
    \end{enumerate}
\end{proof}

\begin{lemma}\label{nodecomps2}
    No decompression $(E, \rho')$ of $Ex_{(3, 0, 1)}^{k+1}$ is an excluded minor.
\end{lemma}
\begin{proof}
    Assume that $\rho'$ an excluded minor. By Lemma \ref{decompressionofpt}, $\rho'(e) = \rho'(f) = \rho'(g) = 1$ and $\rho'(h) = k$. Let $Y = \{e, f, g, h\}$.
    \begin{enumerate}
        \item Let $\rho'(d) = k$. By Lemmas \ref{connectedness} and \ref{beta}, it must be that $2k-1 \leq \rho'(E) \leq 2k+2$. Note that $|E(S(M_{\rho'}^k))| = 2k+3$. 
        \begin{enumerate}
            \item Let $\rho'(E) = 2k+2$. To get a $U_{2,4}$-minor of $M_{\rho'}^k$ we must decrease the size of the ground set by $2k-1$ and decrease the total rank by $2k$. Since $2k-1<2k$, this is not possible.
            \item Let $\rho'(E) = 2k+1$. By Lemma \ref{201E3}, it must be that $k+1 \leq \rho'(Y) \leq k+3$.
            \begin{enumerate}
                \item If $\rho'(Y) = k+3$, then we have a contradiction by the same reasoning as in part (1)(b)(ii) of the proof of Lemma \ref{nodecomps1}.
                \item If $\rho'(Y) = k+2$, then we have a contradiction by the same reasoning as in part (1)(b)(i) of the proof of Lemma \ref{nodecomps1}.
                \item If $\rho'(Y) = k+1$, then $\rho' = (\rho')|_{Y} \oplus (\rho')|_d$, contradiction by Lemma \ref{connectedness}.
            \end{enumerate}
            \item Let $\rho'(E) = 2k$ or $2k-1$. By Lemma \ref{201E3}, it must be that $\rho'(efh) = k+1$ or $\rho'(efh) = k+2$, and that $\rho'(def) = k+1$ or $\rho'(def) = k+2$. Since $\rho'$ is a decompression of $\rho = Ex_{(3, 0, 1)}^k$, and $\rho(eh) = \rho(fh) = k+1$, it must be that $\rho'(eh)$ and $\rho'(fh)$ are at least $k+1$. But since $e$ and $f$ are points, and $h$ is a rank-$k$ element, $\rho'(eh)$ and $\rho'(fh)$ are at most $k+1$. Thus, $\rho'(eh) = \rho'(fh) = k+1$. We have the following for $\rho'\backslash g$:

            \begin{minipage}{\linewidth}
            \begin{tikzcd}[row sep = 10]&&&\underset{\color{blue}defh\color{black}}{ \substack{2k\\2k-1}}\arrow[d,no head]\arrow[rrd,no head,dotted]&&\\&\underset{\color{blue}def\color{black}}{ \substack{k+2\\k+1}}\arrow[rru,no head]\arrow[d,no head]\arrow[rrd,no head,dotted]&\underset{\color{blue}deh\color{black}}{ \substack{2k\\2k-1}}\arrow[ru,no head]\arrow[d,no head]\arrow[rrd,no head,dotted]&\underset{\color{blue}dfh\color{black}}{ \substack{2k\\2k-1}}\arrow[ld,no head]\arrow[rrd,no head,dotted]&&\underset{\color{blue}efh\color{black}}{ \substack{k+2\\k+1}}\arrow[ld,no head]\arrow[d,no head]\\\underset{\color{blue}de\color{black}}{ \substack{k+1\\k}}\arrow[ru,no head]\arrow[rru,no head]&\underset{\color{blue}df\color{black}}{ \substack{k+1\\k}}\arrow[rru,no head]&\underset{\color{blue}dh\color{black}}{ \substack{2k\\2k-1}}\arrow[rrd,no head,dotted]&\underset{\color{blue}ef\color{black}}{ 2}\arrow[rru,no head]&\underset{\color{blue}eh\color{black}}{ k+1}\arrow[d,no head]&\underset{\color{blue}fh\color{black}}{ k+1}\\\underset{\color{blue}d\color{black}}{ k}\arrow[u,no head]\arrow[rru,no head]\arrow[ru,no head]&&\underset{\color{blue}e\color{black}}{ 1}\arrow[ru,no head]\arrow[rru,no head]\arrow[llu,no head,dotted]&\underset{\color{blue}f\color{black}}{ 1}\arrow[u,no head]\arrow[rru,no head]\arrow[llu,no head,dotted]&\underset{\color{blue}h\color{black}}{ k}\arrow[ru,no head]&\\&&\underset{\color{blue}\emptyset\color{black}}{0}\arrow[u,no head]\arrow[ru,no head]\arrow[rru,no head]\arrow[llu,no head,dotted]&&&\end{tikzcd}                
            \end{minipage}
            
            Observe that if $\rho'(defh) = 2k-1$, then $\rho'\backslash g/ef \notin \mathcal{P}_{U_{2,4}}^k$ by Lemmas \ref{alpha}, \ref{beta}, and \ref{epsilon}. So it must be that $\rho'(defh) = 2k$ implying $\rho'(E) = 2k$. Now we have $(\rho'/ef)(dh) = 2k-2$. If $\rho'(efh) = k+2$ or $\rho'(def) = k+2$, then $\rho'\backslash g/ef \notin \mathcal{P}_{U_{2,4}}^k$ by Lemmas \ref{beta} and \ref{epsilon}. Hence, it must be that $\rho'(efh) =\rho'(def) = k+1$. By Lemma \ref{201E3}, it must be that $\rho'(Y) \geq k+1$, and $\rho'(Y)$ cannot exceed the sum of the ranks of the singletons of $Y$, so $\rho'(Y) \leq k+3$. 
            \begin{enumerate}
                \item If $\rho'(Y) = k+3$, then $(\rho'/efh)(g) = \rho'(Y) - \rho'(efh) = 2$, which is not possible since $\rho'(g) = 1$. 
                \item If $\rho'(Y) = k+2$, then $\rho'/efh \notin \mathcal{P}_{U_{2,4}}^k$ by Lemma \ref{gamma}.
                \item If $\rho'(Y) = k+1$, then we have a contradiction by the same reasoning as in part (1)(b)(i) of the proof of Lemma \ref{nodecomps1}.
            \end{enumerate}
        \end{enumerate}
        \item Let $\rho'(d) = k-1$. By Lemmas \ref{epsilon} and \ref{connectedness}, it must be that $2k-1 \leq \rho'(E) \leq 2k+1$. 
        \begin{enumerate}
            \item Let $\rho'(E) = 2k+1$. By Lemma \ref{201E3}, it must be that $k+1 \leq \rho'(Y) \leq k+3$.
            \begin{enumerate}
                \item If $\rho'(Y) = k+3$, then we have a contradiction by the same reasoning as in part (1)(b)(ii) of the proof of Lemma \ref{nodecomps1}.
                \item If $\rho'(Y) = k+2$, then $\rho' = (\rho')|_Y \oplus (\rho')|_d$, a contradiction by Lemma \ref{connectedness}.
                \item If $\rho'(Y) = k+1$, then it is not possible for the total rank $\rho'(E)$ to be $2k+1$ so this case does not give an excluded minor.
            \end{enumerate}
            \item Let $\rho'(E) = 2k$ or $2k-1$. By Lemma \ref{201E3}, it must be that $\rho'(efh) = k+1$ or $\rho'(efh) = k+2$, and by Lemma \ref{210E3}, we have $\rho'(def) = k+1$. We have:

            \begin{minipage}{\linewidth}
                \begin{tikzcd}[column sep = 23, row sep = 10]&&&\underset{\color{blue}defh\color{black}}{ \substack{2k\\2k-1}}\arrow[d,no head]\arrow[rrd,no head,dotted]&&\\&\underset{\color{blue}def\color{black}}{ k+1}\arrow[rru,no head]\arrow[d,no head]\arrow[rrd,no head,dotted]&\underset{\color{blue}deh\color{black}}{ \substack{2k\\2k-1}}\arrow[ru,no head]\arrow[d,no head]\arrow[rrd,no head,dotted]&\underset{\color{blue}dfh\color{black}}{ \substack{2k\\2k-1}}\arrow[ld,no head]\arrow[rrd,no head,dotted]&&\underset{\color{blue}efh\color{black}}{ \substack{k+2\\k+1}}\arrow[ld,no head]\arrow[d,no head]\\\underset{\color{blue}de\color{black}}{ k}\arrow[ru,no head]\arrow[rru,no head]&\underset{\color{blue}df\color{black}}{ k}\arrow[rru,no head]&\underset{\color{blue}dh\color{black}}{ 2k-1}\arrow[rrd,no head,dotted]&\underset{\color{blue}ef\color{black}}{ 2}\arrow[rru,no head]&\underset{\color{blue}eh\color{black}}{ \substack{k+1\\k}}\arrow[d,no head]&\underset{\color{blue}fh\color{black}}{ \substack{k+1\\k}}\\\underset{\color{blue}d\color{black}}{k-1}\arrow[u,no head]\arrow[rru,no head]\arrow[ru,no head]&&\underset{\color{blue}e\color{black}}{ 1}\arrow[ru,no head]\arrow[rru,no head]\arrow[llu,no head,dotted]&\underset{\color{blue}f\color{black}}{ 1}\arrow[u,no head]\arrow[rru,no head]\arrow[llu,no head,dotted]&\underset{\color{blue}h\color{black}}{ k}\arrow[ru,no head]&\\&&\underset{\color{blue}\emptyset\color{black}}{0}\arrow[u,no head]\arrow[ru,no head]\arrow[rru,no head]\arrow[llu,no head,dotted]&&&\end{tikzcd}
            \end{minipage}

            Observe that if $\rho'(defh) = 2k-1$, then $(\rho\backslash g/ef)\notin \mathcal{P}_{U_{2,4}}^k$ by Lemmas \ref{alpha} and \ref{epsilon}, so it must be that $\rho'(defh) = 2k$, implying $\rho'(E) = 2k$ as well. The same arguments from part 1(c) of this proof hold here.
        \end{enumerate}
    \end{enumerate}
\end{proof}

\begin{lemma}\label{nodecomps3}
    No decompression $(E,\rho')$ of $Ex_{(2, 0, 2)}^{2k}$ is an excluded minor.
\end{lemma}
\begin{proof}
        Assume that $\rho'$ is an excluded minor. Let $\rho'(e) = \rho'(f) = 1$ and $\rho'(g) = \rho'(h) = k$. For $S(M_{\rho'}^k)$ to contain a $U_{2,4}$-minor, we would have to decrease its ground set size by $3k-2$ and decrease its rank by $\rho'(E)-2$. This would not be possible if $\rho'(E) \geq 3k+1$ because then $3k-2 < \rho'(E)-2$. Hence, it must be that $\rho'(E) \leq 3k$. Now consider the restriction $\rho'|_{dgh}$; let it have total rank $m$. Since $((\rho')^l_{\downarrow d})(gh) = 2k$, which is the largest rank two rank-$k$ elements can span, it must be that $(\rho')(gh) = 2k$ as well. We have:
            $$\begin{tikzcd}[row sep = 10]&\underset{\color{blue}dgh\color{black}}{ m}&\\\underset{\color{blue}dg\color{black}}{ \substack{2k\\2k-1}}\arrow[ru,no head]&\underset{\color{blue}dh\color{black}}{ \substack{2k\\2k-1}}\arrow[u,no head]&\underset{\color{blue}gh\color{black}}{2k}\arrow[lu,no head]\\\underset{\color{blue}d\color{black}}{ \substack{k\\k-1}}\arrow[u,no head]\arrow[ru,no head]&\underset{\color{blue}g\color{black}}{ k}\arrow[lu,no head]\arrow[ru,no head]&\underset{\color{blue}h\color{black}}{k}\arrow[lu,no head]\arrow[u,no head]\\&\underset{\color{blue}\emptyset\color{black}}{0}\arrow[lu,no head]\arrow[u,no head]\arrow[ru,no head]&\end{tikzcd}$$

    It must be that $3k-1\leq m \leq 3k$; otherwise $(\rho'|_{dgh})/g$ would contain an excluded minor by Lemmas \ref{beta} and \ref{epsilon}. This implies $\rho'(E) \geq 3k-1$. First, assume $\rho'(E) = 3k$. Let $Y = \{e, f, g, h\}$. The possibilities for $\rho'(Y)$ are $2k, 2k+1, 2k+2$. 
    \begin{enumerate}
        \item Assume $\rho'(Y) = 2k$. It is not possible for $\rho'(d) = k-1$ because $\rho'(E) = 3k$. Hence, $\rho'(d) = k$. This implies $\rho'= (\rho')|_Y \oplus (\rho')_d$ so $\rho' \notin \mathcal{P}_{U_{2,4}}^k$ by Lemma \ref{connectedness} and $\rho'$ cannot be an excluded minor.
        \item Assume $\rho'(Y) = 2k+1$.
        \begin{enumerate}
            \item Assume $\rho'(d) = k$. Then we have a contradiction by the same reasoning as in part (1)(b)(i) of the proof of Lemma \ref{nodecomps1}.
            \item Assume $\rho'(d) = k-1$. Then $\rho' = (\rho')|_Y \oplus (\rho')|_d$, contradiction by Lemma \ref{connectedness}.
        \end{enumerate}
        \item Assume $\rho'(Y) = 2k+2$. Then, regardless of whether $\rho'(d) = k$ or $k-1$, we have a contradiction by the same reasoning as in part (1)(b)(ii) of the proof of Lemma \ref{nodecomps1}.
    \end{enumerate}

    Hence, $\rho'(E)=3k-1$. Since $\rho'$ is assumed to be an excluded minor, $(\rho')^*$ must be as well by Lemma \ref{kduality}. By the definition of $k$-dual, we have $(\rho')^*(E) = 2k+1$. By Lemma \ref{compression commutes with dual}, it must be that $((\rho')^l_{\downarrow d})^* = ((\rho')^*)^{k-l}_{\downarrow d}$. Since $Ex_{(2,0,2)}^{2k}$ is self-$k$-dual, and $(\rho')^l_{\downarrow d} = \rho = Ex_{(2,0,2)}^{2k}$, it must be that $((\rho')^*)^{k-l}_{\downarrow d} = Ex_{(2,0,2)}^{2k}$ as well. Therefore, by the same reasoning for $\rho'(E)$, it must be that $(\rho')^*(E) = 3k-1$. But when $k \geq 3$, $3k-1 \neq 2k+1$, so this is a contradiction. Hence, $(\rho')^*$ cannot be an excluded minor. By Lemma \ref{kduality}, it must be that $\rho'$ also cannot be an excluded minor.
\end{proof}
This concludes the proof of Theorem \ref{no decompressions}.
\end{proof}

\appendix
\section{List of excluded minors for \texorpdfstring{$\mathcal{P}_{U_{2,4}}^k$}{}}

\noindent
$|E| = 1$, $k \geq 4$.
\begin{itemize}
    \item $Ex^m$ where $2 \leq m \leq k-2$, consisting of one rank-$m$ element. Its $k$-dual is $Ex^{k-m}$.
\end{itemize}

\noindent
$|E| = 2$, $k = 3$.
\begin{itemize}
    \item $Ex_\gamma^2$: a point lying on a line. Its $k$-dual is $Ex_\epsilon^4$.
    \item $Ex_\epsilon^4$: a line and a plane spanning rank $4$. Its $k$-dual is $Ex_\gamma^2$.
\end{itemize}

\noindent
$|E| = 2$, $k \geq 3$.
\begin{itemize}
    \item $Ex_\alpha^{k-1}$ consisting of two parallel rank-$(k-1)$ elements. Its $k$-dual is $Ex_\beta^{k+1}$.
    \item $Ex_\alpha^k$ consisting of two rank-$(k-1)$ elements spanning rank $k$. $Ex_\alpha^k$ is self-$k$-dual.
    \item $Ex_\beta^k$ consisting of two parallel rank-$k$ elements. $Ex_\beta^k$ is self-$k$-dual.
    \item $Ex_\beta^{k+1}$ consisting of two rank-$k$ elements spanning rank $k+1$. Its $k$-dual is $Ex_\alpha^{k-1}$.
    \item $Ex_\epsilon^k$ consisting of a rank-$(k-1)$ element lying on a rank-$k$ element. $Ex_\epsilon^k$ is identically self-$k$-dual.
\end{itemize}

\noindent
$|E| = 4$.
\begin{itemize}
    \item $Ex_{(4, 0, 0)}^2$ (also known as $U_{2,4}$) consisting of four collinear points. Its $k$-dual is $Ex_{(0, 0, 4)}^{4k-2}$.
    \item $Ex_{(3, 0, 1)}^{k+1}$ consisting of a rank-$k$ element and three collinear points spanning rank $k+1$ with no point lying on the rank-$k$ element. Its $k$-dual is $Ex_{(1, 0, 3)}^{3k-1}$. 
    
\begin{minipage}{\linewidth}
        $$\begin{tikzcd}[row sep = 10]&&&\underset{\color{blue}E\color{black}}{ k+1}\arrow[d,no head]\arrow[rrd,no head,dotted]&&\\&\underset{\color{blue}efg\color{black}}{ 2}\arrow[rru,no head]\arrow[d,no head]\arrow[rrd,no head,dotted]&\underset{\color{blue}efh\color{black}}{ k+1}\arrow[ru,no head]\arrow[d,no head]\arrow[rrd,no head,dotted]&\underset{\color{blue}egh\color{black}}{ k+1}\arrow[ld,no head]\arrow[rrd,no head,dotted]&&\underset{\color{blue}fgh\color{black}}{ k+1}\arrow[ld,no head]\arrow[d,no head]\\\underset{\color{blue}ef\color{black}}{ 2}\arrow[ru,no head]\arrow[rru,no head]&\underset{\color{blue}eg\color{black}}{ 2}\arrow[rru,no head]&\underset{\color{blue}eh\color{black}}{ k+1}\arrow[rrd,no head,dotted]&\underset{\color{blue}fg\color{black}}{ 2}\arrow[rru,no head]&\underset{\color{blue}fh\color{black}}{ k+1}\arrow[d,no head]&\underset{\color{blue}gh\color{black}}{ k+1}\\\underset{\color{blue}e\color{black}}{ 1}\arrow[u,no head]\arrow[rru,no head]\arrow[ru,no head]&&\underset{\color{blue}f\color{black}}{ 1}\arrow[ru,no head]\arrow[rru,no head]\arrow[llu,no head,dotted]&\underset{\color{blue}g\color{black}}{ 1}\arrow[u,no head]\arrow[rru,no head]\arrow[llu,no head,dotted]&\underset{\color{blue}h\color{black}}{ k}\arrow[ru,no head]&\\&&\underset{\color{blue}\emptyset\color{black}}{0}\arrow[u,no head]\arrow[ru,no head]\arrow[rru,no head]\arrow[llu,no head,dotted]&&&\end{tikzcd}$$
    \end{minipage}
    \item $Ex_{(2, 0, 2)}^{2k}$, consisting of two skew rank-$k$ elements and two points in rank $2k$ where each rank-$k$ element spans rank $k+1$ with the pair of points but spans neither point by itself. $Ex_{(2, 0, 2)}^{2k}$ is self-$k$-dual.
\begin{minipage}{\linewidth}
        $$\begin{tikzcd}[row sep = 10]&&&\underset{\color{blue}E\color{black}}{ 2k}\arrow[d,no head]\arrow[rrd,no head,dotted]&&\\&\underset{\color{blue}efg\color{black}}{ k+1}\arrow[rru,no head]\arrow[d,no head]\arrow[rrd,no head,dotted]&\underset{\color{blue}efh\color{black}}{ k+1}\arrow[ru,no head]\arrow[d,no head]\arrow[rrd,no head,dotted]&\underset{\color{blue}egh\color{black}}{ 2k}\arrow[ld,no head]\arrow[rrd,no head,dotted]&&\underset{\color{blue}fgh\color{black}}{ 2k}\arrow[ld,no head]\arrow[d,no head]\\\underset{\color{blue}ef\color{black}}{ 2}\arrow[ru,no head]\arrow[rru,no head]&\underset{\color{blue}eg\color{black}}{ k+1}\arrow[rru,no head]&\underset{\color{blue}eh\color{black}}{ k+1}\arrow[rrd,no head,dotted]&\underset{\color{blue}fg\color{black}}{ k+1}\arrow[rru,no head]&\underset{\color{blue}fh\color{black}}{ k+1}\arrow[d,no head]&\underset{\color{blue}gh\color{black}}{ 2k}\\\underset{\color{blue}e\color{black}}{ 1}\arrow[u,no head]\arrow[rru,no head]\arrow[ru,no head]&&\underset{\color{blue}f\color{black}}{ 1}\arrow[ru,no head]\arrow[rru,no head]\arrow[llu,no head,dotted]&\underset{\color{blue}g\color{black}}{ k}\arrow[u,no head]\arrow[rru,no head]\arrow[llu,no head,dotted]&\underset{\color{blue}h\color{black}}{ k}\arrow[ru,no head]&\\&&\underset{\color{blue}\emptyset\color{black}}{0}\arrow[u,no head]\arrow[ru,no head]\arrow[rru,no head]\arrow[llu,no head,dotted]&&&\end{tikzcd}$$
    \end{minipage}
    \item $Ex_{(1, 0, 3)}^{3k-1}$, whose $k$-dual is $Ex_{(3, 0, 1)}^{k+1}$.
    
\begin{minipage}{\linewidth}
        $$\begin{tikzcd}[row sep = 10]&&&\underset{\color{blue}E\color{black}}{ 3k-1}\arrow[d,no head]\arrow[rrd,no head,dotted]&&\\&\underset{\color{blue}efg\color{black}}{ 2k}\arrow[rru,no head]\arrow[d,no head]\arrow[rrd,no head,dotted]&\underset{\color{blue}efh\color{black}}{ 2k}\arrow[ru,no head]\arrow[d,no head]\arrow[rrd,no head,dotted]&\underset{\color{blue}egh\color{black}}{ 2k}\arrow[ld,no head]\arrow[rrd,no head,dotted]&&\underset{\color{blue}fgh\color{black}}{ 3k-1}\arrow[ld,no head]\arrow[d,no head]\\\underset{\color{blue}ef\color{black}}{ k+1}\arrow[ru,no head]\arrow[rru,no head]&\underset{\color{blue}eg\color{black}}{ k+1}\arrow[rru,no head]&\underset{\color{blue}eh\color{black}}{ k+1}\arrow[rrd,no head,dotted]&\underset{\color{blue}fg\color{black}}{ 2k}\arrow[rru,no head]&\underset{\color{blue}fh\color{black}}{ 2k}\arrow[d,no head]&\underset{\color{blue}gh\color{black}}{ 2k}\\\underset{\color{blue}e\color{black}}{ 1}\arrow[u,no head]\arrow[rru,no head]\arrow[ru,no head]&&\underset{\color{blue}f\color{black}}{ k}\arrow[ru,no head]\arrow[rru,no head]\arrow[llu,no head,dotted]&\underset{\color{blue}g\color{black}}{ k}\arrow[u,no head]\arrow[rru,no head]\arrow[llu,no head,dotted]&\underset{\color{blue}h\color{black}}{ k}\arrow[ru,no head]&\\&&\underset{\color{blue}\emptyset\color{black}}{0}\arrow[u,no head]\arrow[ru,no head]\arrow[rru,no head]\arrow[llu,no head,dotted]&&&\end{tikzcd}$$
    \end{minipage}
    \item $Ex_{(0, 0, 4)}^{4k-2}$, whose $k$-dual is $U_{2,4}$.
    
    \begin{minipage}{\linewidth}
        $$\begin{tikzcd}[row sep = 10]&&&\underset{\color{blue}E\color{black}}{ 4k-2}\arrow[d,no head]\arrow[rrd,no head,dotted]&&\\&\underset{\color{blue}efg\color{black}}{ 3k-1}\arrow[rru,no head]\arrow[d,no head]\arrow[rrd,no head,dotted]&\underset{\color{blue}efh\color{black}}{ 3k-1}\arrow[ru,no head]\arrow[d,no head]\arrow[rrd,no head,dotted]&\underset{\color{blue}egh\color{black}}{ 3k-1}\arrow[ld,no head]\arrow[rrd,no head,dotted]&&\underset{\color{blue}fgh\color{black}}{ 3k-1}\arrow[ld,no head]\arrow[d,no head]\\\underset{\color{blue}ef\color{black}}{ 2k}\arrow[ru,no head]\arrow[rru,no head]&\underset{\color{blue}eg\color{black}}{ 2k}\arrow[rru,no head]&\underset{\color{blue}eh\color{black}}{2k}\arrow[rrd,no head,dotted]&\underset{\color{blue}fg\color{black}}{ 2k}\arrow[rru,no head]&\underset{\color{blue}fh\color{black}}{ 2k}\arrow[d,no head]&\underset{\color{blue}gh\color{black}}{ 2k}\\\underset{\color{blue}e\color{black}}{ k}\arrow[u,no head]\arrow[rru,no head]\arrow[ru,no head]&&\underset{\color{blue}f\color{black}}{ k}\arrow[ru,no head]\arrow[rru,no head]\arrow[llu,no head,dotted]&\underset{\color{blue}g\color{black}}{ k}\arrow[u,no head]\arrow[rru,no head]\arrow[llu,no head,dotted]&\underset{\color{blue}h\color{black}}{ k}\arrow[ru,no head]&\\&&\underset{\color{blue}\emptyset\color{black}}{0}\arrow[u,no head]\arrow[ru,no head]\arrow[rru,no head]\arrow[llu,no head,dotted]&&&\end{tikzcd}$$
    \end{minipage}
\end{itemize}

\label{appendix: excluded minors}

\end{document}